\documentclass[11pt]{article}
\usepackage{amsmath}
\usepackage{amssymb}
\usepackage{amsthm}
\usepackage{amsfonts}
\usepackage{cite}
\usepackage{calc}
\usepackage{epsfig,afterpage}
\advance\textwidth by +1.3in
\advance\textheight by +1.3in
\advance\oddsidemargin by -0.5in
\advance\evensidemargin by -1.5in
\advance\topmargin by -0.6in
\catcode`\@=11

\@addtoreset{equation}{section}
\catcode`\@=12

\def\be{\begin{equation}}
\def\ee{\end{equation}}
\def\bal{\begin{aligned}}
\def\eal{\end{aligned}}
\def\bq{\begin{eqnarray}}
\def\eq{\end{eqnarray}}
\def\beq{\begin{eqnarray*}}
\def\eeq{\end{eqnarray*}}

\def\ba{\begin{array}}
\def\ea{\end{array}}
\def\bth{\begin{theorem}}
\def\eth{\end{theorem}}
\def\blm{\begin{lemma}}
\def\elm{\end{lemma}}
\def\bdf{\begin{definition}}
\def\edf{\end{definition}}
\def\bpr{\begin{proposition}}
\def\epr{\end{proposition}}
\def\brm{\begin{remark}}
\def\erm{\end{remark}}
\def\bnot{\begin{notation}}
\def\enot{\end{notation}}
\def\bobs{\begin{observation}}
\def\eobs{\end{observation}}
\def\bcrl{\begin{corrolary}}
\def\ecrl{\end{corrolary}}

\newcommand{\ab}{\mbox{\boldmath $a$}}
\newcommand{\bb}{\mbox{\boldmath $b$}}

\newcommand{\SSS}{{\bf S}}
\newcommand{\QSL}{{\bf Q\bf S\bf L}}

\newcommand{\QS}{{\bf Q\bf S}}
\newcommand{\IL}{{\bf I\bf L}}
\newcommand{\D}{{\bf D}}

\newcommand{\PP}{{\mathbb P}}
\newcommand{\R}{\mathbb{R}}
\newcommand{\C}{\mathbb{C}}

\newcommand{\Z}{\bf{Z}}
\newcommand{\sign}{\mbox{\rm sign\,}}
\newcommand{\Res}{\mbox{\rm Res\,}}
\newcommand{\Discriminant}{\mbox{\rm Discriminant\,}}
\newcommand{\Jacob}{\mbox{\rm Jacob\,}}
\newcommand{\Hess}{\mbox{\rm Hess\,}}
\newcommand{\Supp}{\mbox{\rm Supp\,}}

\newtheorem{theorem}{Theorem}[section]
\newtheorem{lemma}{Lemma}[section]
\newtheorem{definition}{Definition}[section]
\newtheorem{proposition}[lemma]{Proposition}
\newtheorem{remark}{Remark}[section]
\newtheorem{notation}[remark]{Notation}
\newtheorem{observation}[remark]{Observation}
\newtheorem{corrolary}[remark]{Corrolary}
\newcommand{\BProof}{\noindent{\it Proof:}\ \,}
\newcommand{\EProof}{\ \hfill\rule[-1mm]{1.5mm}{2.0mm}}

\begin{document}

\setcounter{secnumdepth}{5}

\title{ Planar quadratic vector fields with invariant lines  of  total multiplicity
at least five  }

\author{Dana SCHLOMIUK\thanks{Work supported by NSERC and
        by the Quebec Education Ministry}\\
        {\it D\'epartement de Math\'ematiques et de Statistiques}\\
        {\it Universit\'e de Montr\'eal}\\ {\it C.P. 6128, succ. Centre-ville
        Montr\'eal, QC H3C 3J7 Canada}\\
        E-mail: dasch@dms.umontreal.ca
 \and
        Nicolae VULPE\thanks{Partially supported by  NSERC}\\
        {\it Institute of Mathematics and Computer Science}\\
        {\it Academy of Science of Moldova }\\
        {\it str. Academiei 5, Chi\c sin\u au MD-2028, Moldova}\\
         E-mail: nvulpe@math.md}
\date{}
\maketitle
\begin{abstract}

In this article we  consider the action of affine group and time
rescaling on planar quadratic differential systems. We construct a
system of representatives of the orbits of systems with at least
five invariant lines, including the line at infinity and including
multiplicities.  For each orbit we  exhibit  its configuration. We
characterize in terms of algebraic invariants and comitants and
also geometrically, using divisors of the complex projective
plane, the class of quadratic differential systems with at least
five invariant lines. These conditions are such that no matter how
a system may be presented, one can  verify by using them whether
the system has or does not have  at least five invariant lines and
 to check to which orbit (or  family of orbits) it
belongs.

\end{abstract}

%\tableofcontents

%%%%%%%%%%%%%%%%%%%%%%%%%
\section{Introduction } We consider here real planar  differential systems of the form
\be\label{il1}
 (S)\qquad \frac {dx}{dt}= p(x,y),\qquad
  \frac {dy}{dt}= q(x,y),\hphantom{--  }
\ee
where $p,\,q\in \R[x,y]$, i.e. $p,\ q$ are polynomials in $x,\ y$
over $\R$, and their associated vector fields
 \be
     \tilde D=p(x,y)\frac{\partial}{\partial x} + q(x,y)\frac{\partial}{\partial y}.
     \label{il2}
\ee
Each such system generates a complex differential vector field
when the variables range over $\C$. To the complex systems we can
apply the work of Darboux on integrability via invariant algebraic
curves (cf.\cite{Darb}). For a brief introduction to the work of
Darboux we refer to the survey article \cite {Dana1}.  Some
applications of the work of Darboux in connection with the problem
of the center are given in \cite{Dana2}.

For the system (\ref{il1}) we can use the following definition.
 \bdf \label{df1}
 An affine algebraic invariant  curve of a polynomial system {\rm(\ref{il1})}
(or an algebraic particular integral) is a curve $f(x,y)=0$ where
$f\in \C[x,y]$, $\deg(f)\ge1$, such that there exists $k(x,y)\in
\C[x,y]$ satisfying $\tilde Df=fk$ in $\C[x,y]$. We call $k$ the
cofactor of $f$ with respect to the system.
\edf

Poincar\'e was the first to appreciate the work of Darboux
\cite{Darb}, which it called "admirable" (see \cite{Po2}) and
inspired by Darboux's work, Poincar\'e wrote two articles
\cite{Po3},\cite{Po4} where he also stated a problem still open
today.

With this brilliant work Darboux open up a whole new area of
investigations where one studies how the presence of particular
algebraic integrals impacts on global properties of the systems,
for example on global integrability. In recent years there has
been a surge in activity in this area of research and this article
is part of a growing literature in the subject. In particular we
mention here \cite{Cris_Llib}, \cite{Cris_Llib_Pant} and  the
 work of C. Christopher, J.V. Perreira and J. Llibre on the
notion of multiplicity of an invariant algebraic curve of a
differential system \cite{Cris_Llib_Per}.

 In this article, which is based on \cite{Dana_Vlp1}, we study
systematically  the simplest kind of such a structure, i.e.
quadratic systems (\ref{il1}) possessing invariant lines. Some
references on this topic are:
\cite{Sib3,Druzhkova,Art_Llib2,Lyubim1,
Lyubim2,Popa_Sib1,Popa2,Sokulski,ZX,Lib_Vul}.

To a line $f(x,y)=ux+vy+w=0$ we associate its projective
completion $F(X,Y,Z)=uX+vY+wZ=0$ under the embedding
$\C^2\hookrightarrow \PP_2(\C)$, $(x,y)\mapsto [x:y:1]$. The line
$Z=0$ is   called the line at infinity of the system (\ref{il1}).
It follows from the work of Darboux that each system of
differential equations of the form (\ref{il1}) yields a
differential equation on the complex projective plane which is the
compactification of the complex system (\ref{il1}) on $\PP_2(\C)$
(cf. Section \ref{proj_eq}). The line $Z=0$ is an
 invariant manifold of this complex differential equation.
\bnot
  Let us denote by
\beq
  \QS &=& \left\{\ S\ \left|\ba{l}\ S\ \hbox {is a system {\rm (\ref{il1})} such that}\
       \gcd(p(x,y),q(x,y))=1\\ \ \ \hbox{and}\quad \max\big(\deg(p(x,y)),\deg(p(x,y))\big)=2
            \ea \right.\right\}; \\
  \QSL &=& \left\{\,S\in \QS \left|\ba{l} \ \hbox{ $S$ possesses at
                                       least one invariant affine line or}\\
                                     \ \, \hbox{the line at infinity with multiplicity
                                       at least two}
                       \ea \right.\right\}.
\eeq
 \enot
For the multiplicity of the line at infinity the reader is
refereed to Section \ref{proj_eq}.

We shall call \textit{degenerate quadratic differential system} a
system (\ref{il1}) with $\deg \gcd(p,q)\ge1$ and
$\max\big(\deg(p),\deg(q)\big)=2$.
 \bnot
 To a  system {\rm(\ref{il1})} in $\QS$ we can associate  a point in
 $ \R^{12}$, the
ordered tuple of the coefficients of $p(x,y)$, $q(x,y)$  and this
correspondence is an injection
\be\label{bigection}
\bal
{\cal B}:\quad \QS &\hookrightarrow \R^{12}\\
S\ \ & \mapsto\ \ \ab = {\cal B}(S)
 \eal
 \ee
The topology of $\R^{12}$ yields an induced topology on {\bf QS}.
 \enot

\bdf\label{def:multipl} We say that an invariant straight line ${\cal
L}(x,y)=ux+vy+w=0$, $(u,v)\ne(0,0)$, $(u,v,w)\in \C^3$ for a
quadratic vector field $\tilde D$ has multiplicity $m$ if there
exists a sequence of real quadratic vector fields $\tilde D_k$
converging to $\tilde D$, such that each $\tilde D_k$ has $m$
distinct (complex) invariant straight lines ${\cal
L}^1_k=0,\ldots, {\cal L}^m_k=0$, converging to ${\cal L}=0$ as
$k\to\infty$ (with the topology of their coefficients), and this
does not occur for $m+1$.
 \edf

\bpr\label{pr:m_il} {\rm\!\!\cite{Art_Llib2}} The
maximum number of invariant lines (including the line at infinity
and including multiplicities) which a quadratic system could have
is six.\epr

\bdf
We call configuration of invariant lines of a system $(S)$ in
$\QSL$ the set of all its invariant lines (real or/and complex),
each endowed with its own multiplicity and  together with all the
real singular points of $(S)$ located on these   lines, each one
endowed with its own multiplicity.
\edf

We associate to each system in {\bf QSL} its   configuration of
invariant lines. In analogous manner to how we view the phase
portraits of the systems on the Poincar\'e disc (see for example,
 \cite{Lib_DS}), we can also view the configurations of real
lines on the disc. To help imagining the full configurations, we
complete the picture by drawing dashed lines whenever these are
complex.

On the class of quadratic systems acts the group of real affine
transformations and time rescaling. Since quadratic systems depend
on 12 parameters and since this group depends on 7 parameters, the
class of quadratic systems modulo this group action, actually
depends on five parameters.

It is clear that the configuration of invariant lines of a system
is an affine invariant. The notion of multiplicity  defined by
Definition \ref{def:multipl} is invariant under the group action,
i.e. if a quadratic system $S$ has an invariant line $l$ of
multiplicity $m$, then each system $\tilde S$ in the orbit of $S$
under the group action has an invariant line $l~$ of the same
multiplicity $m$.

In this article  we shall consider the case when the system (\ref{il1}) has at
least five invariant lines considered with their multiplicities.

The problems which we solve in this article are the following:

I) Construct a system of representatives of the orbits of systems with at least
five invariant lines, including the line at infinity and including
multiplicities.
 For each orbit  exhibit  its configuration.

II) Characterize in terms of algebraic invariants and comitants
and also geometrically, using divisors or zero-cycles of the
complex projective plane, the class of quadratic differential
systems with at least five invariant lines. These conditions
should be such that no matter how a system may be presented to us,
we should be able to verify by using them whether the system has
or does not have  at least five invariant lines and to check to
which orbit or perhaps family of orbits it belongs.

Our main results are formulated in Theorems 5.1 and 6.1. Theorem
5.1 gives a total of 11 distinct orbits of systems with a
configuration with exactly six invariant lines including the line
at infinity and including multiplicities. Theorem 6.1 gives a
system of representatives  for 17 distinct orbits  of systems with
exactly five invariant lines including the line at infinity and
including multiplicities. Furthermore we give a complete list of
representatives of the remaining  orbits which are classified in
13 one-parameter families. We characterize each one of these 13
families  in terms of algebraic invariants and comitants and
geometrically.  As the calculation of invariants and comitants can
be implemented on a computer, this verification can be done by a
computer.

All quadratic  systems with at least five invariant lines
including the line at infinity and including multiplicities are
algebraically integrable, i.e. they all have the rational first
integrals and  the phase portraits of these systems can easily be
drawn. We leave the discussion of issues related to integrability,
as well as the drawing of the phase portraits of the systems we
consider here, in a follow up paper of this work.

The invariants and comitants of differential equations used in the
classification theorems (Theorems \ref{th_mil_6} and \ref{th_mil_5}) are
obtained following the theory established by K.Sibirsky and his disciples (cf.
\cite{Sib1}, \cite{Sib2}, \cite{Vlp1}, \cite{Popa4}).

\section{ Differential equations in $\PP_2(\C)$ of first degree and
  first order and their invariant projective curves}
\label{proj_eq}
In \cite{Darb} Darboux considered differential  equations  of first degree and
  first order of the complex projective plane. These are equations of the form
\be
\left|\begin{array}{ccc}
L & M & N\\[0mm]
X & Y & Z\\
dX & dY & dZ\\
\end{array}\right|=0 \tag{$CF$}
\ee
where $L$, $M$, $N$ are homogeneous polynomials of the same degree
$m$. These are called equations in Clebsch form $(CF)$. \footnote{
Darboux
  used the notion of Clebsch connex to define them.}

We remark that we can have an infinity of such equations yielding
the same integral curves. Indeed, for any ordered triple $L,M,N$
of homogeneous polynomials of the same degree $m$  and for any
homogeneous polynomial $A$ of degree $m-1$, the $(CF)$-equation
corresponding to
\be\label{equiv}
L'=L+AX,\quad M'=M+AY,\quad N'=N+AZ
\ee
has the same integral curves as the equation $(CF)$. Two equations
$(CF)$ determined by polynomials $L,M,N$ and $L',M',N'$ satisfying
\eqref{equiv} are said to be equivalent.
\bth[Darboux, \cite{Darb}]\label{th:Darb1}
Let  $L,\,M,\,N$ be homogeneous polynomials of the same degree $m$
over $\C$. Then there exists a unique $A$, more precisely
$$
A=-\frac{1}{m+2}\Big(\frac{\partial L }{\partial X}+\frac{\partial
M }{\partial Y}+\frac{\partial N }{\partial Z}\Big)
$$
such that  if $L',M',N'$ satisfy (\ref{equiv}) for this $A$ then
$$
\frac{\partial L' }{\partial X}+\frac{\partial M' }{\partial Y}+\frac{\partial
N' }{\partial Z}\equiv0.
$$
\eth
\bth[Darboux, \cite{Darb}]\label{th:Darb2} Every equation (CF) with
$m=\deg(L)=\\deg(M)=\deg(N)$ is equivalent to an equation
\be\label{ABC}
A\,dX +B\,dY+C\,dZ=0
\ee
where $A,\ B,\ C$ are homogeneous polynomials of degree $m+1$ subject
to the identity
\be\label{ABC=0}
A\,X +B\,Y+C\,Z=0
\ee
\eth
\bdf[Darboux, \cite{Darb}]\label{def:Darb} An  algebraic invariant
curve for an equation $(CF)$ is a projective curve $F(X,Y,Z)=0$
where $F$ is a homogeneous polynomial over $\C$ such that $F\mid
\hat DF$ where $\hat D$ is the differential operator
$$
\hat D=L\frac{\partial }{\partial X}+M\frac{\partial }{\partial
Y}+N\frac{\partial }{\partial Z}
$$
i.e. $\exists K\in \C[X,Y,Z]$ such that $\hat DF=FK$. $K$ is
called the {\it cofactor} of $F$ with respect to the equation
$(CF)$.
\edf

We now show that this definition is in agreement with Definition
\ref{df1}, i.e. it includes as a particular Definition \ref{df1}.

 To a system (\ref{il1}) we can associate an equation
(\ref{ABC}) subject to the identity (\ref{ABC=0}). We first
associate to the systems (\ref{il1}) the differential  form
$$
\omega_1=q(x,y)dx-p(x,y)dy
$$
and its associated  differential equation $\omega_1=0$.

We consider the map $j: \C^3 \setminus \{Z = 0\} \to \C^2$, given
by $i(X,Y,Z)= (X/Z,Y/Z)=(x,y)$ and suppose that $\max\big(
\deg(p),\deg(q)\big)= m>0$. Since $x=X/Z$ and $y=Y/Z$ we have:
\[
dx= (ZdX-XdZ)/Z^2 \ , \qquad dy= (ZdY-YdZ)/Z^2 \ ,
\]
the pull--back form $j^*(\omega_1)$ has poles at $Z=0$ and its
associated equation $j^*(\omega_1)=0$  can be written as
\[
j^*(\omega_1)= q(X/Z,Y/Z) (ZdX-XdZ)/Z^2 - p(X/Z,Y/Z) (ZdY-YdZ)/Z^2
= 0 .
\]
Then the $1$--form $\omega= Z^{m+2} j^*(\omega_1)$ in
$\C^3\setminus \{Z\ne 0\}$ has homogeneous polynomial coefficients
of degree $m+1$, and for $Z\ne0$ the equations $\omega=0$ and
$j^*(\omega_1)=0$ have the same solutions. Therefore the
differential equation $\omega=0$ can be written as (\ref{ABC})
where
\bq
A(X,Y,Z) &=& Z Q(X,Y,Z)= Z^{m+1} q(X/Z,Y/Z)\,, \nonumber \\
B(X,Y,Z) &=& -Z P(X,Y,Z)=-Z^{m+1} p(X/Z,Y/Z)\,, \label{24} \\
C(X,Y,Z) &=& Y P(X,Y,Z)- X Q(X,Y,Z)\, \nonumber
\eq
and $P(X,Y,Z)=Z^m p(X/Z,Y/Z),\quad Q(X,Y,Z)=Z^m q(X/Z,Y/Z)$.
 Clearly $A$,
$B$ and $C$ are homogeneous polynomials of degree $m+1$ satisfying
(\ref{ABC=0}).

The equation (\ref{ABC}) becomes in this case
$$
P(YdZ-ZdY) +Q(ZdX-XdZ)=0
$$
or equivalently
\begin{equation}\label{CFp}
\left|\begin{array}{ccc}
P& Q & 0\\[0mm]
X & Y & Z\\
dX & dY & dZ\\
\end{array}\right|=0.
\end{equation}
We observe that $Z=0$ is an algebraic invariant curve of this
equation according to Definition \ref{def:Darb}, with cofactor
\mbox{$K=0$}. We shall also say that $Z=0$ is an invariant line
for the systems~(\ref{il1}).

To an affine algebraic curve $f(x,y)=0$, $\deg f=n$, we can
associate its projective completion $F(X,Y,Z)=0$ where
$F(Z,Y,Z)=Z^{n}f(X/Z,Y/Z)$. From the indicated above the
correspondence between a system (\ref{il1}) and equation
(\ref{CFp}) follows the next proposition.
\bpr\label{pr:f=0}
Let $f=0$ ($\deg f=n$) be an invariant  algebraic curve of
(\ref{il1}) according to Definition \ref{df1}, with cofactor
$k(x,y)$. Then its associated projective completion $F(X,Y,Z)=0$
where $F(Z,Y,Z)=Z^{n}f(X/Z,Y/Z)$ is an invariant algebraic curve
according to Definition \ref{def:Darb} for the equation
(\ref{CFp}), with cofactor $K(X,Y,Z)=Z^{m-1}k(X/Z,Y/Z)$.
\epr
Conversely, starting now with an  equation in Clebsch form $(CF)$
we can consider its restriction on the affine chart $Z=1$ and
associate a differential system:
\be\label{syst:Z=1}
\left|\begin{array}{ccc}
L & M & N\\[0mm]
X & Y & Z\\
dX & dY & dZ\\
\end{array}\right|=0\
 \longrightarrow \
 (\hat M -y\hat N)dx-(\hat L -x\hat N)dy=0\
 \longrightarrow\ \left\{\ba{l} \dot x=\hat L -x\hat N\\ \dot y= \hat M
 -y\hat  N,\ea\right.
\ee
where $\hat L=L(x,y,1)$, $\hat M=M(x,y,1)$, $\hat N=N(x,y,1)$. The
following proposition follows easily by using Euler's formula
$XF'_X+YF'_Y+ZF'_Z= nF$ for a homogeneous polynomial $F(Z,Y,Z)$ of
degree n.
\bpr
Let $F(Z,Y,Z)=0$ ($\deg F=n$)  be an invariant algebraic curve
(according to Definition \ref{def:Darb}) for the equation $(CF)$
with cofactor $K(X,Y,Z)$, such that $Z\nmid F$. Then
$f(x,y)=F(x,y,1)=0$  is an invariant affine algebraic curve
(according to Definition \ref{df1}) of the differential system in
\eqref{syst:Z=1} corresponding to $(CF)$,  with cofactor $k(x,y)=K(x,y,1)-nN(x,y,1)$.
\epr
\bdf\label{def:Z-mult} We say that $Z=0$  is an invariant line of multiplicity
$m$  for a system $(S)$ of the form (\ref{il1}) if and only if
there exists a sequence of systems $(S_i)$ of the form (\ref{il1})
such that $(S_i)$ tend to $(S)$ when $i\to \infty$ and the systems
$(S_i)$ have $m-1$ distinct invariant affine lines ${\cal
L}^j_i=u^j_ix+v^j_iy+w^j_i=0$, $(u^j_i,v^j_i)\ne(0,0)$,
$(u^j_i,v^j_i,w^j_i )\in \C^3$ $(j=1,\ldots,m-1)$ such that\quad
for every $ j$\
$\displaystyle{\lim_{i\to\infty}(u^j_i,v^j_i,w^j_i)=(0,0,1)}$.
\edf

\section{ Divisors associated to invariant lines configurations }
\label{divisors} Consider real differential systems of the form:
\be\label{2l1}
(S)\qquad \left\{ \ba{ll}
\displaystyle  \frac {dx}{dt}&=p_0+ p_1(x,y)+\,p_2(x,y)\equiv p(x,y), \\[2mm]
\displaystyle  \frac {dx}{dt}&=q_0+ q_1(x,y)+\,q_2(x,y)\equiv
q(x,y)
\ea\right.
\ee
with
$$
\bal
&p_0=a_{00},\quad p_1(x,y)=  a_{10}x+ a_{01}y,\quad  p_2(x,y)=
a_{20}x^2 +2
a_{11}xy + a_{02}y^2,\\
&q_0=b_{00},\quad q_1(x,y)=  b_{10}x+ b_{01}y,\quad\  q_2(x,y)=
b_{20}x^2 +2
b_{11}xy + b_{02}y^2.\\
\eal
$$
 Let $
a=(a_{00},a_{10},a_{01},a_{20},a_{11},a_{02},b_{00},b_{10},b_{01},b_{20},b_{11},b_{02})$
be the 12-tuple of the coefficients of system \eqref{2l1} and
denote $\R[
a,x,y]=\R[a_{00},a_{10},a_{01},a_{20},a_{11},a_{02},b_{00},b_{10},b_{01},b_{20},b_{11},b_{02},x,y]$.
\bnot Let us denote by $\ab=(\ab_{00},\ab_{10}\ldots,\bb_{02})$ a
point in $\R^{12}$. Each  particular system \eqref{2l1} yields an
ordered 12-tuple $\ab$ of its coefficients.
\enot
\bnot Let
$$
\bal
  P(X,Y,Z)=& p_0(\ab)Z^2+ p_1(\ab,X,Y)Z+\,p_2(\ab,X,Y)=0,\\
  Q(X,Y,Z)=& q_0(\ab)Z^2+ q_1(\ab,X,Y)Z+\,q_2(\ab,X,Y)=0.\\
\eal
$$
We denote $\quad \sigma(P,Q)= \{w\in\PP_2(\C)\ |\ P(w)= Q(w)=0\}$.
\enot

\bdf\label{df3_1a} We consider formal expressions of the form $\D = \sum
n(w)w$\ where $w \in\PP_2(\C)$ or $w$ is an irreducible curve of
$\PP_2(\C)$ and $n(w)$ is an integer and only a finite number of
$n(w)$ are not zero. Such an expression  will be called a
zero-cycle of $\PP_2(\C)$ (respectively a divisor of $Z=0$ or a
divisor of $\PP_2(\C)$ ) if $w \in\PP_2(\C)$  (respectively, $w$
belongs to the line $Z=0$, or $w$ is an irreducible curve of
$\PP_2(\C)$). We call degree of the expression $\D$ the integer
$\deg(\D) = \sum n(w)$. We call support of   $\D$ the set
$\Supp(\D)$ of points $w$ such that $n(w)\ne0$.
\edf
In this section, for systems \eqref{2l1}  we shall assume the
conditions $\max(\deg(p),\deg(q))=2$ and $\gcd(p,q)=1$.
\bdf Let $C(X,Y,Z)=YP(X,Y,Z) - XQ(X,Y,Z)$.
\beq
    \D_{{}_S}(P,Q) &=& \sum_{w\in\sigma(P,Q)}I_w(P,Q)w;\\
    \D_{{}_S}(C,Z) &=& \sum_{w\in\{ Z = 0\}}I_w(C,Z)w \quad \mbox{if}\quad Z\nmid C(X,Y,Z);\\
    \D_{{}_S}(P,Q;Z) &=& \sum_{w\in\{ Z = 0\}}I_w(P,Q)w;\\
    \widehat\D_{{}_S}(P,Q,Z) &=& \sum_{w\in\{ Z = 0\}}\Big(I_w(C,Z),\,I_w(P,Q)\Big)w,\\
\eeq
where $I_w(F,G)$ is the intersection number (see, \cite{Fult}) of
the  curves defined by homogeneous polynomials  $F,\ G\in
\C[X,Y,Z]$\ and $\deg(F),\deg(G)\ge1$.

\edf
\bnot
\be\label{invar1}
\bal
  n_{{}_\R}^{\!{}^\infty} =&\# \{ w\in Supp\, \D_{{}_S}(C,Z)\,\big|\, w\in \PP_2(\R)\};\\
  d_{\sigma}^{{}^\infty} =&\deg \D_{{}_S}(P,Q;Z).
\eal
\ee
\enot

A complex projective line $uX+vY+wZ=0$ is invariant for the system
$(S)$ if either it coincides with $Z=0$ or it is the projective
completion of an invariant affine line $ux+vy+w=0$.

\bnot
 Let  $S\in \QSL$. Let us  denote
$$
\bal
    \IL(S)=&\left\{\ \ l\ \ \left|\ba{ll} & l\  \hbox {is a line in $\PP_2(\C)$ such }\\
                                           & \hbox{that}\  l\ \hbox{is invariant for}\ (S)\\
                      \ea\ \right\}\right.;\\
     M(l)=& \ \hbox{the multiplicity of the invariant line $l$ of
     $(S)$}.
\eal
$$
\enot
\brm
We note that the line  $l_\infty: Z=0$ is included  in $\IL(S)$
for any $S\in \QSL$.
\erm
 Let $l_i\,:\ f_i(x,y)=0$, $i=1,\ldots,k$,  be all the distinct invariant
affine lines (real or complex) of a system $S\in \QSL$. Let
$l'_i\,:\ {\cal F}_i(X,Y,Z)=0 $ be the complex projective
completion of $l_i$.
\bnot We denote
\beq
    && {\cal G}\,:\quad \ \prod_{i}{\cal F}_i(X,Y,Z)\,Z=0;\quad
        Sing\, {\cal G}=\left\{w\in {\cal G}|\ w\ \mbox{is a singular point of}\ {\cal G}\right\};  \\
    &&\nu(w)= \ \hbox{the multiplicity of the point $w$, as a point of}\
               {\cal G}.
\eeq
\enot
\bdf
\beq
    \D_{{}_\IL}(S)&=&\sum_{l\in\IL(S)} M(l)l,\quad (S)\in\QSL;\\
    Supp\, \D_{{}_\IL}(S)& =& \{\, l\ |\ l\in \IL(S)\}.
\eeq
\edf
\bnot
\be\label{invar}
\bal
  M_{{}_\IL}=&\deg\D_{{}_\IL}(S);\\
  N_{{}_\C} =&\# Supp\, \D_{{}_\IL};\\
  N_{{}_\R} =&\# \{   l\in Supp\, \D_{{}_\IL}\,\big|\,l\in  \PP_2(\R) \};\\
 n^{{}^\R}_{{}_{{\cal G},\,\sigma}}=&\# \{\omega\in Supp\, \D_{{}_S}(P,Q)\, |\,
                    \omega\in  {\cal G}\raisebox{-0.3em}[0pt][0pt]{$\big|_{\R^2}$}\};\\
 d^{\,{}^\R}_{{}_{{\cal G},\,\sigma}}=&\sum_{\omega\in{\cal G}\raisebox{-0.2em}[0pt][0pt]{$|_{\R^2}$}}I_\omega(P,Q);\\
  m_{{}_{\cal G}}=& \max\{\nu(\omega)\, |\, \omega\in Sing\, {\cal G}\raisebox{-0.2em}[0pt][0pt]{$|_{\C^2}$}\};\\
  m^{{}^\R}_{{}_{\cal G}}=& \max\{\nu(\omega)\, |\, \omega\in Sing\, {\cal
  G}\raisebox{-0.2em}[0pt][0pt]{$|_{\R^2}$}\}.
\eal
\ee
\enot
%
%%%%%%%%%   New  Section %%%%%%%%%%%
%%%%%%%%%%%%%%%%%%%%%%%%%%%%%%%%%%
\section{The main $T$-comitants associated to configurations of
invariant lines}\label{T-comit} On the set $\widehat\QS$
  of all  differential systems of the
form \eqref{2l1} acts the group $Aff(2,\R)$ of affine
transformation on the plane. Indeed for every $g\in Aff(2,\R)$,
$g:\ \R^{2}\longrightarrow \R^{2}$ we have:
$$
  g:\ \left(\ba{c} \tilde x\\ \tilde y \ea\right) =
      M\left(\ba{c} x\\ y \ea\right) +B;\qquad
  g^{-1}:\ \left(\ba{c} x\\ y \ea\right) =
      M^{-1}\left(\ba{c}\tilde x\\\tilde y \ea\right) -M^{-1}B,
$$
where $M=|| M_{ij} || $ is a $2\times 2$ nonsingular matrix, $B$
is a $2\times 1$ matrix over $\R$. For every $S\in \widehat\QS $
we can form its transformed system $\tilde S=g S$:
$$
  \frac{\partial\tilde x}{\partial t} =\tilde p(\tilde x,\tilde y),\qquad\quad
  \frac{\partial\tilde y}{\partial t} =\tilde q(\tilde x,\tilde y),\eqno(\tilde S)
$$
where
$$
  \left(\ba{c} \tilde p(\tilde x,\tilde y)\\ \tilde q(\tilde x,\tilde y)\ea\right) =
      M\left(\ba{c} (p\,{\mbox{\footnotesize$\circ$}}\, {g^{-1}})(\tilde x,\tilde y)\\
                   (q\,{\mbox{\footnotesize$\circ$}}\, {g^{-1}})(\tilde x,\tilde y) \ea\right).
$$
The map
\beq
&&\\[-8mm]
   Aff(2,\R)\times \widehat\QS\ &\longrightarrow &\ \widehat\QS \\
   (g,\ \, S)\ \ \, &\longrightarrow &\ \tilde S=g\,{\mbox{\footnotesize$\circ$}}\,S\\
&&\\[-12mm]
\eeq
verifies the axioms for a left group action.  For every  subgroup
$G\subseteq Aff(2,\R)$ we have an  induced action of $G$ on
$\widehat\QS$. We can identify the set $\widehat\QS$ of systems
\eqref{2l1} with $\R^{12}$ via the map  $\widehat\QS$ $\longrightarrow \R^{12}$  which
associates to
 each system \eqref{2l1} the 12-tuple $\ab=(\ab_{00},\ab_{10}\ldots,\bb_{02})$ of
 its coefficients.

  The action of  $Aff(2,\R)$ on $\widehat\QS$ yields an action
of this group on $\R^{12}$. For every $ g\in Aff(2,\R)$ let $r_g:\
\R^{12}\longrightarrow \R^{12}$, $r_g(\ab)=\tilde \ab$ where
$\tilde \ab$ is the 12-tuple of coefficients of $\tilde S$. We
know \mbox{(cf. \cite {Sib1})} that $r_g$ is linear and that the
map $r:\ Aff(2,\R)\longrightarrow GL(12,\R)$ thus obtained is a
group homomorphism. For every  subgroup $G$ of $Aff(2,\R)$, $r$
induces a representation of $G$ onto a subgroup $\cal G$ of
$GL(12,\R)$.
 \bdf \label{def:comit} A polynomial $U(
a\,,x,y)\in \R[a,x,y]$ is called a comitant  of systems
\eqref{2l1} with respect to a subgroup  $G$ of $Aff(2,\R)$, if
there exists $\chi\in \Z$ such that for every $({\frak g},\, \ab
)\in G\times\R^{12}\ $ and for every $(x,y)\in \R^2$ the following
relation holds:
$$
U(r_{\frak g}(\ab) ,\ {\frak g}(x,y)\,)\equiv\ (\det\,{\frak
g})^{-\chi}\, U( \ab,x,y),
$$
where $\det {\frak g}=\det M$. If the polynomial $U$ does not
explicitly depend  on $x$ and $y$ then it is called invariant. The
number $ \chi\in \Z $  is called the weight of the
 comitant $ U( a ,x,y)$.
If $G=GL(2,\R)$ (or $G=Aff(2,\R)$\,) then the comitant $U( a
,x,y)$ of systems \eqref{2l1} is called $GL$-comitant
(respectively, affine  comitant). \edf
%%%%%%%%%%%%%%
\bdf\label{def_G}
A subset $X\subset \R^{12}$ will be called $G$-invariant,\ if \
for every $ {\frak g}\in G$\ we have\ \mbox{$r_{\frak
g}(X)\subseteq X$}.
\edf
 Let us consider the polynomials
\beq
   C_i(a,x,y)&=&yp_i(a,x,y)-xq_i(a,x,y)\in \R[a,x,y],\ i=0,1,2, \\
  D_i(a,x,y)&=&\frac{\partial}{\partial x}p_i(a,x,y)+
        \frac{\partial}{\partial y}q_i(a,x,y)\in \R[a,x,y],\ i=1,2.
\eeq
As it was shown in \cite{Sib1} the polynomials
\be\label{C_i:D_i}
\Big\{\ C_0(a,x,y),\quad C_1(a,x,y),\quad C_2(a,x,y),\quad D_1(a), \quad
D_2(a,x,y)\ \Big\}
\ee
of degree one in the coefficients of systems \eqref{2l1} are
$GL$-comitants   of these systems.
\bnot Let $f,$
$g\in$ $\R[a,x,y]$ and
\be \label{trsv}
  (f,g)^{(k)}=
   \sum_{h=0}^k (-1)^h {k\choose h}
   \frac{\partial^k f}{ \partial x^{k-h}\partial y^h}\
   \frac{\partial^k g}{ \partial x^h\partial y^{k-h}}.
\ee
$(f,g)^{(k)}\in \R[a,x,y] $ is called the transvectant of index $k$ of $(f,g)$
(cf. {\rm\cite{Gr_Yng}, \cite{Olver}})
\enot

\bth\label{th:Vlp} {\rm\cite{Vlp1}} Any $GL$-comitant  of systems \eqref{2l1} can
be constructed from the elements of the set \eqref{C_i:D_i} by
using the operations: $+,\, -,\,\times,$   and by applying the
differential operation $(f,g)^{(k)}$.
\eth

 Let $T(2,\R)$ be the subgroup of $Aff(2,\R)$ formed by
translations. Consider the linear representation of $T(2,\R)$ into
its corresponding subgroup ${\cal T}\subset GL(12,\R)$, i.e. for
every $ \tau\in T(2,\R)$,\ $\tau:\ x=\tilde x+\alpha, y=\tilde
y+\beta$ we consider as above $r_\tau:\ \R^{12}\longrightarrow
\R^{12}$.
%%%%%%%%%%
\bdf\label{def:T-com}
Consider a polynomial\quad $U(a,x,y)=\sum_{j=0}^{d} U_{i}( a
)x^{d-j}y^j\in \R[a,x,y]$ which is a  $GL$-comitant of  systems
\eqref{2l1}. We say that this polynomial is a  $T$-comitant of  systems
\eqref{2l1}
if for every $(\tau,\, \ab )\in T(2,\R)\times \R^{12}$ \
$U_j(r_\tau(\ab))\, =\, U_j( \ab )$, $\forall$ $j=0,1,\ldots,d$.
\edf
Consider $s$ polynomials\quad
$U_i(a,x,y)=\sum_{j=0}^{d_i} U_{ij}( a )x^{d_i-j}y^j\in
\R[a,x,y]$, $i=1,\ldots,s$\quad and assume that the polynomials
$U_i$ are  $GL$-comitants of  systems
\eqref{2l1} where $d_i$ denotes the degree of the binary form $U_i( a,x,y)$ in
$x$ and $y$ with coefficients in $\R[a]$. We denote  by
$$
{\cal U}=\left\{\,U_{ij}( a )\in \R[a]\ |\ i=1,\ldots,s,\
j=0,1,\ldots,d_i\,\right\},
$$
 the set of the coefficients in $\R[a]$ of the $GL$-comitants $U_i( a ,x,y)$,
$i=1,\ldots,s$  and by $V(\cal U)$ its zero set:
$$
  V({\cal U})=\left\{\,\ab\in \R^{12}\ |\
  U_{ij}(\ab)=0, \ \forall\ U_{ij}( a )\in \cal U\,\right\}.
$$
\bdf Let $U_1,U_2,\ldots, U_s$ be $GL$-comitants of  systems
\eqref{2l1} and homogeneous polynomials in the coefficients of
these systems.
 A $GL$-comitant  $U(a ,x,y)$ of  systems \eqref{2l1}
 is called a conditional $\ T$-comitant (or
$CT$-comitant)  modulo $\left<U_1,U_2,...,U_s\right>$ (i.e. modulo
the ideal  generated by $U_{ij}(a)$ $(i=1,\ldots,s; j=0,1,\ldots,
d_i)$ in the ring $\R[a]$)  if the following two conditions are
satisfied:

 (i) the algebraic subset $V({\cal U})\subset \R^{12}$ is $Aff(2,\R)$-invariant
(see Definition \ref{def_G});

 (ii) for every ($\tau,\  \ab )\in T(2,\R)\times V(\cal U)$\
 we have
$ U(r_\tau( \ab) ,\ \tilde x,\,\tilde y)=
 U( \ab ,\ \tilde x,\,\tilde y)\
\mbox{in}\ \R[\tilde x,\,\tilde y]. $
\edf
\bdf
 A polynomial $U(a ,x,y)\in \R[a ,x,y]$, homogeneous of even degree in $x$, $y$
 has  well determined sign on \mbox{$V\subset\R^{12}$} with respect to $x,\,y$
 if for every $ \ab\in V$, the binary form $u(x,y)=U(\ab,x,y)$ yields a
 function of constant sign on $\R^2\setminus\{u=0\}$.
\edf
%%%%%%%%%%%%%
\bobs We draw the attention to the fact, that if a $CT$-comitant
$U( a ,x,y)$ of systems \eqref{2l1} of even weight is a binary
form of even degree in $x$ and $y$ and of  even degree in $a$ and
also has well determined sign on some $Aff(2,\R)$-invariant
algebraic subset $V$, then this sign is conserved after  an affine
transformation of the plane $x,y$  and time rescaling. \eobs
%%%%%%%%%%%%%%
We now construct  polynomials $D(a,x,y)$ and $H(a,x,y)$ which will
be shown in Lemma \ref{Table:Propreties} to be $T$-comitants.
\bnot\label{not:1}  Consider the
polynomial  $\Phi_{\alpha,\beta}=\alpha P+\beta Q\in
\R[a,X,Y,Z,\alpha,\beta]$ where $P=Z^2p(X/Z,Y/Z),$
$Q=Z^2q(X/Z,Y/Z)$, $p,$ $q\in \R[a,x,y]$ and $\max
(\deg_{(x,y)}p,\deg_{(x,y)}q)=2$. Then
$$
\bal
\Phi_{\alpha,\beta}=&\ c_{11}(a,\alpha,\beta)X^2  +2
c_{12}(a,\alpha,\beta)XY
       + c_{22}(a,\alpha,\beta)Y^2+
   2c_{13}(a,\alpha,\beta)XZ\\ &
   +2c_{23}(a,\alpha,\beta)YZ
 +c_{a,33}(\alpha,\beta)Z^2,\qquad
\Delta(a,\alpha,\beta) =\ \det\left|\left|c_{ij}(a,\alpha,\beta)
    \right|\right|_{i,j\in\{1,2,3\}},\\
&  D(a,x,y) = 4\Delta(a,-y,x),\qquad H(a,x,y)  =
4\big[\det\left|\left|c_{ij}(a,-y,x)
    \right|\right|_{i,j\in\{1,2\}}\big].\\
\eal
$$
\enot
%%%%%%%%%%%%%%%%
\bpr \label{pr:D_C}
 Consider $m\le3$ distinct directions in the affine
plane, where by direction we mean a point $(u,v)\in
\C^2\setminus(0,0)$. For the existence of an invariant straight
line of a system $S$ of coefficients $\ab$ corresponding to each
one of these directions it is necessary that there exist $m$
distinct common factors of the polynomials $C_2(\ab,x,y)$ and
$D(\ab,x,y)$ over $\C$.
\epr
\BProof Suppose that ${\cal L}(x,y)\equiv ux+vy+w=0$ is an invariant line for a
quadratic  system corresponding to $\ab\in \R^{12}$. Then we must
have $r,s,t\in \C$ such that
\be\label{2l2}
      \frac{\partial {\cal L}}{\partial x}p(x,y)+
      \frac{\partial {\cal L}}{\partial y}q(x,y)={\cal L}(x,y)(rx+sy+t).
\ee
Hence
$$
      up(x,y)+v q(x,y)=(ux+vy+w)(rx+sy+t).
$$
So $\Phi_{u,v}(\ab,x,y)=0$ is a reducible conic which occurs if
and only if the respective determinant   $\Delta(\ab,u,v)=0$. But
$4\Delta(\ab,u,v)=D(\ab,-v,u)=0$. The point at infinity of ${\cal
L}=0$ is $[-v:u:0]$ and so $C_2(\ab,-v,u)=0$. Hence, the two
homogeneous polynomials of degree 3 in $x$, $y$ must have the
common factor $ux+vy$. \EProof
\brm\label{rm:H}
Consider two parallel invariant affine lines \mbox{${\cal
L}_i(x,y)\equiv ux+vy+w_i=0$}, $(u,v)\ne(0,0)$, ${\cal
L}_i(x,y)\in \C[x,y],\ (i=1,2)$ of a quadratic system $S$ of
coefficients $\ab$. Then \mbox{$H(\ab,\!-v,u)\!=\!0$}, i.e. the
T-comitant $H(a,x,y)$ can be used for determining the directions
of parallel invariant lines of systems
\eqref{2l1}.
\erm
Indeed, according to (\ref{2l2}) from the hypothesis we must have
$$
      up(x,y)+ vq(x,y)=(ux+vy+w_1)(ux+vy+w_2).
$$
Therefore for the quadratic form in $x$ and $y$: $F_2(\ab,x,y)
=up_2(\ab,x,y)+vq_2(\ab,x,y)$ we obtain $F_2=(ux+vy)^2$ and hence
 $\Discriminant(F_2)=0$. Then calculations yield:
 $\Discriminant(F_2(\ab,x,y))=-H(\ab,-v,u)$ and hence
$H(\ab,-v,u)=0$.

We construct the  following polynomials which will be shown in
Lemma \ref{Table:Propreties} to be $T$-comitants:
\bnot\label{not1}
\be\label{Comit:Bi}
\bal
&B_3(a,x,y)=(C_2,D)^{(1)}=Jacob\left( C_2,D\right),\\
&B_2(a,x,y)=\left(B_3,B_3\right)^{(2)} - 6B_3(C_2,D)^{(3)},\\
&B_1(a)=\Res_x\left( C_2,D\right)/y^9=-2^{-9}3^{-8}\left(B_2,B_3\right)^{(4)}.\\
\eal
\ee
\enot
\bpr \label{pr:BGI}
    Suppose  $\tilde d=\deg \gcd\left(C_2(\ab,x,y),D(\ab,x,y)\right)$.
Then:
$$
  \bal
    \tilde d=0\ &\ \Leftrightarrow\ \ B_1(\ab)\ne0;\\
    \tilde d=1\ &\ \Leftrightarrow\ \ B_1(\ab)=0,\ B_2(\ab,x,y)\ne0;\\
    \tilde d=2\ &\ \Leftrightarrow\ \ B_2(\ab,x,y)=0,\ B_3(\ab,x,y)\ne0;\\
    \tilde d=3\ &\ \Leftrightarrow\ \ B_3(\ab,x,y)=0.
  \eal
$$
\epr
\BProof Since the  polynomial $B_3(a)$ is  the  Jacobian of the cubic binary
forms $C_2(a,x,y)$ and $D(a,x,y)$ we conclude that  $\tilde d=3$
if and only if  $B_3(\ab,x,y)=0$. We assume  $B_3(\ab,x,y)\ne0$
(i.e. $\tilde d\le2$) and consider the two subcases:
$B_2(\ab,x,y)=0$ and $B_2(\ab,x,y)\ne0$.

{\bf 1)} Assuming that $B_2(\ab,x,y)=0$ then $\tilde d=2$. Indeed,
suppose $\tilde d<2$. From
\eqref{Comit:Bi} the condition $B_2=0$ yields $B_1=0$ and since
the  polynomial $B_1(a)$ is  the resultant of the binary forms
$C_2(a,x,y)$ and $D(a,x,y)$ we get $\tilde d=1$, i.e. these
polynomials have a common linear factor $ax+by$. We may assume
$b=0$ (the case $b\ne0$ can be reduced to this one    via the
transformation $x_1=ax+b$, $y_1=x$). Then
$$
  C_2=x(a_1x^2+b_1xy +c_1y^2)\equiv x\tilde A(x,y),\quad
  D=x(a_2x^2+b_2xy +c_2y^2)\equiv x\tilde B(x,y).
$$
Considering \eqref{Comit:Bi}, calculations yield $\quad
    B_2(\ab,x,y)=  3x^4\cdot \Res_x (\tilde A,\tilde B)/y^4
$\quad and we obtain a contradiction:  since $B_2=0$ according to
\cite{Walker} (see Theorem 10.7 on page 29) the polynomials
$\tilde A$ and $\tilde B$ have a common nonconstant factor, i.e.
$\tilde d>1$. Conversely, suppose that $\tilde d=2$. Then clearly
we have
$$
  C_2=(ax+by)\tilde C,\qquad  D=(cx+dy)\tilde C
$$
and taking into account \eqref{Comit:Bi}  calculations yield
$B_2=0$.

 {\bf 2)} Let us assume now that the condition  $B_2(\ab,x,y)\ne0$
holds. Then $\tilde d\le1$ and since the  polynomial $B_1(a)$ is
the resultant of the binary forms $C_2(a,x,y)$ and $D(a,x,y)$ we
get $\tilde d=1$ if and only if $B_1(a)=0$.  \EProof

  From the Propositions \ref{pr:D_C} and  \ref{pr:BGI}\  the next result follows:
\bcrl\label{lm:BGI} For the existence of an invariant straight line in one
(respectively 2 or 3 distinct ) directions in the affine plane it is necessary
that $B_1=0$ (respectively $B_2=0$ or $B_3=0$).
\ecrl

Let us apply a translation $x=x'+x_0$, $y=y'+y_0$ to the
polynomials $p(a,x,y)$ and $q(a,x,y)$. We obtain $ \tilde p(\tilde
a(a,x_0,y_0),x',y')=p(a, x'+x_0, y'+y_0),$ $\quad \tilde q(\tilde
a(a,x_0,y_0),x',y')=q(a, x'+x_0, y'+y_0).$ Let us construct the
following polynomials
$$
\bal
\Gamma_i(a,x_0,y_0)& \equiv  \Res_{x'}
    \Big(C_i\big(\tilde a(a,x_0,y_0),x',y'\big),C_0\big(\tilde
    a(a,x_0,y_0),x',y'\big)\Big)/(y')^{i+1},\\
             &   \Gamma_i(a,x_0,y_0) \in \R[a,x_0,y_0],\ (i=1,2).\\
\eal
$$
\bnot\label{not2}
\be\label{2l4a}
   \tilde{\cal E}_i(a,x,y)=\left.\Gamma_i(a,x_0,y_0)\right|_{\{x_0=x,\ y_0=y\}}\in \R[a,x,y]
    \ \ (i=1,2).
\ee
\enot
\bobs\label{obs_} It can easily be checked using the Definition
\ref{def:comit}  that the constructed polynomials $\tilde{\cal
E}_1(a,x,y)$ and $\tilde{\cal E}_2(a,x,y) $ are affine comitants
of systems
\eqref{2l1} and are homogeneous polynomials in coefficients
$a_{00},\ldots, b_{02}$ and non-homogeneous in $x,y$ and $\
  \deg_{a} \tilde{\cal E}_1=3,\  \deg_{\,(x,y)} \tilde{\cal
  E}_1=5,\ \
    \deg_{a} \tilde{\cal E}_2=4,\  \deg_{\,(x,y)} \tilde{\cal E}_2=6.
$ \eobs
\bnot\label{GCD:Ei}
Let  ${\cal E}_i(a,X,Y,Z)$ $(i=1,2)$ be the homogenization of
$\tilde{\cal E}_i(a,x,y)$, i.e.
$$
{\cal E}_1(a,X,Y,Z)=Z^5\tilde {\cal E}_1(a,X/Z,Y/Z),\qquad {\cal
E}_2(a,X,Y,Z)=Z^6\tilde {\cal E}_1(a,X/Z,Y/Z)
$$
and $ \qquad
 {\cal H}(a,X,Y,Z)=\gcd\Big({\cal E}_1(a,X,Y,Z),\
 {\cal E}_2(a,X,Y,Z)\Big).
$
\enot
\indent In what follows  we shall examine the geometrical meaning
of  these affine comitants. We shall prove the following theorem:
\bth\label{theor:E1,E2}
 The straight line ${\cal L}(x,y)\equiv ux+vy+w=0$ $u,v,w\in \C$, $(u,v)\ne(0,0)$
is an invariant line for a system \eqref{2l1} in $\QS$
corresponding to a point $\ab\in\R^{12}$ if and only if the
polynomial ${\cal L}$ is a common factor of the polynomials
$\tilde{\cal E}_1(\ab,x,y)$ and $\tilde{\cal E}_2(\ab,x,y)$ over
$\C$,
 i.e.
$$
\tilde{\cal E}_i(\ab,x,y)=(ux+vy+w)\widetilde W_i(x,y)\in
\C[x,y]\quad (i=1,2).
$$
\eth

To prove this Theorem we first prove the following lemma:
\blm\label{lem:C0,C1,C2}
 The straight line $\tilde {\cal L}(x,y)\equiv ux+vy=0$
is an invariant line of a  system \eqref{2l1}  of coefficients
$\ab$ with $\ab_{00}^2+\bb_{00}^2\ne0$ if and only if $
C_0(\ab,-v,u)=0,$ $C_1(\ab,-v,u)=0,$ and $ C_2(\ab,-v,u)=0$. These
condition are equivalent to the following ones:
\be
      \Res_x(C_0(a,x,y),\,C_1(a,x,y))/y^2\Big|_{(\ab)}=0=\Res_x(C_0(a,x,y),C_2(a,x,y))/y^3\Big|_{(\ab)}.
\label{2l3}
\ee
\elm
\BProof According to Definition \ref{df1} the line $\tilde {\cal L}(x,y)$=0 is
a particular algebraic integral for a system \eqref{2l1} if and
only if the identity (\ref{2l2}) holds for this system and this
line. So in this case
$$
      u(p_0(\ab)+ p_1(\ab,x,y)+p_2(\ab,x,y))+v(q_0(\ab)+q_1(\ab,x,y)+q_2(\ab,x,y)) = (ux+vy)(S_0+S_1(x,y)),
$$
for some  $S_0\in\C$ and $S_1\in\C[x,y]$. Herein we obtain:
\beq
     &(i)\ \ & up_0(\ab)+vq_0(\ab)=0;  \\
     &(ii)\ & up_1(\ab,x,y)+vq_1(\ab,x,y)=(ux+vy)S_0(\ab);  \\
     &(iii) & up_2(\ab,x,y)+vq_2(\ab,x,y)=(ux+vy)S_1(\ab,x,y).
\eeq
We observe that, if $x=-v$ and $y=u$ then the left-hand sides of
$(i)$, $(ii)$ and $(iii)$  become $C_0(\ab,-v,u)$,\
$C_1(\ab,-v,u)$ and $C_2(\ab,-v,u)$, respectively. At the same
time the right-hand sides of these identities vanish. Therefore
the following equations are obtained:
\be
  C_0(\ab,-v,u)=0,\ C_1(\ab,-v,u)=0,\ C_2(\ab,-v,u)=0.\label{2l4}
\ee
As the degree of  $C_0(a,x,y)$ is one,  the relations (\ref{2l4})
hold. \EProof

\textit{ Proof of the Theorem \ref{theor:E1,E2}:} Consider the
straight line ${\cal L}(x,y)=0$. Let $(x_0,y_0)\in \R^2$ be any
fixed non-singular point of the systems
\eqref{2l1} (i.e. $p(x_0,y_0)^2+q(x_0,y_0)^2\ne0$) which lies on
the line ${\cal L}(x,y)=0$, i.e.
 $ux_0+vy_0+w=0$. Let $\tau_0$ be the translation $x=x'+x_0$, $y=y'+y_0$, $\tau_0(x',y')=(x,y)$.
Then
$$
{\cal L}(x,y)= {\cal L}(x'+x_0,y'+y_0)=ux'+vy'\equiv \tilde{\cal
L}(x',y')
$$
and consider the line $ux'+vy'=0$. By Lemma \ref{lem:C0,C1,C2} the
straight line $\tilde {\cal L}(x',y')=0$ will be an invariant line
of systems ($\ref{2l1}{}^{\tau_0}$) if and only if the conditions
(\ref{2l3}) are satisfied for these systems, i.e.\ $
  \Gamma_1(\ab,x_0, y_0))=\Gamma_2(\ab,x_0,y_0)=0
$\ for each point $(x_0,y_0)$ situated on the line ${\cal L}(x,y)\equiv
ux+vy+w=0$, since the relation  $ux_0+vy_0+w=0$ is satisfied.

Thus we have\quad $ \Gamma_i(\ab,x_0, y_0)= (ux_0+vy_0+w)\tilde
\Gamma_i(a,x_0, y_0) \ \ (i=1,2)$.\quad Taking into account the
notations (\ref{2l4a}) we conclude that the statement of Theorem
\ref{theor:E1,E2} is true.\EProof

We now  consider the possibility for a  straight line  to be a multiple
invariant  line.

\blm\label{lm3} If ${\cal L}(x,y)\equiv ux+vy+w=0$, $u,v,w\in \C$, $(u,v)\ne(0,0)$ is an invariant
straight line of multiplicity $k$ for a quadratic  system
\eqref{2l1} then $[{\cal L}(x,y)]^k\mid \gcd(\tilde {\cal
E}_1,\tilde {\cal E}_2)$, i.e. there exist $W_i(\ab,x,y)\in
\C[x,y]$\ $(i=1,2)$ such that
\be
\tilde {\cal E}_i(\ab,x,y)= (ux+vy+w)^k W_i(\ab,x,y),\quad i=1,2.
\label{2l5}
\ee
\elm
\BProof  Suppose that line ${\cal L}(x,y)\equiv ux+vy+w=0$ is an invariant line
of multiplicity $k$ for  a system \eqref{2l1} which corresponds to
point $\ab\in \R^{12}$. Let  us denote by $\ab_\varepsilon\in
\R^{12}$ the point corresponding to the perturbed system
$(\ref{2l1}_\varepsilon)$, which has $k$ distinct invariant lines:
${\cal L}_{i\varepsilon}(x,y)$ $(i=1,2,...k)$.

  According to Theorem \ref{theor:E1,E2} for systems $(\ref{2l1}_\varepsilon)$ the following relations
  are valid:
$$
\tilde {\cal E}_{j\varepsilon}(\ab_\varepsilon,x,y) = {\cal
L}_{1\varepsilon}\cdot {\cal L}_{2\varepsilon}...\cdot {\cal
L}_{k\varepsilon}\widetilde W_j(\ab_\varepsilon,x,y),\quad j=1,2,
$$
and according to Definition \ref{def:multipl} when perturbation
$\varepsilon\to 0$ then $ {\cal L}_{i\varepsilon}(x,y) \to {\cal
L}(x,y),\ \forall i=1,..k. $\ At the same time $
 \tilde {\cal E}_{j\varepsilon} \to \tilde {\cal E}_j = {\cal L}(x,y)^k W_j,\quad
 j=1,2.$\quad
 Lemma \ref{lm3} is proved. \EProof
\bcrl\label{Mult:Z=0}
 If the line $l_\infty:Z=0$ is of multiplicity $k>1$ then
$Z^{k-1}\mid \gcd({\cal E}_1, {\cal E}_2)$.
\ecrl
Indeed, suppose that the line $l_\infty:Z=0$ is of multiplicity
$k>1$ for  a system $S$ which corresponds to a point $\ab\in
\R^{12}$. Then by Definition \ref{def:Z-mult} there exist a
perturbed system $S_\varepsilon $ corresponding to the point
$\ab_\varepsilon\in \R^{12}$ which  has $k-1$ distinct invariant
affine straight lines: ${\cal
L}_{i\varepsilon}(x,y)=u_{i\varepsilon} x+ v_{i\varepsilon} y+
w_{i\varepsilon}$, $(u_{i\varepsilon},v_{i\varepsilon})\ne(0,0)$,
$(u_{i\varepsilon},v_{i\varepsilon},w_{i\varepsilon})\in\C^3$
$(i=1,2,...k-1)$
 such that for every $i$:\ $\displaystyle{\lim_{\varepsilon\to 0}
 (u_{i\varepsilon},v_{i\varepsilon},w_{i\varepsilon})=(0,0,1)}.$

By  Lemma \ref{lm3}  each of the $k-1$ affine lines ${\cal
L}_{i\varepsilon}$ must be a factor of the polynomial ${\cal
H}(\ab_\varepsilon,X,Y,Z)=\gcd\left({\cal
E}_1(\ab_\varepsilon,X,Y,Z),{\cal
E}_2(\ab_\varepsilon,X,Y,Z)\right)$. Therefore we conclude that
for the system $S$ we have\ $ Z^{k-1}\mid {\cal H}(\ab,X,Y,Z)$.

 As a next step we shall determine necessary conditions for the
existence of parallel invariant lines.  Let us consider the
following $GL$-comitants of systems \eqref{2l1}:
\bnot\label{not3}
$$
\ba{ll}
   M(a,x,y) = 2\,\Hess\big(C_2(x,y)\big), &
    \eta(a) = \Discriminant\big(C_2(x,y)\big),\\
   K(a,x,y) = \Jacob\big(p_2(x,y),q_2(x,y)\big),\qquad &
   \mu(a) =  \Discriminant\big(K(a,x,y)\big),\\
   N(a,x,y) =  K(a,x,y) + H(a,x,y), &
   \theta(a)  =   \Discriminant\big(N(a,x,y)\big),
\ea
$$
\enot
\noindent the geometrical meaning of which is revealed in the next
3 lemmas below.
\blm\label{lem:K,mu} Let $S\in \QS$ and let $\ab\in \R^{12}$ be its 12-tuple of
coefficients.  The common points of $P=0$ and $Q=0$ on the line $Z=0$ are given
by the common linear factors over $\C$ of $p_2$ and $q_2$. This yields the
geometrical meaning of  the T-comitants $\mu(a)$ and $K(a,x,y)$:
\beq
\deg\gcd(p_2(x,y),q_2(x,y)) =\left\{\begin{array}{lcl}
           0 & iff &\mu(\ab)\ne0;\\
           1 & iff &\mu(\ab)=0,\ K(\ab,x,y)\not=0;\\
           2 & iff &K(\ab,x,y)=0.
                         \end{array}\right.
\eeq
\elm
The proof follows from the fact that $K$ is the Jacobian of
$p_2(x,y)$ and $q_2(x,y)$ (i.e. $p_2$ and $q_2$ are proportional
if and only if $K(\ab,x,y)=0$ in $\R[x,y]$) and
$\mu=\Res_x(p_2,q_2)/y^4$.

 We shall prove the following assertion:
\blm\label{lm4}
 A necessary condition for  the existence of one
couple (respectively, two couples) of parallel invariant straight
lines of a systems \eqref{2l1} corresponding to $\ab\in\R^{12}$ is
the condition $\theta(\ab) =0$ (respectively, $N(\ab,x,y)=0$).
\elm
\BProof Let\quad
$ {\cal L}_i(x,y)\equiv ux+vy+w_i=0$, $(u,v)\ne(0,0)$,
$(u,v,w_i)\in\C^3$ $(i=1,2) $\ be two distinct $(w_1\ne w_2)$
parallel invariant lines for a quadratic system
\eqref{2l1}. Then by (\ref{2l2}) we have
$$
      up(x,y)+vq(x,y)=\xi(ux+vy+w_1)(ux+vy+w_2)
$$
and via a time rescaling we may assume  $\xi$= 1.  Therefore for
the quadratic homogeneities we obtain
\be
  (u\,a_{20} +v\,b_{20} )x^2+2(u\,a_{11} +v\,b_{11})xy +(u\,a_{02} +v\,b_{02})y^2=
  (u x +v y)^2, \label{(2l6}
\ee
and hence, for the existence of parallel invariant lines the
solvability of the following systems of quadratic equations with
respect to parameters $u$ and $v$ is necessary:
\be\label{2l7}
   (A_1) \    u\, a_{20} + v\, b_{20}= u^2;\qquad
   (A_2) \     u\,a_{11} + v\,b_{11} = uv;\qquad
   (A_3) \    u\,a_{02} + v\,b_{02} = v^2.
\ee
Without loss of generality we may consider $uv\ne0$, otherwise a
rotation of phase plane can be done. We now consider $vA_1 - u
A_2$ and $uA_3 - v A_2$:
$$
\bal
 &vA_1 - u A_2: \quad &  -a_{11}u^2 + (a_{20}-b_{11})uv + b_{20}v^2=0, \\
& uA_3 - v A_2:   &\quad   a_{02}u^2 + (b_{02}-a_{11})uv -
b_{11}v^2=0.
\eal
$$
Let $F_1(u,v)$ and $F_2(u,v)$ be the left hand sides of the above
equations. Clearly, for the existence of two directions
$(u_1,v_1)$ and $(u_2,v_2)$ such that in each of them there are 2
parallel invariant straight lines of a system
\eqref{2l1} it is necessary that the ${\rm rank}(U)=1$, where
$$
    U=\left(\ba{ccc}  -a_{11} &\ a_{20}\!-\!b_{11}\ & b_{20} \\
                     a_{02} &\ b_{02}\!-\!a_{11}\ & -b_{11} \\
           \ea\right).
$$
Hence, it is necessary
$$
\tilde A=\left|\!\!\ba{cc}  -a_{11} &\ a_{20}\!-\!b_{11}\  \\
                         a_{02} &\ b_{02}\!-\!a_{11}\  \\
           \ea\!\!\right|=0,\qquad
\tilde B=\left|\!\!\ba{cc}  -a_{11} &\ b_{20}  \\
                         a_{02} &\ -b_{11}  \\
           \ea\!\!\right|=0,\qquad
\tilde C=\left|\!\!\ba{cc}    a_{20}\!-\!b_{11}\ & b_{20}  \\
                          b_{02}\!-\!a_{11}\ & -b_{11} \\
           \ea\!\!\right|=0.
$$
Since the resultant of the binary forms $F_1(u,v)$ and $F_2(u,v)$
is $ \Res_u(F_1,F_2)/v^4=\tilde B^2-\tilde A\tilde C, $ we
conclude that for the existence of  one couple of parallel
invariant lines it is necessary that $\tilde B^2-\tilde A\tilde
C=0$. On the other hand calculations yield\qquad $
  N(\ab,x,y)= \tilde Cx^2 + 2 \tilde B xy + \tilde Ay^2,$ \qquad
 $ \theta = 4({\tilde B}^2-\tilde A\tilde C)$ \quad
and this completes the proof of lemma.
 \EProof

\blm\label{lm_3:2}
The type of the divisor $D_S(C,Z)$ for systems (\ref{il1}) is
determined by the corresponding conditions indicated in Table 1,
where we write $w_1^c+w_2^c+w_3$ if two of the points, i.e.
$w_1^c, w_2^c$, are complex but not real.
        Moreover, for each type of the divisor $D_S(C,Z)$ given
by Table 1 the quadratic systems (\ref{il1}) can be brought via a
real linear transformation to one of the following  canonical
systems $(\SSS_{I})-(\SSS_{V})$ corresponding to their behavior at
infi\-ni\-ty.
%\newpage
\begin{table}[!htb]
\begin{center}
\begin{tabular}{|c|c|c|c|}
\multicolumn{3}{r}{\bf Table  1}\\[1mm]
\hline
  \raisebox{-0.7em}[0pt][0pt]{Case}  & \raisebox{-0.7em}[0pt][0pt]{Type of $D_S(C,Z)$}
      & Necessary and sufficient   \\[-1mm]
         & & conditions on the comitants \\
 \hline\hline
 \rule{0pt}{1.2em} $1$ & $w_1+w_2+w_3 $ &  $\eta>0 $ \\[1mm]
\hline
 \rule{0pt}{1.2em}$2$  & $w_1^c+w_2^c+w_3 $ &  $\eta<0$ \\[1mm]
\hline
 \rule{0pt}{1.2em}  $3$ & $2w_1+w_2 $ &  $\eta=0,\quad M\ne0$ \\[1mm]
\hline
 \rule{0pt}{1.2em} $4$ & $3w $ &  $ M=0,\quad C_2\ne0$ \\[1mm]
\hline
 \rule{0pt}{1.2em} $5$ & $D_S(C,Z)$ undefined  &  $ C_2=0$ \\[1mm]
\hline
\end{tabular}
\end{center}
$$
\bal
&\left\{\ba{rcl}
 \displaystyle \frac{dx}{dt}&=&k+cx+dy+gx^2+(h-1)xy,\\[2mm]
 \displaystyle \frac{dy}{dt}&=& l+ex+fy+(g-1)xy+hy^2;
\ea\right. &\qquad (\SSS_I)\\[3mm]
&\left\{\ba{rcl}
 \displaystyle \frac{dx}{dt}&=&k+cx+dy+gx^2+(h+1)xy,\\[2mm]
 \displaystyle \frac{dy}{dt}&=& l+ex+fy-x^2+gxy+hy^2;
\ea\right.& \hspace{2cm}(\SSS_{I\!I})\\[3mm]
&\left\{\ba{rcl}
 \displaystyle \frac{dx}{dt}&=&k+cx+dy+gx^2+hxy,\\[2mm]
 \displaystyle \frac{dy}{dt}&=& l+ex+fy+(g-1)xy+hy^2;
\ea\right.&\qquad (\SSS_{I\!I\!I})\\[3mm]
&\left\{\ba{rcl}
 \displaystyle \frac{dx}{dt}&=&k+cx+dy+gx^2+hxy,\\[2mm]
 \displaystyle \frac{dy}{dt}&=& l+ex+fy-x^2+gxy+hy^2,
\ea\right.&\qquad (\SSS_{I\!V})\\[3mm]
&\left\{\ba{rcl}
 \displaystyle \frac{dx}{dt}&=&k+cx+dy+x^2,\\[2mm]
 \displaystyle \frac{dy}{dt}&=& l+ex+fy+xy.
\ea\right.&\qquad (\SSS_{V})
\eal
$$

\end{table}
\smallskip\noindent
\elm
%%%%%%%%%
\BProof The Table 1 follows easily from the definitions of  $\eta(a)$ and  $
M(a,x,y)$ in Notation \ref{not3}.
It is well known that a cubic binary form in $x,y$ over $\R$ can
be brought via a real linear transformation of the plane $(x,y)$:
${\frak g}(x,y) = (\tilde x,\tilde y)$ to one of the following
four canonical forms forms
\be\label{can_forms}
I.\ \tilde y(\tilde x-\tilde y); \qquad II.\ \tilde x(\tilde
x^2+\tilde y^2); \qquad III.\ \tilde x^2\tilde y; \qquad IV.\
\tilde x^3; \qquad V.\ 0.
\ee
 Let us consider a system  (\ref{il1}) corresponding to a point $\ab\in\R^{12}$
and let us consider the $GL$-comitant ${C_2}( \ab
,x,y)=yp_2(\ab,x,y)-xq_2(\ab,x,y)$ simply as a cubic binary form
in $x$ and $y$. Then the transformed binary form $g{C_2}(
\ab,x,y)=C_2(\ab,{\frak g}^{-1}(\tilde x,\tilde y))$ is one of the
canonical forms
\eqref{can_forms}  corresponding to cases indicated in Table 1.

  On the other hand, according to the
Definition \ref{def:comit} of a $GL$-comitant, for ${C_2}(\ab
,x,y)$ whose weight $\chi=-1$, we have for the same linear
transformation ${\frak g}\in GL(2,\R)$
$$
  C_2(r_{\frak g}(\ab),\, {\frak g}(x,y))= \det(g)\, C_2( \ab,\, x,y).
$$
Using ${\frak g}(x,y) = (\tilde x,\tilde y)$ we obtain $
C_2(r_{\frak g}(\ab),\, \tilde x,\tilde y)=\det(g) C_2( \ab ,\,
{\frak g}^{-1}(\tilde x,\tilde y)),$ \  where we may assume
$\det(g)=1$ via the rescaling:
 $\tilde x\to \tilde x/\det(g)$, \ $\tilde y\to\tilde y/\det(g)$. Thus, recalling that
$$
p_2(\tilde x,\tilde y)=\tilde a_{20}\tilde x^2+2\tilde
a_{11}\tilde x\tilde y+\tilde a_{02}\tilde y^2,\qquad q_2(\tilde
x,\tilde y)=\tilde b_{20}\tilde x^2+2\tilde b_{11}\tilde x\tilde
y+\tilde b_{02}\tilde y^2,
$$
  for the first canonical form in (\ref{can_forms}) we have
$$
C_2(\tilde \ab,\tilde x,\tilde y)=-\tilde b_{20}\tilde x^3+(\tilde
a_{20}-2\tilde b_{11})\tilde x^2\tilde y+(2\tilde a_{11}-\tilde
b_{02})\tilde x\tilde y^2+\tilde a_{02}\tilde y^3= \tilde x\tilde
y(\tilde x-\tilde y).
$$
Identifying the coefficients of the above identity we get the
canonical form $(\SSS_{I})$.

 Analogously for the cases $II,\ III$
and $IV$ we obtain the canonical form $(\SSS_{I\!I})$,
$(\SSS_{I\!I\!I})$ and $(\SSS_{I\!V})$ associated to the
respective polynomials in (\ref{can_forms}).

Let us consider the case $V$, i.e. $C_2(\ab,x,y)=0$ in $\R[x,y]$. Then we
obtain the systems
$$
 \frac{dx}{dt}=k+cx+dy+gx^2+hxy,\quad
  \frac{dy}{dt}= l+ex+fy+gxy+hy^2
$$
 with $g^2+h^2\ne0$. By interchanging $x$ and $y$ we may assume $g\ne0$ and
then via the linear transformation $\tilde x=g x+h y,$ $\tilde
y=y$ we obtain the systems $(S_{V})$. \EProof

In order to determine the existence of a common factor of the
polynomials ${\cal E}_1(\ab,X,Y,Z)$ and ${\cal E}_2(\ab,X,Y,Z)$ we
shall use the notion of the  resultant   of two polynomials with
respect to a given indeterminate (see for instance,
\cite{Walker}).

Let us  consider two polynomials $f,g\in R[x_1,x_2,\ldots,x_r]$
where $R$ is a unique factorization domain. Then we can regard the
polynomials $f$ and  $g$  as polynomials in $x_r$ over the ring
${\cal R}= R[x_1,x_2,\ldots,x_{r-1}]$, i.e.
$$
\bal
&f(x_1,x_2,\ldots,x_r)=a_0+a_1x_r+\ldots+a_nx_r^n,\\
&g(x_1,x_2,\ldots,x_r)=b_0+a_1x_r+\ldots+b_mx_r^m \quad
a_i,b_i\in{\cal R}.
\eal
$$
\blm\label{Trudi:2}{\rm \cite{Walker}} Assuming $n,m>0$,
$a_nb_m\ne0$
 the resultant $\Res_{x_r}(f,g)$ of the polynomials $f$
and $g$ with respect to $x_r$ is a polynomial in
$R[x_1,x_2,\ldots,x_{r-1}]$ which is zero if and only if $f$ and
$g$ have a common factor involving $x_r$.
 \elm
 We also shall use the
following remark:
\brm\label{rem:transf} Assume $s,\, \gamma\in \R$, $\gamma>0$. Then the transformation
$x=\gamma^{s}x_1$, $y=\gamma^{s}y_1$ and \mbox{$t=\gamma^{-s}t_1$}
does not change the coefficients of the  quadratic part of a
quadratic system, whereas each coefficient of the linear
(respectively
 constant ) part will be multiplied by  $\gamma^{-s}$
 (respectively
by $\gamma^{-2s}$).
\erm
%%%%%%%%%%%%%%%%%%%%%%%%%%%%%%%%%%%
%%%%%%%   New Section   %%%%%%%%%%%
%%%%%%%%%%%%%%%%%%%%%%%%%%%%%%%%%%%
\vspace{-3mm}
\section{The configurations  of invariant lines  of quadratic\\ differential
systems with $ M_{{}_\IL}=6$ }\label{Sec:m_il:6}%%% {Config}
%%%%%%%%%%%%%%%%
\bnot
We denote by $\QSL_{\bf6}$  the class of  all  quadratic
differential systems
\eqref{2l1} with $p,$ $q$ relatively prime $((p,q)=1)$, $Z\nmid C$ and  possessing a
configuration of 6 invariant straight lines including the line at
infinity and including possible multiplicities.
\enot
\blm\label{lm_3:1} For a quadratic system $S$ in $\QSL_{\bf6}$
the conditions $N(\ab,x,y)=0$  and $B_3(\ab,x,y)=0$ in $\R[x,y]$,
are satisfied.
\elm
\BProof Indeed, if for a system \eqref{2l1} the condition $M_{{}_\IL}=6$ is
satisfied, then taking into account the Definition \ref{def:multipl} we
conclude that there exists a perturbation  of the coefficients of the
 system~\eqref{2l1} within the class of quadratic systems such that
the perturbed systems has $6$ distinct invariant lines (real or
complex, including the  line $Z=0$). Hence, the   perturbed
systems must possess 2 couples of parallel lines with distinct
directions and an  additional line in a third direction. Then, by
continuity and according to Lemma  \ref{lm4} and Corollary
\ref{lm:BGI} we have $B_3(\ab,x,y)=0$ and   $N(\ab,x,y)=0$.
\EProof

 By Theorem \ref {theor:E1,E2} and Lemma \ref{lm3}  we obtain the
following result:
\blm\label{gcd:5} If $M_{{}_\IL}=6$ then\quad $\deg\gcd\big({\cal
E}_1,(\ab,X,Y,Z), {\cal E}_2(\ab,X,Y,Z)\big)=5$,\ i.e.
 ${\cal E}_1\mid {\cal E}_2$.
\elm

\bth\label{th_mil_6}
(i) The class $\QSL_{\bf6}$ splits into 11 distinct subclasses
indicated in {\bf Diagram 1} with the corresponding Configurations
6.1-6.11 where the complex invariant straight lines are indicated
by dashed lines. If an invariant straight line has multiplicity
$k>1$, then the number $k$ appears near the corresponding straight
line and this line is in bold face. We indicate next to the real
singular points their multiplicities as follows:
$\left(I_w(p,q)\right)$ if $w$ is a finite singularity,
$\left(I_w(C,Z),\ I_w(P,Q)\right)$ if $w$ is an infinite
singularity with $I_w(P,Q)\ne0$ and $\left(I_w(C,Z)\right)$ if $w$
is an infinite singularity with $I_w(P,Q)=0$.

 (ii) We consider the orbits of the
class $\QSL_{\bf6}$ under  the action of the real affine group and
time rescaling. The systems {\sl(VI.1)} up to {\sl(VI.11)} from
the Table 2 form a system of representatives of these orbits under
this action. A differential system  $(S)$  in $\QSL_{\bf6}$ is in
the orbit of  a system  belonging to $(VI.i)$ if and only if
$B_3(\ab,x,y)=0=N(\ab,x,y)$ and the corresponding conditions in
the middle column (where the polynomials $H_i$ $(i=1,2,3)$ and
$N_j$ $(j=1,\ldots,4)$ are   $CT$-comitants  to be introduced
below) is verified for this system $(S)$.  The conditions
indicated in the middle column are affinely invariant.

Wherever we have a case with invariant straight lines of
multiplicity greater than one,  we indicate the corresponding
perturbations proving this  in the Table 3.
\begin{figure}%[h]
\centerline{\bf Diagram 1 $(M_{{}_\IL}=6)$} \vspace{3mm}
\centerline{\psfig{figure=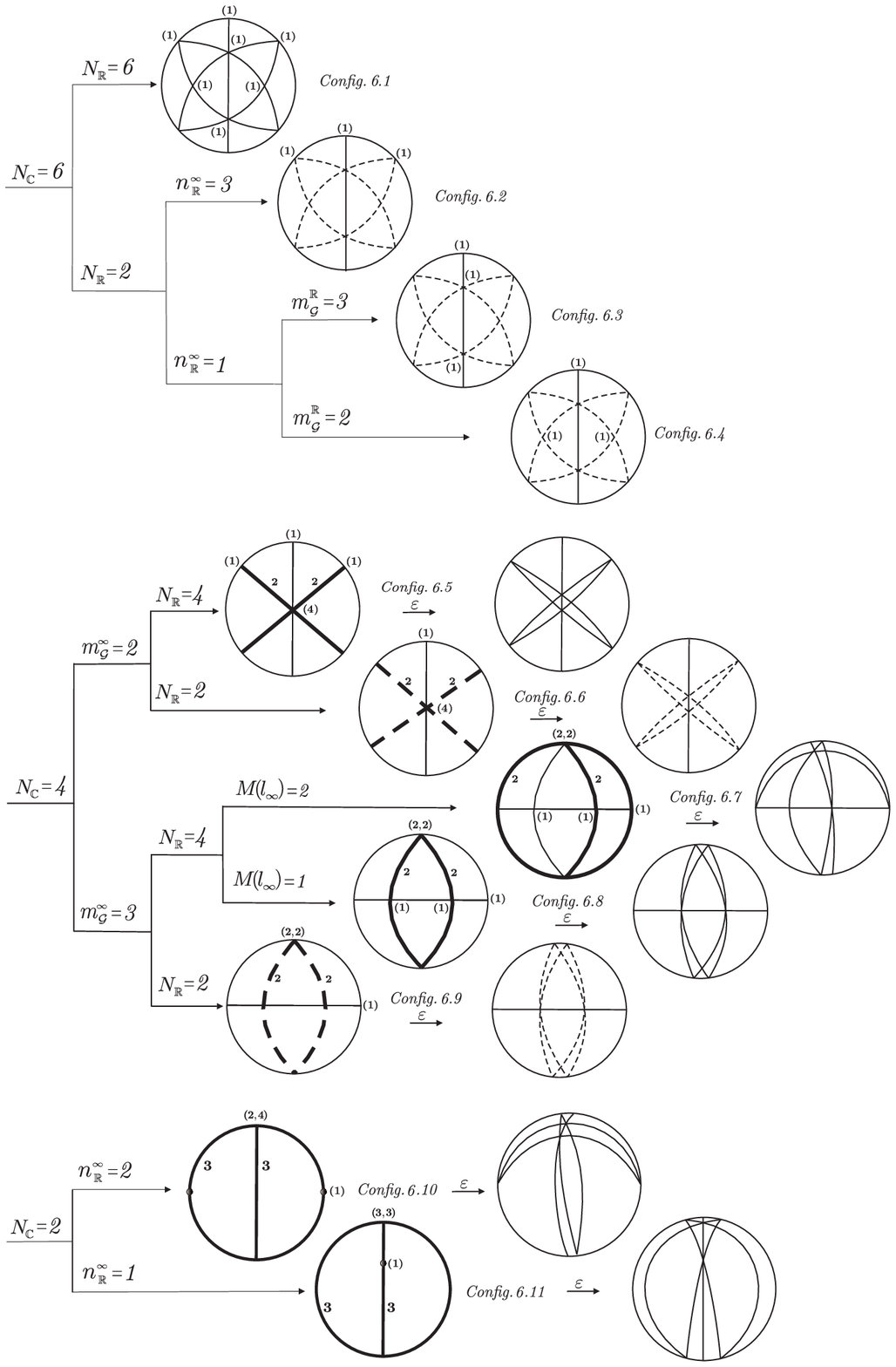}}
\vspace{2mm}
\end{figure}
\eth
\begin{table}[!htb]
\begin{center}
\begin{tabular}{|l|c|c|}
\multicolumn{3}{r}{\bf Table 2}\\[2mm]
\hline \raisebox{-0.7em}[0pt][0pt]{\qquad Orbit representative}  &
Necessary and sufficient
           & \raisebox{-0.7em}[0pt][0pt]{Configuration} \\
                             & conditions & \\
 \hline\hline
({\it VI.1})\,:\,$\ba{l} \dot x=x^2-1,\
                \dot y=y^2-1  \ea$ & $\ba{l}\eta>0,\ H_1>0 \ea $ & {\it Conf\mbox{}ig.\ 6.1}\\[1mm]
\hline ({\it VI.2})\,:\,$\ba{l} \dot x= x^2+1,\
                \dot y= y^2+1  \ea$ & $\ba{l}\eta>0,\ H_1<0 \ea\ $ & {\it Conf\mbox{}ig.\ 6.2}\\[1mm]
\hline ({\it VI.3})\,:\,$\ba{l} \dot x= 2xy,\
                \dot y=y^2-x^2-1\!\!  \ea$ & $ \ba{l}\eta<0,\ H_1<0  \ea $ & {\it Conf\mbox{}ig.\ 6.3}\\[1mm]
\hline ({\it VI.4})\,:\,$\ba{l} \dot x= 2xy,\
                \dot y= y^2-x^2+ 1\!\!  \ea$ & $ \ba{l}\eta<0,\ H_1>0 \ea $ & {\it Conf\mbox{}ig.\ 6.4}\\[1mm]
\hline ({\it VI.5})\,:\,$\ba{l} \dot x= x^2,\
                \dot y= y^2 \ea$ & $ \ba{l}\eta>0,\ H_1=0 \ea $ & {\it Conf\mbox{}ig.\ 6.5} \\[1mm]
\hline ({\it VI.6})\,:\,$\ba{l} \dot x= 2xy,\
                \dot y= y^2-x^2 \ea$ & $ \ba{l}\eta<0,\ H_1=0 \ea $ &  {\it Conf\mbox{}ig.\ 6.6}  \\[1mm]
\hline ({\it VI.7})\,:\,$\ba{l} \dot x=x^2-1,\
                \dot y= 2y \ea$ & $  MD\ne0,  \eta\!=\!H\!=\!N_1\!=\!N_2=0 $ & {\it Conf\mbox{}ig.\ 6.7 }\\[1mm]
\hline ({\it VI.8})\,:\,$\ba{l} \dot x=1 - x^2,\
                \dot y= -2xy  \ea$ & $ MH\ne0,\eta=H_2=0, H_3>0 $ &  {\it Conf\mbox{}ig.\ 6.8} \\[1mm]
\hline ({\it VI.9})\,:\,$\ba{l} \dot x=-1 - x^2,\
                \dot y= -2xy  \ea$ & $  MH\ne0, \eta=H_2=0, H_3<0  $ &  {\it Conf\mbox{}ig.\ 6.9} \\[1mm]
\hline (\!{\it VI.10})\,:\,$\ba{l} \dot x= x^2,\
                \dot y= 1  \ea $ & $  M\!\ne\!0,\eta\!= \!H\!=\!D\!=\!N_1\!=\!N_2\!=\!0 $ &  {\it Conf\mbox{}ig.\ 6.10} \\[1mm]
\hline (\!{\it VI.11})\,:\,$\ba{l} \dot x= x,\
                \dot y= y-x^2 \ea $ & $ \ba{l}\eta= M= N_3=N_4=0\ea $ &  {\it Conf\mbox{}ig.\ 6.11} \\[1mm]
\hline
\end{tabular}
\end{center}
\end{table}
\begin{table}[!htb]
\begin{center}
\begin{tabular}{|l|l|}
\multicolumn{2}{r}{\bf Table 3}\\[2mm]
\hline
 \hfil Perturbations\hfil &\hfil Invariant straight lines\hfil \\
 \hline\hline
({\it VI.5${}_\varepsilon$}) :\ $\ba{l} \dot x=
x^2-\varepsilon^2,\
                \dot y= y^2-\varepsilon^2 \ea$ & $ \ba{l} x=\pm \varepsilon,\ y=\pm \varepsilon,\ y=x \ea $ \\[1mm]
\hline ({\it VI.6${}_\varepsilon$}) :\ $\ba{l} \dot x= 2xy,\
                \dot y=\varepsilon^2 -x^2+ y^2 \ea $ & $ \ba{l} x=0,\ y\pm ix=\varepsilon,\ y\pm ix=-\varepsilon \ea $ \\[1mm]
\hline ({\it VI.7${}_\varepsilon$}) :\ $\ba{l} \dot x=-1+ x^2,\
                \dot y= 2y(\varepsilon y+1) \ea$ & $ \ba{l} x=\pm 1,\, \varepsilon
                                          y=-1,\, y=0,\, x-2 \varepsilon y=1\!\!\ea $ \\[1mm]
\hline ({\it VI.8${}_\varepsilon$}) :\ $\ba{l} \dot x=1 - x^2,\
                \dot y= -2xy-\varepsilon y^2  \ea$& $ \ba{l} y=0,\ x=\pm 1, \
                                            x+ \varepsilon y=\pm 1\ea $ \\[1mm]
\hline ({\it VI.9${}_\varepsilon$}) :\ $\ba{l} \dot x=-1 - x^2,\
                \dot y= -2xy-\varepsilon y^2  \ea$& $ \ba{l} y=0,\ x=\pm i, \
                                            x+ \varepsilon y=\pm i\ea $ \\[1mm]
\hline \!\!({\it VI.10${}_\varepsilon$})\,:\ $\bigg\{\ba{l} \dot
x=(1-\varepsilon)^2x^2-\varepsilon^2,\\
                \dot y= \left (2\,\varepsilon ^2 y+1\right )
  \left (2\,\varepsilon y+1\right ) \ea $ & $ \ba{l}
                              (1-\varepsilon)x=\pm
                              \varepsilon,\, 2\varepsilon y=-1,\, 2\varepsilon^2 y=-1,\!\!\!\\
                           (\varepsilon -1)^2 x-4\varepsilon^3 y-\varepsilon (\varepsilon +1)=0\ea $ \\[1mm]
\hline \!\!({\it VI.11${}_\varepsilon$})\,:\ $\bigg\{\ba{l} \dot
x=x+\varepsilon x^2,\\
                \dot y= y-x^2-2\varepsilon xy-2\varepsilon^2y^2 \ea $ & $ \ba{l}
                              x=0,\ \varepsilon x+1=0,\ x+ \varepsilon  y=0,\\
                           x+ 2\varepsilon  y=0,\ \varepsilon x+2\varepsilon^2y-1=0\ea $ \\[1mm]
\hline
\end{tabular}
\end{center}
\end{table}
\noindent {\it Proof of the Theorem \ref{th_mil_6}}:  According to
Table 1 we shall consider the subcases corresponding to distinct
types of the divisor $D_S(C,Z)$). Since  we only discuss the case
$C_2\ne0$, in what follows it suffices to consider only the
canonical forms $(\SSS_I)$ to $(\SSS_{I\,V})$.  The idea of the
proof is to perform a case by case discussion for each one of
these canonical forms, for which   according   to
Lemma~\ref{lm_3:1}  the conditions $B_3=0=N$ must be fulfilled.
These conditions   yield  specific conditions on the parameters.
The discussion proceeds further by breaking these cases in more
subcases determined by more restrictions on the parameters.
Finally we construct  new invariants or T-comitants which put
these conditions in invariant form.

For constructing the invariant polynomials included in the
statement of Theorem~\ref{th_mil_6} we shall use the $T$-comitants
$D(a,x,y)$ and $H(a,x,y)$ indicated  before as well as the
\mbox{$GL$-comitants} (\ref{C_i:D_i}).

%%%%%%%%%%%%%
\vspace{-3mm}
\subsection{Systems with the divisor\ $D_S(C,Z)=1\cdot w_1+1\cdot w_2+1\cdot w_3$}
\vspace{-3mm} For this case we shall later  need the following polynomial which
is shown to be an affine invariant in  Lemma
\ref{Table:Propreties}.
\bnot\label{not_H1} Let us denote\quad
$ H_1(a)= -\big((C_2,C_2)^{(2)},C_2)^{(1)},D\big)^{(3)}. $
\enot
\indent According to Lemma \ref{lm_3:2} a system with this type of
divisor can be brought by linear transformations  to the canonical
form $(\SSS_I)$ for which we have:
%%%%%%%%%%%%%%%
\beq
N(\ab,x,y)&=& (g^2-1)x^2+2(g-1)(h-1)xy+(h^2-1)y^2.
\eeq
Hence the condition $N=0$ yields $(g-1)(h-1)=g^2-1=h^2-1=0$ and we obtain 3
possibilities: $(a)\ g=1=h;$ $(b)\ g=1=-h;$ $(c)\ g=-1=-h$. The cases $(b)$ and
$(c)$ can be brought by linear transformations to the case $(a)$. Hence the
resulting polynomials are: $p_2(x,y)=x^2$  and $q_2(x,y)=y^2$. Then the term in
$x$ of the first equation and the term in $y$ in the second equation can be
eliminated via a translation. Thus we obtain the systems
\be\label{s4.1}
  \dot x=k + dy + x^2,\quad  \dot y=l + ex +y^2
\ee
for which we have\ $ B_3= 3[-e^2x^4+2e^2x^3y +4(l-k)x^2y^2
-2d^2xy^3+d^2y^4]. $\ Hence, the condition $B_3=0$ yields
$d=e=k-l=0$ and we get the systems of the form:
\be\label{CF_1}
  \dot x=l + x^2,\quad  \dot y=l + y^2.
\ee
By Remark \ref{rem:transf} ( $\gamma=|l|,\ s=1/2$) for systems (\ref{CF_1}) we
can consider  $l\in \{-1,0,1\}$. Clearly these systems possess the invariant
straight lines $x=\pm\sqrt{-l}$, $y=\pm\sqrt{-l}$, $y=x$. Therefore, we obtain
Config. 6.1 (respectively, Config. 6.2) for $l<0$ (respectively, for $l>0$) and
Config. 6.5 for $l=0$. For systems~(\ref{CF_1}) the affine invariant
$H_1(\ab)=-2^93^3l$ and, hence,  $\sign(l)=-\sign(H_1(\ab))$.
%%%%%%%%%%%%%%%
%%%%%%%%%%%%%%%
\vspace{-3mm}
\subsection{Systems with the divisor\ $D_S(C,Z)=1\cdot w^c_1+1\cdot w^c_2+1\cdot w_3$}
\vspace{-3mm} In this case  by Lemma \ref{lm_3:2} the systems
\eqref{2l1} can be brought by linear transformations to the
canonical form $(\SSS_{I\!I})$ for which we have: \beq
N(\ab,x,y)&=&(g^2-2h+2)x^2 +2g(h+1)xy+(h^2-1)y^2. \eeq Hence the
condition  $N=0$ yields $g=h-1=0$ and we may consider $c=d=0$ due
to the translation $x=x_1-d/2$, $y=y_1-c/2$. We thus obtain the
systems
\be\label{S2:N=0}
  \dot x= k +2xy,\qquad
  \dot y= l + ex +fy -x^2 + y^2
\ee
for which\
$
 B_3=6\,\left[(ef-2k)x^4+(f^2-e^2)x^3y-(4k+ef)x^2y^2-2ky^4\right].
$\ Hence, the condition $B_3=0$ yields $k=e=f=0$ and we obtain the
following form
\be\label{CF_3}
  \dot x= 2xy,\qquad
  \dot y= l -x^2 + y^2
\ee
where $l\in \{-1,0,1\}$  by the Remark \ref{rem:transf} ( $\gamma=|l|,\
s=1/2$).
 It is not
difficult to convince ourselves  that these systems possess  as invariant
straight lines the components  over $\C$ of:
$$
  x=0,\quad x^2+2\,i\, xy-y^2-l=0,\quad x^2-2\,i\,xy-y^2-l=0,
$$
with  the intersection points:\quad $
  p_{1,2}=(0,\pm\sqrt{-l}),\ \  p_{3,4}=(\pm\sqrt{l},0).
$\
 On the other hand for systems (\ref{CF_3}) we have
$ H_1= 2^{10}\,3^2\,l.$\ Therefore, if $H_1\ne0$  we get Config.
6.3 for $H_1<0$ and
 Config. 6.4 for $H_1>0$, whereas for $H_1=0$
we obtain Config. 6.6.
%%%%%%%%%%%%%%%%%%%%%%%%%%%%
%%%%%%%%%%%%%%%%%%%%%%%%%%%%
\vspace{-3mm}
\subsection{Systems with the divisor\ $D_S(C,Z)=2\cdot w_1+1\cdot w_2$}
\vspace{-3mm}  For this case we shall later  need the following polynomials which
are shown to be $CT$-comitants in  Lemma \ref{Table:Propreties}.
%%%%%%%% Notation %%%%%%%%
\bnot\label{not:H2,3-N1,2}
Let us denote
$$
\bal
&H_2(a,x,y)=(C_1,\ 2H\!-\!N)^{(1)}\!-\!2D_1N,\quad
N_1(a,x,y)=C_1(C_2,C_2)^{(2)} \!-\!2C_2(C_1,C_2)^{(2)},\\
&H_3(a,x,y)=(C_2,D)^{(2)}, \qquad
N_2(a,x,y)=D_1(C_1,C_2)^{(2)}\!-\!\Big((C_2,C_2)^{(2)},C_0\Big)^{(1)}.
\eal
$$
\enot
%%%%%%%%%  Observation  %%%%%%%%%%

We are in the case of the canonical form $(\SSS_{I\!I\!I})$ for which we have:
\be\label{eq_H}
 \bal
 N(a,x,y)&=(g^2-1)x^2 +2h(g-1)xy+h^2y^2,\\
 H(a,x,y)&=-(g-1)^2x^2 -2h(g+1)xy-h^2y^2.\\
\eal
\ee
The condition $N=0$ yields $h=g^2-1=0$ and we shall examine two
subcases: \mbox{$H(a,x,y)\not=0$} and $H(a,x,y)=0$.
%%%%%%%%%%%%
\vspace{-3mm}
\subsubsection{The case $H(a,x,y)\not=0$}
\vspace{-3mm}
In this case for $h=0$ we have
$H(a,x,y)=-(g-1)^2x^2\not=0$ and hence the condition $N=0$ yields
$g=-1$. Moreover, for systems $(\SSS_{I\!I\!I})$ we can consider
$e=f=0$ due to the translation of the origin of coordinates to the
point $(f/2,e/2)$. Thus, the systems $(\SSS_{I\!I\!I})$ can be
brought to the form
\be\label{S3_NM_}
  \dot x =k +cx +dy -x^2,\qquad
  \dot y =l - 2xy,
\ee
for which\quad $
 B_3= 6x(-2lx^3 +cd xy^2+d^2y^3).
$\quad  So, the condition $B_3=0$ yields $l=d=0$ and we obtain the
systems
\be\label{s4.4}
  \dot x =k +cx -x^2,\qquad
  \dot y = - 2xy
\ee
with $k\ne0$ (otherwise  systems (\ref{s4.4}) become degenerate).

So far we have only used the necessary conditions $N=0$ and
$B_3=0$ for this particular case. These are not sufficient for
having 6 invariant lines. According to Lemma \ref{gcd:5} we must
have ${\cal E}_1\mid {\cal E}_2$ (see Notation~\ref{GCD:Ei}).
 Calculations yield :
$$
 {\cal E}_1= (kZ^2-X^2){\cal H},\quad
 {\cal E}_2 =X(X^2-cXZ-kZ^2){\cal H},
\quad {\cal H}=\gcd\left({\cal E}_1,{\cal E}_2\right)=
2Y\left(kZ^2+cXZ-X^2\right).
$$
 Since $k\ne0$ according to  Lemma \ref{Trudi:2} we obtain the
 condition
$ \Res_Z( {\cal E}_1/{\cal H},\ {\cal E}_2/{\cal
H})=-c^2kX^6\equiv0 $ must hold.  This yields $c=0$ and the
systems (\ref{s4.4}) become
\be\label{CF_7}
  \dot x =k - x^2,\qquad
  \dot y = - 2xy.
\ee
By Remark \ref{rem:transf} ($\gamma=|k|,\ s=1/2$) we may assume
$k\in\{-1,1\}$.

For the systems (\ref{CF_7}) we have\ $ {\cal H}= \gcd\left({\cal
E}_1,{\cal E}_2\right)=2Y\left(kZ^2-X^2\right)^2$\ and according
to Lemma~\ref{lm3} each one of the two invariant lines
$x=\pm\sqrt{k}$ of the systems (\ref{CF_7}) could be of the
multiplicity two. And they  are indeed of  multiplicity two as it
is shown by the perturbations ({\it VI.8${}_\varepsilon$}) (for
$k=1$) and ({\it VI.9${}_\varepsilon$}) (for $k=-1$) from Table 3.
 Thus, we obtain Config. 6.8 for $k=1$ and Config. 6.9 for $k=-1$.

On the other hand for the systems (\ref{CF_7}) we have $
  H_2 = 16cx^2,$ $H_3 =32kx^2.$ Hence the
$T$-comitants $H_2 $ and $H_3 $ capture exactly the conditions
$c=0$ and  $k>0$ or $c=0$ and  $k<0$ and this leads to the
corresponding conditions in Table 2.
%%%%%%%%%%%%
\vspace{-3mm}
\subsubsection{The case $H(a,x,y)=0$}
\vspace{-3mm} According to (\ref{eq_H}) the conditions $N=H=0$
yield $h=0,$ $g=1$ and translating the origin of coordinates to
the point $(-c/2,0)$ the systems $(\SSS_{I\!I\!I})$ can be brought
to the form
\be\label{S3_NH_0}
\dot x=k+dy+ x^2,\qquad
\dot y= l+ex+fy.
\ee
 For these systems we have\
$
  B_3= 6dxy^2(fx-dy)
$ and the condition $B_3=0$ yields $d=0$. So, we obtain the  systems
\be\label{s4.5}
\dot x=k+ x^2,\qquad
\dot y= l+ex+fy
\ee
for which we have $D(\ab,x,y)=-f^2x^2y$.
%%%%%%%%%%%%%%%%

 {\bf 1)} If $D\ne0$  then $f\ne0$ and by Remark \ref{rem:transf}
($\gamma=f/2,\ s=1$)   we can consider $f=2$. Then via  the
translation we may assume $l=0$ and  we obtain the systems
\be\label{s4.6}
\dot x=k+ x^2,\qquad \dot y= ex+2y,
\ee
for which calculations yield
\be\label{val:Ei}
\bal
&{\cal E}_1= \Big[4 Y(X - Z) + e(X^2 - 2 X Z - kZ^2)\Big]{\cal H},
\quad
 {\cal E}_2=3(e X + 2Y)(X^2 + kZ^2){\cal H},\\
& {\cal H}= 2Z\left(X^2+kZ^2\right),\quad \Res_Y({\cal E}_1/{\cal
H},\ {\cal E}_2/{\cal H})=-2e(X^2 +  kZ^2)^2.
\eal
\ee
Hence  for ${\cal E}_1\mid {\cal E}_2$  the condition
$\Res_Y({\cal E}_1/{\cal H},\ {\cal E}_2/{\cal H})=0$ must be
fulfilled in $\R[X,Z]$. This  yields $e=0$ and  then  we obtain:\
$ \Res_X\Big( ({\cal E}_1/{\cal H})|_{e=0},({\cal E}_2/{\cal
H})|_{e=0}\Big)= 32(k+1)Y^3Z^2=0$. Hence $k+1=0$ and for $e=k+1=0$
we obtain the system
\be\label{CF_9}
\dot x= x^2-1,\qquad
\dot y= 2y,
\ee
for which ${\cal H}=\gcd\left({\cal E}_1,{\cal
E}_2\right)=YZ(X-Z)^2(X+Z)$. This system possesses the invariant
affine lines\ $
   x=\pm1,\quad y=0.
$ \ Moreover, taking into account the polynomial ${\cal H}$, by
Lemma \ref{lm3} and Corollary \ref{Mult:Z=0} the line $x=1$ as
well as the line $l_\infty:Z=0$ could be of multiplicity two. This
is confirmed by the  perturbations ({\it VI.7${}_\varepsilon$})
from Table 3. Since this system possesses only two finite
singularities $(\pm1,0)$ which are simple, we conclude that the
configuration of the invariant lines  of the system (\ref{CF_9})
is Config. 6.7.

It remains to observe that the conditions $e=0=k+1$ are equivalent
to $N_1=N_2=0$, as for systems (\ref{s4.5}) we have $ N_1=
8e\,x^4,\quad N_2=16(k+1)x $.

  {\bf 2)} The condition $D=0$ implies $f=0$ and we obtain the systems
\be\label{s4.7}
\dot x=k+ x^2,\qquad \dot y= l+ex.
\ee
\noindent Calculations yield:
\be\label{val:Eia}
\bal
&{\cal E}_1= \left[2\,l X Z + e(X^2 - k Z^2)\right]{\cal H},\quad
{\cal E}_2=(e X + l Z)(X^2 + kZ^2)\,{\cal H},
\eal
\ee
where $ {\cal H}= Z\left(X^2+kZ^2\right).$\  Hence for ${\cal
E}_1\mid {\cal E}_2$ according to Lemma \ref{Trudi:2} at least one
of the following conditions must hold:
$$
\Res_X({\cal E}_1/{\cal H},\ {\cal E}_2/{\cal H})=-4\, ek(e^2k +
l^2)^2Z^6=0,\quad \Res_Z({\cal E}_1/{\cal H},\ {\cal E}_2/{\cal
H})=-4\, ek(e^2k + l^2)^2X^6=0,
$$
and we obtain that either  $ek=0$ or $e^2k + l^2=0$. Since the
second case yields a degenerate system we obtain the necessary
condition $ek=0$. It is easy to observe that for $e^2+k^2\ne0$ we
obtain ${\cal E}_1\nmid {\cal E}_2$. Therefore $k=e=0$ (then
$l\ne0$) and
 via the additional  rescaling  $y\to l\,y$  we obtain the system:
\be\label{CF_10}
\dot x= x^2,\qquad \dot y= 1
\ee
for which ${\cal H}=\gcd\left({\cal E}_1,{\cal
E}_2\right)=X^3Z^2$. By Lemma \ref{lm3} and Corollary
\ref{Mult:Z=0} the line $x=0$ as well as the line $Z=0$ could be
of multiplicity three. This is confirmed by the   perturbations
({\it VI.10${}_\varepsilon$}) from Table 3. It remains to note
that for systems (\ref{s4.7}) we obtain \ $  N_1=8e\,x^4,$\ $
N_2=16kx $\ \ and, hence in this case we obtain Config. 6.10 if
and only if \mbox{$N_1=0=N_2.$}
%%%%%%%%%%%%%%%%%
\vspace{-3mm}
\subsection{Systems with the divisor\ $D_S(C,Z)=3\cdot w$}
\vspace{-3mm}  For this case we shall later  need the following
polynomials which are shown to be  $CT$-comitants in Lemma
\ref{Table:Propreties}.
\bnot\label{not_HC3} Let us denote\quad  $
  N_3=   \left(C_2,C_1\right)^{(1)},\ \
 N_4=   4\left(C_2,C_0\right)^{(1)} - 3C_1D_1.$
\enot
We are in the case of the canonical form $(\SSS_{I\!V})$ for which
we have:
$$
  N=(g^2-2h)x^2+2ghxy+h^2y^2.
$$
So, the condition $N=0$ yields $h=g=0$ and  due to the translation
$x=x_1+e/2,$ $y=y_1$  we may assume $e=0$. Hence the systems
$(\SSS_{I\!V})$ become
\be\label{S4_N0}
\dot x=k+cx+dy,\qquad
\dot y= l+fy-x^2,
\ee
for which\quad
$
  B_3= 6dx^3(fx-dy).
$\quad The condition $B_3=0$ yields $d=0$ and we shall examine the
systems of the form
\be\label{s4.8}
\dot x=k+cx,\qquad
\dot y= l+fy-x^2.
\ee
Calculations yield
\be\label{val_Eib}
\bal
 &{\cal E}_1= \left[(c+f)X^2+2kXZ+ (c-f)(fY+lZ)Z\right]{\cal H},\\
 &{\cal E}_2=Z(c X + k Z)^2\,{\cal H},\quad
 {\cal H}= Z^2(cX+kZ).\\
\eal
\ee
Since the polynomial ${\cal E}_2/{\cal H}$ depends only on $X$ and
$Z$ for ${\cal E}_1\mid {\cal E}_2$  the following condition must
hold: $f(c-f)=0$. We claim that for $f=0$ we cannot have ${\cal
E}_1\mid {\cal E}_2$. Indeed, assuming $f=0$ we obtain the
quadratic form ${\cal E}_1/{\cal H}=cX^2+2kXZ+ clZ^2$ in $X$ and
$Z$, which must divide $Z(cX+kZ)^2$. This clearly implies that the
discriminant of this form must be zero, i.e. $4(k^2-c^2l)=0$.
However  this leads to degenerate systems.

Therefore we must have $c-f=0$ and for the systems (\ref{s4.8})
with $f=c$ calculations yield: $ {\cal E}_1=X\,\tilde{\cal H},$  $
{\cal E}_2=Z(cX+kZ)\,\tilde{\cal H}, $ where $\tilde{\cal H}=Z^2(c
X + k Z)^2$. Therefore ${\cal E}_1\mid {\cal E}_2$ if and only if
$k=0$ and we obtain the systems\quad $ \dot x= cx,\quad \dot y=
l+cy-x^2.$\quad with $c\ne0$.  We may assume $c=1$ by  Remark
\ref{rem:transf} ( $\gamma=c,\ s=1$)  and via the translation of
the origin of coordinates to the point $(0,-l)$ we obtain $l=0$.
This leads to the following system
\be\label{CF_11}
  \dot x= x,\qquad \dot y= y -x^2,
\ee
with ${\cal H}=\gcd\left({\cal E}_1,{\cal E}_2\right)=X^3Z^2$ and
by Lemma \ref{lm3} and Corollary \ref{Mult:Z=0} each one of the
invariant lines $x=0$ and $Z=0$ is of multiplicity 3. This is
confirmed by the  perturbed systems ({\it VI.11${}_\varepsilon$})
from Table 3.\ On the other hand for  systems (\ref{s4.8})\quad $
 N_3=3(c-f)x^3,\qquad N_4=3x[4kx+(f^2-c^2)y]
$\ and hence, the conditions $c-f=k=0$ are equivalent to
$N_3=N_4=0$. Taking into account the existence of the simple
singular point $(0,0)$ placed on the line $x=0$ we obtain
Config.~6.11.

\smallskip
All the cases in  Theorem \ref{th_mil_6} are thus examined. To
finish the proof of the Theorem \ref{th_mil_6} it remains to show
that the conditions  occurring in the middle column of Table 2 are
affinely invariant. This follows from the proof of Lemma
\ref{Table:Propreties}.
  \EProof

%%%%
%%%%%%%  NEW SECTION  %%%%%%%%%
%%%%
\vspace{-3mm}
\section{ The configurations  of invariant lines  of quadratic differential
 systems with  $M_{{}_\IL}=5$ } \label{Sec:m_il:5}
\vspace{-3mm}
\bnot
We denote by $\QSL_{\bf5}$  the class of  all  quadratic
differential systems
\eqref{2l1} with $p,$ $q$ relatively prime $((p,q)=1)$, $Z\nmid C$ and
possessing a configuration of five invariant straight lines
including the line at infinity and including possible
multiplicities.
\enot
\blm\label{lm_NB2} If for a quadratic system $(S)$  $M_{{}_\IL}=5$, then for this system
one of the two following conditions are satisfied:\\[-5mm]
$$
\bal
&(i)\quad N(\ab,x,y)=0=B_2(\ab,x,y)\ \mbox{in}\ \R[x,y];\quad (ii)\quad
\theta(\ab)=0= B_3(\ab,x,y) \ \mbox{in}\ \R[x,y].
\eal
$$
\elm
\BProof Indeed, if for a system \eqref{2l1} the condition $M_{{}_\IL}=5$ is
satisfied then taking into account the Definition \ref{def:multipl} we conclude
that there exists a perturbation  of the coefficients of the
 system~\eqref{2l1} within the class of quadratic systems such that
the perturbed systems have five distinct invariant lines (real or
imaginary, including the  line $Z=0$). Hence, the   perturbed
systems must possess either 2 couples of parallel lines with
distinct directions or one couple of parallel lines and 2
additional lines with distinct directions.  Then, by continuity
and according to Lemma \ref{lm:BGI} and Corollary  \ref{lm4} we
respectively have
 either the conditions $(i)$ or $(ii)$. \EProof

 By Theorem \ref{theor:E1,E2} and Lemmas \ref{lm3}  and
\ref{lm_3:1} we obtain the following result:
\blm \label{gcd:4}
(a) If for a system $(S)$ of coefficients $\ab\in \R^{12}$,\
$M_{{}_\IL}=5$ then\\\hphantom{m} $\deg\gcd\big({\cal
E}_1(\ab,X,Y,Z), {\cal E}_2(\ab,X,Y,Z)\big)=4;$ \hphantom{mm} (b)
If $N(\ab,x,y)\not\equiv0$ then $M_{{}_\IL}\le5$.
\elm
 \bth\label{th_mil_5}
(i) The class $\QSL_{\bf5}$ splits into 30 distinct subclasses
indicated in {\bf Diagram 2} with the corresponding Configurations
5.1-5.30 where the complex invariant straight lines are indicated
by dashed lines. If an invariant straight line has multiplicity
$k>1$, then the number $k$ appears near the corresponding straight
line and this line in bold face. We indicate next to the singular
points their multiplicities as follows: $\left(I_w(p,q)\right)$ if
$w$ is a finite singularity, $\left(I_w(C,Z),\ I_w(P,Q)\right)$ if
$w$ is an infinite singularity with $I_w(P,Q)\ne0$ and
$\left(I_w(C,Z)\right)$ if $w$ is an infinite singularity with
$I_w(P,Q)=0$.

(ii) We consider the orbits of the class $\QSL_{\bf5}$ under  the
action of the real affine group and time rescaling. The systems
{\sl(V.1)} up to {\sl(V.30)} from the Table 4 form a system of
representatives of these orbits under this action. A differential
system $(S)$ in $\QSL_{\bf5}$ is in the orbit of   a system
belonging to  $(V.i)$ if and only if the corresponding conditions
in the middle column (where the polynomials $H_i$ $(i=7,\ldots,
11)$ and $N_j$ $(j=5,6)$ are   $CT$-comitants  to be introduced
below) are verified for this system $(S)$.  The conditions
indicated in the middle column are affinely invariant.

Wherever we have a case with invariant straight lines of
multiplicity greater than one,  we indicate the corresponding
perturbations in the Table 5.
\eth
\brm\label{rem:H>0}
We observe that in the middle column of the Table 5 (and of the
Table 2) there occur conditions of the form $\mathcal{M}(a,x,y)=0$
in $\R[x,y]$ or of the form $\mathcal{M}(a,x,y)>0$ (or $<0$),
where $\mathcal{M}(a,x,y)$ is a homogeneous polynomial in $a$ and
separately in $x$ an $y$, which is a $CT$-comitant. All
polynomials occurring in conditions of the second type are of even
weight, of even degree in $a_{00},\ldots,b_{02}$ and have a well
determined sign on the corresponding variety indicated in the
Lemma \ref{Table:Propreties}.
\erm

\begin{figure}%[h]
\centerline{\bf \hfil Diagram 2  $(M_{{}_\IL}=5)$}
\vspace{3mm}
\centerline{\psfig{figure=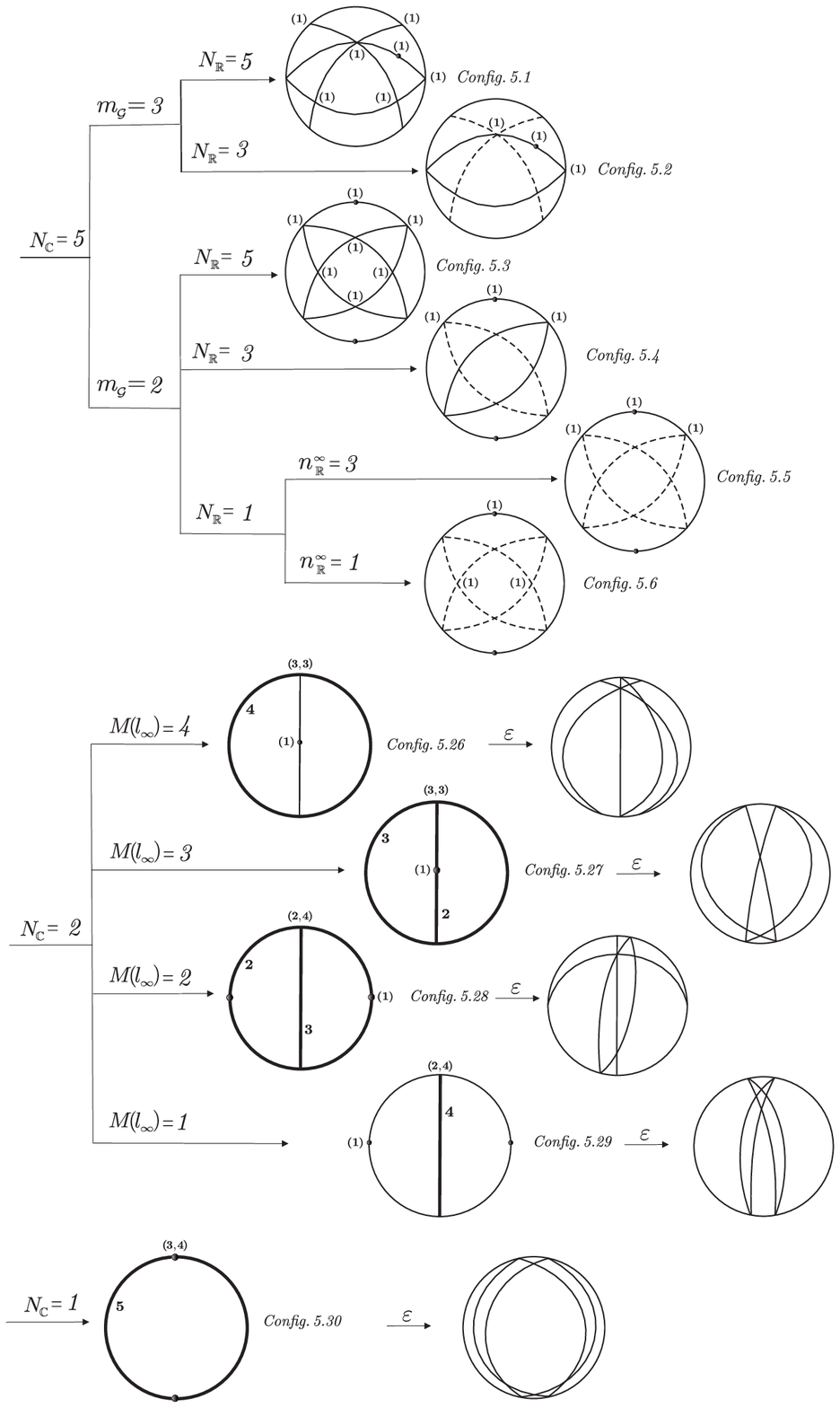}}
\end{figure}

\begin{figure}%[h]
\centerline{ \bf Diagram 2  $(M_{{}_\IL}=5)$\ {\it(continued)}}
\vspace{2mm} \psfig{figure=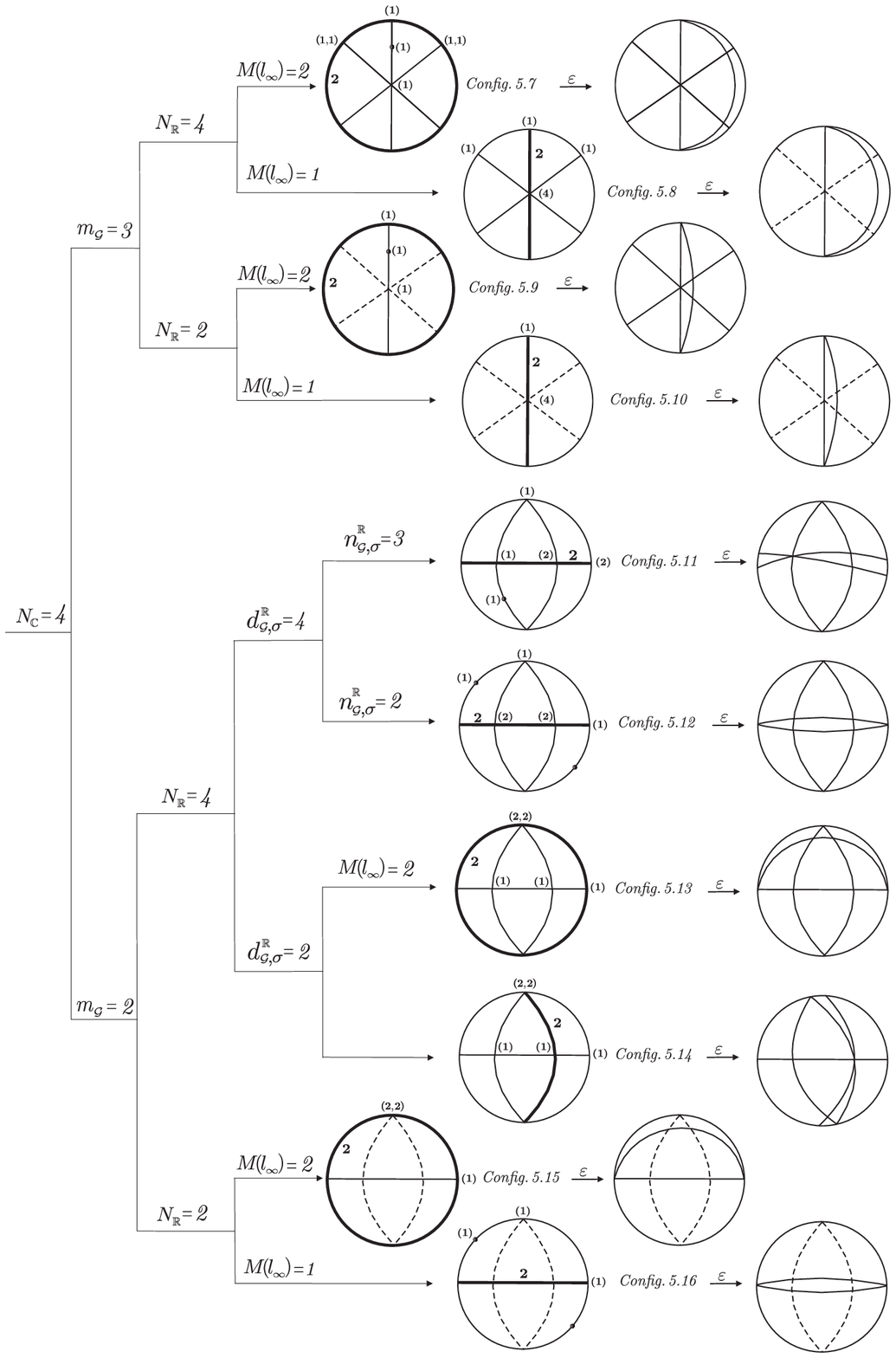}
\end{figure}

\begin{figure}%[h]
\centerline{\bf Diagram 2  $(M_{{}_\IL}=5)$\ {\it(continued)}}
\vspace{6mm} \psfig{figure=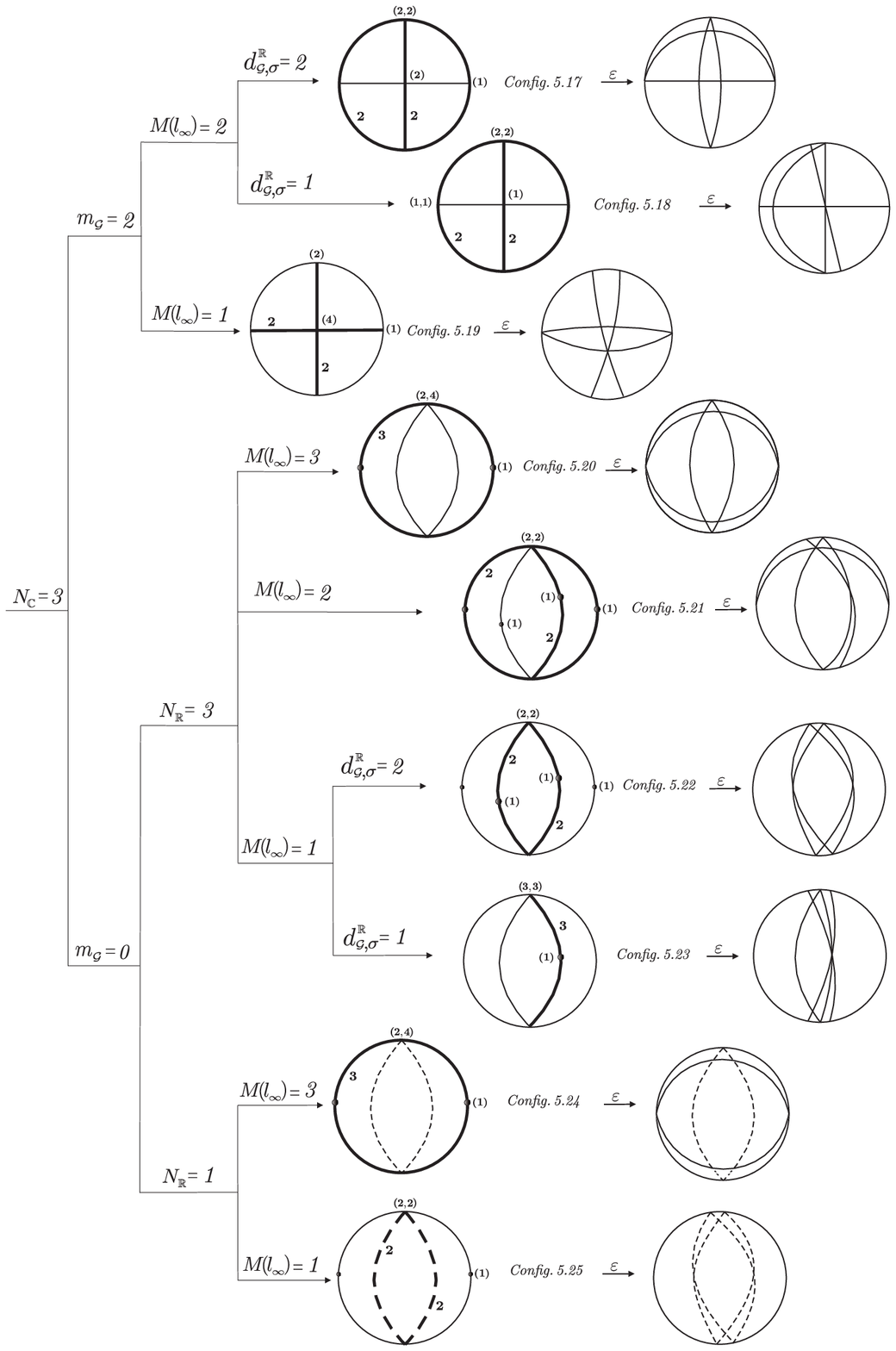}
\end{figure}

\begin{table}
\begin{center}
\begin{tabular}{|l|c|c|}
\multicolumn{3}{r}{\bf Table 4}\\[1mm]
\hline \raisebox{-0.7em}[0pt][0pt]{\qquad Orbit representative}  & Necessary
and sufficient
           & \raisebox{-0.7em}[0pt][0pt]{Configuration} \\[-1mm]
                             & conditions & \\
 \hline\hline\rule[0mm]{0mm}{10.0mm}
 ({\it V.1})\ $\left\{\!\!\ba{l} \dot x=(x+1)(gx+1), \\[-0.8mm]
                \dot y=(g-1)xy+y^2,  \\[-0.8mm]
                \hspace{8mm} g(g^2-1)\ne0 \ea\!\!\right. $ & $ \ba{c}\eta>0,\ B_3=\theta=0,\\
                                          N\ne0,\, \mu\ne0,\, H_1\ne0\ea $ & Config.\ 5.1\\[5.9mm]
\hline\rule[0mm]{0mm}{10.0mm}
({\it V.2})\ $\left\{\!\!\ba{l} \dot x=g(x^2-4),\ g\ne0\\[-0.8mm]
                \dot y=(g^2-\!4)\!+\!(g^2+\!4)x \\[-0.8mm]
               \hspace{7mm} -x^2+gxy-y^2   \ea\!\!\!\!\right. $ & $ \ba{c}\eta<0,\, B_3=\theta=0,\\
                        N\ne0,\, \mu\ne0,\, H_1\ne0\ea $ & Config.\ 5.2\\[5.8mm]
\hline%\rule[0mm]{0mm}
({\it V.3})\ $\bigg\{\!\!\ba{l} \dot x=-1+x^2,\\
                \dot y=g(y^2-1),\ g\ne0 \ea $ & $ \ba{c}\eta>0,\, B_2=N=0,\, B_3\ne0,\\
                                             H_1>0, \  H_4=0,\  H_5>0 \ea $ & Config.\ 5.3 \\%[4mm]
\hline%\rule[0mm]{0mm}
({\it V.4})\ $\bigg\{\!\!\ba{l} \dot x=-1+x^2,\\
                \dot y=g(1+y^2),\ \ g\ne0 \ea $ & $ \ba{c}\eta>0,\, B_2=N=0,\, B_3\ne0,\\
                                                 H_4=0,\  H_5<0 \ea $ & Config.\ 5.4\\%[4mm]
\hline%\rule[0mm]{0mm}
({\it V.5})\ $\bigg\{\!\!\ba{l} \dot x=1+x^2,\\
                \dot y=g(1+y^2),\ \ g\ne0 \ea $ & $ \ba{c}\eta>0,\, B_2=N=0,\, B_3\ne0,\\
                                              H_1<0,\   H_4=0,\  H_5>0 \ea $ & Config.\ 5.5 \\%[4mm]
\hline%\rule[0mm]{0mm}
({\it V.6})\ $\bigg\{\!\!\ba{l} \dot x= 1+2xy,\\[-0.8mm]
                \dot y=g-x^2+y^2,\ \ g\in\R \ea $ & $ \ba{c}\eta<0,\, B_3\ne0,\,
                                                       B_2=N=0 \ea $ & Config.\ 5.6 \\%[0mm]
\hline%\rule[0mm]{0mm}
({\it V.7})\ $\bigg\{\!\!\ba{l} \dot x= 1+x,\\[-0.8mm]
                \dot y=-xy+y^2  \ea $ & $ \ba{c}\eta>0,\, B_3=\theta=0,\\[-0.8mm]
            N\ne0,\, \mu= H_6\!=\!0\ea $ & $ \ba{l} \mbox{Config.\ 5.7}  \ea $ \\%[4mm]
\hline%\rule[0mm]{0mm}
({\it V.8})\ $\bigg\{\!\!\ba{l} \dot x= gx^2,\ \ g(g^2-1)\ne0\\
                \dot y=(g-1)xy+y^2   \ea $ & $ \ba{c} \eta>0,\, B_3=\theta=0,\\
            N\ne0,\, \mu\ne0,\, H_1=0 \ea $ & $ \ba{l} \mbox{Config.\ 5.8}\ea $ \\%[4mm]
\hline%\rule[0mm]{0mm}
({\it V.9})\ $\bigg\{\!\!\ba{l} \dot x= 2x,\\[-0.8mm]
                \dot y=1-x^2-y^2   \ea $ & $ \ba{c}\eta<0,\, B_3=\theta=0,\\[-0.8mm]
            N\ne0, \, \mu= H_6\!=\!0\ea $ & $ \ba{l} \mbox{Config.\ 5.9}\ea $ \\%[4mm]
\hline\!\!\!
({\it V.10})\ $\bigg\{\!\!\ba{l} \dot x= gx^2,\ \ g\ne0 \\
                \dot y=-x^2+gxy-y^2   \ea $ & $ \ba{c} \eta<0,\, B_3=\theta=0,\\
            N\ne0,\, \mu\ne0,\, H_1=0 \ea $ & $ \ba{l} \mbox{Config.\ 5.10}\ea $ \\%[4mm]
\hline\!\!\!
({\it V.11})\ $\bigg\{\!\!\ba{l} \dot x=x^2+xy,\\
                \dot y=y+y^2 \ea\!\! $ & $\!\! \ba{c} \eta=0,\,M\ne0,\, B_3=\theta=0,\\
             \mu\ne0,\, N\ne0,\,D\ne0 \ea\!\! $ & $ \ba{l} \mbox{Config.\ 5.11}\ea $ \\[4mm]
\hline\!\!\!
({\it V.12})\ $\bigg\{\!\!\ba{l} \dot x=-1+x^2,\\
                \dot y=y^2  \ea $ & $\!\! \ba{l} \eta>0,\, B_2=N=0,\,
             B_3\ne0,\\ H_1>0,\,
             H_4= H_5=0 \ea\!\! $ & $ \ba{l} \mbox{Config.\ 5.12}\ea $ \\%[4mm]
\hline\!\!\!
({\it V.13})\ $\bigg\{\!\!\ba{l} \dot x= g(x^2-1),\\
                \dot y=2y,\  \ g(g^2-1)\ne0   \ea $ & $\!\! \ba{c} \eta=0,\,M\ne0, B_3=N=0,\\
              H\!=\!N_1\!=\!0,\, N_2D\ne0,\, N_5>0 \ea\!\! $ & $ \ba{l} \mbox{Config.\ 5.13}\ea $ \\%[4mm]
\hline\!\!\!
({\it V.14})\ $\bigg\{\!\!\ba{l} \dot x=(x+1)(gx+1),\\
                \dot y=(g\!-\!1)xy, \
                 g(g^2\!-\!1)\!\ne\!0   \ea\!\! $ & $\!\! \ba{c} \eta=0,\,M\ne0, B_3=\theta=0,\\
             NK\ne0,\ \mu=H_6=0
                \ea\!\! $ & $ \ba{l} \mbox{Config.\ 5.14}\ea $ \\%[4mm]
\hline\!\!\!
({\it V.15})\ $\bigg\{\!\!\ba{l} \dot x= g(x^2+1),\\
                \dot y=2y,\ \ \ g\ne0   \ea $ & $\!\! \ba{c} \eta=0,\,M\ne0, B_3=N=0,\\
              H\!=\!N_1\!=\!0,\, N_2D\ne0,\, N_5<0 \ea\!\! $ & $ \ba{l} \mbox{Config.\ 5.15}\ea $ \\%[4mm]
\hline
\end{tabular}
\end{center}
\end{table}
%%%%%%%%%%
\begin{table}
\begin{center}
\begin{tabular}{|l|c|c|}
\multicolumn{3}{r}{\bf Table 4} {\it (continued)}\\[2mm]
\hline \raisebox{-0.7em}[0pt][0pt]{\qquad Orbit representative}  & Necessary
and sufficient
           & \raisebox{-0.7em}[0pt][0pt]{Configuration} \\
                             & conditions & \\
\hline
({\it V.16})\ $\bigg\{\!\!\ba{l} \dot x=1+x^2,\\
                \dot y=y^2  \ea $ & $\!\! \ba{c} \eta>0,\, B_2=N=0,\,
             B_3\ne0,\\ H_1<0,\,
             H_4= H_5=0 \ea\!\! $ & $ \ba{l} \mbox{Config.\ 5.16}\ea $ \\%[4mm]
\hline
({\it V.17})\ $\bigg\{\!\!\ba{l} \dot x= x^2,\\
                \dot y=2y \ea\!\! $ & $\!\! \ba{c} \eta=0,\,M\ne0, B_3=N=0,\\
              H\!=\!N_1\!=\!N_5\!=\!0,\,
            N_2D\ne0 \ea\!\! $ & $ \ba{l} \mbox{Config.\ 5.17}\ea $ \\[4mm]
\hline
({\it V.18})\ $\bigg\{\!\!\ba{l} \dot x=1+x,\\
                \dot y=-xy \ea\!\! $ & $\!\! \ba{c} \eta=0,\,M\ne0, B_3=\theta=0,\\
             N\!\ne\!0,\, \mu\!=\!K\!=\!H_6\!=\!0
                \ea\!\! $ & $ \ba{l} \mbox{Config.\ 5.18}\ea $ \\[4mm]
\hline
({\it V.19})\ $\bigg\{\!\!\ba{l} \dot x=x^2+xy,\\
                \dot y=y^2\ea\!\! $ & $\!\! \ba{c} \eta=0,\,M\ne0, B_3=\theta=0,\\
              \mu\ne0,\, N\ne0,\,D=0
                \ea\!\! $ & $ \ba{l} \mbox{Config.\ 5.19}\ea $ \\[4mm]
\hline
({\it V.20})\ $\bigg\{\!\!\ba{l} \dot x=-1+x^2,\\
                \dot y=1 \ea\!\! $ & $\!\! \ba{c} \eta\!=\!0,\,M\!\ne\!0, B_3\!=\!N\!=\!H\!=\!0,\\
              D=N_1=0,  N_2\ne0, N_5>0
                \ea\!\! $ & $ \ba{l} \mbox{Config.\ 5.20}\ea $ \\[4mm]
\hline
({\it V.21})\ $\bigg\{\!\!\ba{l} \dot x=-1+x^2,\\
                \dot y=x+2y  \ea\!\! $ & $\!\! \ba{c} \eta=0,\,M\ne0, B_3=N=0,\\
              H\!=\!N_2\!=\!0, D\ne0,\,N_1\ne0
                \ea\!\! $ & $ \ba{l} \mbox{Config.\ 5.21}\ea $ \\[4mm]
\hline
({\it V.22})\ $\bigg\{\!\!\ba{l} \dot x= 1-x^2,\\
                \dot y=1-2xy  \ea\!\! $ & $\!\! \ba{c} \eta=0,\,M\ne0, B_2=N=0,\\
              B_3\ne0,\,H_2=0,\,H_3>0
                \ea\!\! $ & $ \ba{l} \mbox{Config.\ 5.22}\ea $ \\[4mm]
\hline
({\it V.23})\ $\bigg\{\!\!\ba{l}\dot x=-1+x^2,\\
                \dot y=-3+y-x^2+xy  \ea\!\! $ & $\!\! \ba{c} \eta=M=0,\, N\ne0,\\
              B_3=\theta= N_6=0
                \ea\!\! $ & $ \ba{l} \mbox{Config.\ 5.23} \ea $ \\[4mm]
\hline
({\it V.24})\ $\bigg\{\!\!\ba{l} \dot x= 1+x^2,\\
                \dot y=1 \ea\!\! $ & $\!\! \ba{c} \eta\!=\!0,\,M\!\ne\!0, B_3\!=\!N\!=\!H\!=\!0,\\
              D=N_1=0,  N_2\ne0, N_5<0
                \ea\!\! $ & $ \ba{l} \mbox{Config.\ 5.24}\ea $ \\[4mm]
\hline
({\it V.25})\ $\bigg\{\!\!\ba{l}\dot x= -1-x^2,\\
                \dot y=1-2xy  \ea\!\! $ & $\!\! \ba{c} \eta=0,\,M\ne0, B_2=N=0,\\
              B_3\ne0,\,H_2=0,\,H_3<0
                \ea\!\! $ & $ \ba{l} \mbox{Config.\ 5.25}\ea $ \\[4mm]
\hline
({\it V.26})\ $\bigg\{\!\!\ba{l}\dot x=g-x,\ \ g\in\{0,1\}, \\
                \dot y=y-x^2  \ea\!\! $ & $\!\! \ba{c} \eta=M=0,\, N_3\ne0,\\
              B_3=N= D_1=0
                \ea\!\! $ & $ \ba{l} \mbox{Config.\ 5.26}\ea $ \\[4mm]
\hline
({\it V.27})\ $\bigg\{\!\!\ba{l}\dot x=1+x, \\
                \dot y=y-x^2   \ea\!\! $ & $\!\! \ba{c} \eta=M=0,\, \,N_4\ne0,\\
              B_3\!=\!N\!=\! N_3\!=\!0,\,D_1\ne0
                \ea\!\! $ & $ \ba{l} \mbox{Config.\ 5.27} \ea $ \\[4mm]
\hline
({\it V.28})\ $\bigg\{\!\!\ba{l} \dot x= x^2,\\
                \dot y=1+x  \ea\!\! $ & $\!\! \ba{c} \eta=0,\,M\ne0, B_3=N=0,\\
              H=D=N_2=0, N_1\ne0
                \ea\!\! $ & $ \ba{l} \mbox{Config.\ 5.28}\ea $ \\[4mm]
\hline
({\it V.29})\ $\bigg\{\!\!\ba{l}\dot x=-x^2,\\
                \dot y=1-2xy \ea\!\! $ & $\!\! \ba{c} \eta=0,\,M\ne0,\,
                B_2=N=0,\\    B_3\ne0,\, H_2=H_3=0
                \ea\!\! $ & $ \ba{l} \mbox{Config.\ 5.29}\ea $ \\[4mm]
\hline
({\it V.30})\ $\bigg\{\!\!\ba{l}\dot x=1,\ \ g\in\{-1,0,1\}, \\
                \dot y=g-x^2 \ea\!\! $ & $\!\! \ba{c} \eta=M=0,\, N_4\ne0,\\
              B_3=N =  N_3=D_1 =0
                \ea\!\! $ & $ \ba{l} \mbox{Config.\ 5.30} \ea $ \\[4mm]
\hline
\end{tabular}
\end{center}
\end{table}
\begin{table}[!htb]
\begin{center}
\begin{tabular}{|l|l|}
\multicolumn{2}{r}{\bf Table 5}\\[2mm]
\hline
 \hfil Perturbations\hfil & \hfil Invariant straight lines \hfil \\
 \hline\hline
({\it V.7${}_\varepsilon$})\,:$\ba{l} \dot x= (x+1)(\varepsilon
x+1),\   \dot y= (\varepsilon-1)xy+y^2\!\! \ea$ & $ \ba{l}y=0,\, x=-1,\, y-x=1,\,  \varepsilon x=-1 \ea $ \\[0mm]
 \hline
({\it V.8${}_\varepsilon$})\,:$\ba{l} \dot x=
(x+\varepsilon)(gx+\varepsilon),\   \dot y= (g-1)xy+y^2\!\! \ea$
 & $ \ba{l}y=0,\, x=-\varepsilon,\, y-x=\varepsilon,\,  gx=-\varepsilon  \ea $ \\[0mm]
 \hline
({\it V.9${}_\varepsilon$})\,: $\Big\{\ba{l} \dot
x=2x(\varepsilon x+1),\\[-1.4mm]   \dot y= 1+2\varepsilon x -x^2+2\varepsilon xy-y^2\!\! \ea$
& $ \ba{l}x=0,\ \varepsilon x=-1,\ y\pm i x+1=0 \ea $ \\[0mm]
 \hline
\!\!\!\! (\!{\it V.10${}_\varepsilon$})\,: $\Big\{\ba{l} \dot
x=4g\varepsilon^2+\varepsilon(g^2+4)x+gx^2,\\[-1.1mm]
\dot y= \varepsilon^2(4-g^2)-x^2+g xy-y^2\!\! \ea$
& $ \ba{l}x=-g\varepsilon,\  gx=-4\varepsilon,\\[-1.1mm] x+g\varepsilon=\pm i(y+2\varepsilon) \ea $ \\[0mm]
 \hline
\!\!\!\! (\!{\it V.11${}_\varepsilon$})\,:$\ba{l} \dot
x=\varepsilon x+x^2+(1+\varepsilon)xy,\ \dot y= y+ y^2\!\! \ea$
& $ \ba{l} x=0,\,  y=-1,\, y=0,\,  x+ \varepsilon y=-\varepsilon \ea $ \\[0mm]
 \hline
\!\!\!\! (\!{\it V.12${}_\varepsilon$})\,: $\ba{l} \dot x=x^2-1,\
\dot y= y^2-\varepsilon^2 \ea$
& $ \ba{l} x=\pm 1,\  y=\pm \varepsilon \ea $ \\[0mm]
 \hline
\!\!\!\! (\!{\it V.13${}_\varepsilon$})\,: $\ba{l} \dot
x=g(x^2-1),\ \dot y= 2y(\varepsilon y+1) \ea$
& $ \ba{l}y=0,\ x=\pm1,\  \varepsilon y=-1 \ea $ \\[0mm]
 \hline
\!\!\!\! (\!{\it V.14${}_\varepsilon$})\,:$\ba{l} \dot
x=(x\!+\!1)(g x\!+\!1),\, \dot y= (g\!-\!1)xy\!-\!\varepsilon
y^2\!\! \ea$
& $ \ba{l} x\!=\!-1,\, g x\!=\!-1,\,  y\!=\!0,\,   x\!+\! \varepsilon y\!=\!-1 \!\ea $ \\[0mm]
 \hline
\!\!\!\! (\!{\it V.15${}_\varepsilon$})\,: $\ba{l} \dot
x=g(x^2+1),\ \dot y= 2y(\varepsilon y+1) \ea$
& $ \ba{l}y=0,\ x=\pm i,\  \varepsilon y=-1 \ea $ \\[0mm]
 \hline
\!\!\!\! (\!{\it V.16${}_\varepsilon$})\,: $\ba{l} \dot x=x^2+1,\
\dot y= y^2-\varepsilon^2 \ea$
& $ \ba{l} x=\pm i,\  y=\pm \varepsilon \ea $ \\[0mm]
 \hline
\!\!\!\! (\!{\it V.17${}_\varepsilon$})\,:$\ba{l} \dot
x=x^2-\varepsilon^2 ,\ \dot y= 2y(\varepsilon y+1) \ea$
& $ \ba{l}y=0,\ x=\pm\varepsilon,\  \varepsilon y=-1 \ea $ \\[0mm]
 \hline
\!\!\!\! (\!{\it V.18${}_\varepsilon$})\,:$\ba{l} \dot
x=(x\!+\!1)(\varepsilon x\!+\!1),\, \dot y=
(\varepsilon\!-\!1)xy\!-\!\varepsilon y^2\!\! \ea$
& $ \ba{l} x\!=\!-1,\, \varepsilon x\!=\!-1,\,  y\!=\!0,\,   x\!+\! \varepsilon y\!=\!-1 \!\ea $ \\[0mm]
 \hline
\!\!\!\! (\!{\it V.19${}_\varepsilon$})\,: $\ba{l} \dot
x=\varepsilon^2 x+ x^2+(1+\varepsilon)xy,\ \dot y= \varepsilon y+
y^2\!\! \ea$
& $ \ba{l} x=0,\,  y=0,\, y=-\varepsilon,\,  x+ \varepsilon y=-\varepsilon^2 \!\!\ea $ \\[0mm]
 \hline
\!\!\!\! (\!{\it V.20${}_\varepsilon$})\,: $\ba{l} \dot x=x^2-1,\
\dot y= 1-\varepsilon^2 y^2\ea$
& $ \ba{l} x=\pm 1,\  \varepsilon y=\pm 1 \ea $ \\[0mm]
 \hline
\!\!\!\! (\!{\it V.21${}_\varepsilon$})\,: $\Big\{\ba{l} \dot
x=(x+1)(x+4\varepsilon x-1),\\[-1.1mm]
\dot y= (x+2y)(1+4\varepsilon y)\!\! \ea$
& $ \ba{l} x=-1,\  x(1+4\varepsilon)=1,\\[-1.1mm] 4\, \varepsilon y=-1,\  x- 8\, \varepsilon y=1 \ea $ \\[0mm]
 \hline
\!\!\!\! (\!{\it V.22${}_\varepsilon$})\,: $\ba{l} \dot x=1-x^2,\
\dot y= 1-2xy-\varepsilon y^2 \ea$
& $ \ba{l} x=\pm\, 1,\  x+ \varepsilon y=\pm\sqrt{1+\varepsilon} \ea $ \\[0mm]
 \hline
\!\!\!\! (\!{\it V.23${}_\varepsilon$})\,: $\Bigg\{\ba{l} \dot
x=(1+\varepsilon)(x-1+2\varepsilon)(x+1-2\varepsilon),\\[-1.1mm]
\dot y= (4\varepsilon^2-3)+(1+2\varepsilon)y-x^2\\[-1.1mm]
\qquad +(1-2\varepsilon)xy-2\varepsilon^2 y^2\!\!
\ea$
& $ \ba{l} x=\pm(1-2\varepsilon),\  x+\varepsilon y=1,\\   x+ 2\varepsilon y=1+2\varepsilon \ea $ \\[0mm]
 \hline
\!\!\!\! (\!{\it V.24${}_\varepsilon$})\,: $\ba{l} \dot x=x^2+1,\
\dot y= 1-\varepsilon^2 y^2\ea$
& $ \ba{l} x=\pm i,\  \varepsilon y=\pm 1 \ea $ \\[0mm]
 \hline
\!\!\!\! (\!{\it V.25${}_\varepsilon$})\,: $\ba{l} \dot x=-1-x^2,\
\dot y= 1-2xy-\varepsilon y^2 \ea$
& $ \ba{l} x=\pm\, i,\  x+ \varepsilon y=\pm\, i\sqrt{1-\varepsilon} \ea $ \\[0mm]
 \hline
\!\!\!\! (\!{\it V.26${}_\varepsilon$})\,: $\Bigg\{\ba{l} \dot
x=g+(2g\varepsilon-1)x-2\varepsilon x^2,\\[-1.1mm]
\dot y= (2g\varepsilon+1)y-x^2-6g\varepsilon^2xy\\[-1.1mm]
\qquad +3\varepsilon^2(1+2g\varepsilon-3g^2\varepsilon^2) y^2\!\!
\ea$
& $ \ba{l} x=g,\  3\varepsilon[x+\varepsilon(3g\varepsilon+1)y]=-1,\\
  2\varepsilon x=-1,\ \varepsilon[x+3\varepsilon(g\varepsilon-1)y]=1 \ea $ \\[0mm]
 \hline
\!\!\!\! (\!{\it V.27${}_\varepsilon$})\,: $\Big\{\ba{l} \dot
x=1+ x+ \varepsilon x^2,\\[-1.1mm]
\dot y= y-x^2 -2\varepsilon xy-2\varepsilon^2y^2\!\! \ea$
& $ \ba{l}1+x+ \varepsilon^2x^2=0,\\[-1.1mm] \varepsilon(x+2\varepsilon y)^2-(x+2\varepsilon y)=1 \ea $ \\[0mm]
 \hline
\!\!\!\! (\!{\it V.28${}_\varepsilon$})\,: $\Bigg\{\ba{l} \dot
x=(\varepsilon-1)\varepsilon^2+2\varepsilon^3 x\\[-1.1mm]
\qquad+(1-\varepsilon)(1-2\varepsilon+3\varepsilon^2)x^2,\\[-1.1mm]
\dot y= (1-\varepsilon)(2\varepsilon^2y+1)(x+2\varepsilon y+1)
\!\!
\ea$
& $ \ba{l}(\varepsilon-1)x=\varepsilon,\ 2\varepsilon^2y=-1,\\[-1.1mm]
 (1-2\varepsilon+3\varepsilon^2)x=\varepsilon(1-\varepsilon),\\[-1.1mm]
  (\varepsilon -1)^2x-4\varepsilon^3 y=\varepsilon(\varepsilon+1) \ea $ \\[0mm]
 \hline
\!\!\!\! (\!{\it V.29${}_\varepsilon$})\,: $\ba{l} \dot
x=\varepsilon^2-x^2,\ \dot y= 1-2xy-\varepsilon y^2 \ea$
& $ \ba{l} x=\pm\, \varepsilon,\  x+ \varepsilon y=\pm\sqrt{\varepsilon^2+\varepsilon} \ea $ \\[0mm]
 \hline
\!\!\!\! (\!{\it V.30${}_\varepsilon$})\,: $\Big\{\ba{l} \dot
x=1+ \varepsilon x+ \varepsilon^3 x^2,\\[-1.1mm]
\dot y= g+ \varepsilon y-x^2 -2\varepsilon^3
xy-2\varepsilon^6y^2\!\!
\ea$
& $ \ba{l}1+\varepsilon x+ \varepsilon^3x^2=0,\\[-1.1mm]
\varepsilon^3(x\!+\!2\varepsilon^3 y)^2\!-\!\varepsilon(x\!+\!2\varepsilon^3 y)\!=\!1\!+\!2g\varepsilon^3\!\!\! \ea $ \\[0mm]
\hline
\end{tabular}
\end{center}
\end{table}

\textit{Proof of Theorem \ref{th_mil_5}:}  Since  we only discuss
the case $C_2\ne0$, in what follows it suffices to consider only
the canonical forms $(\SSS_I)$ to $(\SSS_{I\,V})$.  The idea of
the proof is the same as in the proof of the Theorem
\ref{th_mil_6}. We shall  perform a case by case discussion for
each one of these canonical forms, for which according   to
Lemma~\ref{lm_NB2}  we must examine two subcases: $(i)\ N=B_2=0$
and  $ (ii)\ N\ne0$, $\theta=B_3=0$. Each one of these conditions
yields  specific conditions on the parameters. The discussion
proceeds further by breaking these cases in more subcases
determined by more restrictions on the parameters. Finally we
construct  invariants or T-comitants which put these conditions in
invariant form.

 %%%%%%%%  New subsection %%%%%%%
%%%%%%

\vspace{-3mm}
\subsection{Systems with the divisor\ $D_S(C,Z)=1\cdot w_1+1\cdot w_2+1\cdot w_3$}
\vspace{-3mm}
For this case we shall later  need the following polynomials which
are shown to be $T$-comitants in  Lemma \ref{Table:Propreties}.
\bnot\label{not_H4,5} Let us denote
$$
\bal
&H_4(a)=\big((C_2,D)^{(2)},(C_2,D_2)^{(1)}\big)^{(2)},\
  H_6(a,x,y)= 16N^2(C_2,D)^{(2)}+ H_2^2(C_2,C_2)^{(2)},\\
  &H_5(a)=\big((C_2,C_2)^{(2)},(D,D)^{(2)}\big)^{(2)}+
    8\big((C_2,D)^{(2)},(D,D_2)^{(1)}\big)^{(2)}.
\eal
$$
\enot
\vspace{-5mm}
\subsubsection{The case $N=0=B_2$}
\vspace{-3mm} It was shown above (see page~\pageref{s4.1}) that
the systems $(S_I)$  with $N(\ab,x,y)=0$ can be brought by an
affine transformation to the systems~(\ref{s4.1}) for which we
have
\beq
&&B_2=648\big[e^2(4k-4l-e^2)x^4+2d^2e^2(2x^2-3xy+2y^2) -d^2(4k-4l+d^2)y^4\big].
\eeq
Hence the condition $B_2=0$ yields
$de=e(4k-4l-e^2)=d(4k-4l+d^2)=0$. According to Lemma~\ref{gcd:4},
in order to have $M_{{}_\IL}=5$ we must satisfy the condition $\
\deg\gcd({\cal E}_1, {\cal E}_2)=4. $\  We claim, that for this it
is necessary that $d=e=0$. Indeed, let us suppose, that $de=0$ but
$d^2+e^2\ne0$.  Then by interchanging $x$ and $y$ we may assume
$d=0,\ e=2$ via   Remark \ref{rem:transf} ($\gamma=e/2,\ s=1$).
Then we obtain the systems
\be \label{ss2_1}
  \dot x=k +  x^2,\quad  \dot y=l +2x + y^2,
\ee
for which the condition $B_2=2^73^4(k-l-1)x^4=0$   yields $k=l+1$.
Then for  the systems~(\ref{ss2_1}) with $k=l+1$  we obtain
$$
\bal
&{\cal E}_1= -2\big[Y^2 - YZ + Z( X + Z + l Z)\big]{\cal H},\quad
{\cal E}_2=-(X + Y - Z)\big[Y^2 + Z(2 X + l Z)\big]{\cal H},\\
\eal
$$
where $ {\cal H}= (Y-X+Z)(X^2+Z^2+lZ^2)$. Thus, $\deg{\cal H}=3$
and we shall  show that for all values given to the parameter  $l$
the degree of $\gcd({\cal E}_1,{\cal E}_2)$ remains three. Indeed,
 since the common factor of the polynomials ${\cal E}_1/{\cal H}$
and ${\cal E}_2/{\cal H}$ must depend on $Y$, according to Lemma
\ref{Trudi:2} it is sufficient to observe that $ \Res_Y({\cal
E}_1/{\cal H},\ {\cal E}_2/{\cal H})=-8 Z^2[X^2 + (1 +
l)Z^2]^2\ne0$. This   proves our claim and hence, the condition
$d=e=0$  must hold. Since for  systems (\ref{s4.1}) we have
$H_4=96(d^2+e^2)$ this condition is equivalent to $H_4=0$.

 Assuming $H_4=0$  (i.e. $d=e=0$) the systems (\ref{s4.1}) become
\be \label{2CF_1}
  \dot x=k +  x^2,\quad  \dot y=l + y^2,
\ee
and calculations yield $\quad
 {\cal E}_1=  2(X - Y)\,{\cal H},\quad
{\cal E}_2= \big[X^2 - Y^2 + (k - l)Z^2\big]{\cal H},$ where\quad
${\cal H}= (X^2 + k Z)(Y^2 + l Z^2). $
 Hence by Theorem \ref{theor:E1,E2}  each  system in the family \eqref{2CF_1} possesses
 four invariant affine lines which means that for these systems
 $M_{{}_\IL}\ge5$. We observe that to have an additional common
factor of ${\cal E}_1$ and ${\cal E}_2$ it is necessary and
sufficient that $k-l=0$. So, to have $M_{{}_\IL}=5$ the condition
$k-l\ne0$  must be satisfied. This condition  is equivalent to
$B_3\ne0$, since for the systems~(\ref{2CF_1}) we have $
B_3=4(l-k)x^2y^2$.

Systems (\ref{2CF_1}) possess the invariant  lines, components of
the conics:\   $x^2+k=0,$\ $ y^2+l=0, $\  and then we obtain the
following configurations of invariant straight lines ({\bf Diagram
2}):\\
 $(i)$\  Config. 5.3  for
$k<0$ and $l<0$;\ \ $(ii)$\  Config. 5.4  for $kl<0$;\ \ $(iii)$\
Config. 5.5 for $k>0$ and $l>0$; \ \ $(iv)$\ Config.~5.12 for
$kl=0$ and $k+l<0$; \ \ $(v)$\  Config.~5.16 for $kl=0$ and
$k+l>0$.

On the other hand for the systems (\ref{2CF_1}) we have:\quad $
 H_1= -2^73^2\,(k+l), \ \  H_5= 2^{11}3\,kl.
$\quad Herein we conclude, that these two $T$-comitants capture in
invariant form  exactly the conditions for  distinguishing the
Configurations 5.3--5.5, 5.12 and 5.16 as it is indicated in Table
4.

We observe, that if for the systems~(\ref{2CF_1}) the condition $H_5\le0$ holds
(i.e. $kl\le0$) then interchanging $x$ and $y$ we may assume $l\ge0$. Moreover,
by Remark \ref{rem:transf}  ( $\gamma=|k|,\ s=1/2$) we may assume
 $k\in\{-1,1\}$. We also note that for $l<0$ (respectively, $l>0$) we may set $l=-g^2$
(respectively, $l=g^2$) and due to the substitution $y\to gy$ we
obtain the canonical system ({\sl V.3}) (respectively, ({\sl V.4})
and ({\sl V.5}) ) from Table 4.
%%%%%%%%  New Subsubsection %%%%%%%
%%%%%

\vspace{-3mm}
\subsubsection{The case $N\ne0$, $\theta=0=B_3$}
\vspace{-3mm} For the canonical systems $(\SSS_I)$ we calculate\ \
$
 \theta=  -8(h-1)(g-1)(g+h).
$\ \ Hence the condition $\theta=0$ yields $(h-1)(g-1)(g+h)=0$
and without loss of generality we can consider $h=1$. Indeed, if
$g=1$ (respectively, $g+h=0$) we can apply the linear
transformation which will replace the straight line $x=0$ with $
y=0$ (respectively, $x=0 $ with $y=x$) reducing this case to
$h=1$.
 Assume $h=1$. Then
$N=(g^2-1)x^2\ne0$ and we may assume $e=f=0$  via a translation.\
 Thus the systems $(\SSS_I)$  become
\bq\label{s1_kap}
  \dot x=k + cx+dy +g x^2,\quad
  \dot y=l  + (g-1)xy + y^2
\eq
and calculations yield $ \mu =\ 32g^2$ and $B_3 =
-3l(g-1)^2x^4+6l(g-1)^2x^3y-6d^2gxy^3+3d^2gy^4\
 +3\left[(4gl-k(g+1)^2+c^2+cd-cdg\right]x^2y^2.$\quad
The condition  $B_3$=0 implies $dg=0$. We shall examine two
subcases: $\mu\ne0$ and $\mu=0$.
%%%%%%%%%%%%%%%%%%%%%%%%%%%%%%%%

\smallskip\noindent
{\bf The subcase $\mu\ne0$.} In this case we obtain
$g\ne0$, and from $g-1\ne0$ the condition $B_3=0$  for the
systems~(\ref{s1_kap}) yields  $d= l= c^2-k(g+1)^2=0$. Since
$N\ne0$ then $g+1\ne0$ and we may set $c=u(g+1) $ where $u$ is a
new parameter. Then $k=u^2$ and we obtain the
systems\\[-6mm]
\be\label{s1_kap_}
\dot x=u^2 + u(g+1)x + gx^2,  \quad \dot y=(g-1)xy +y^2,
\ee
for which $H_1=576u^2(g-1)^2$.

{\bf 1)} If $ H_1\ne0$ then $u\ne0$ and we may assume $u=1$ via
Remark~\ref{rem:transf} ($\gamma=u$, $s=1$). This leads to the
systems
\be\label{2CF_18}
  \dot x=(x+1)(gx+1),  \quad \dot y=(g-1)xy +y^2,
\ee
for which $g(g^2-1)\ne0$ and calculations yield:
\be\label{GCD:1}
 {\cal H}=\gcd\left({\cal E}_1,{\cal E}_2\right)=  Y (Y-X - Z) (X + Z) (g X + Z).
\ee
Hence, $\deg {\cal H}=4$. By hypothesis $N\ne0$ and hence,
according to Lemma \ref{gcd:4} for every $ g$ such that
$g(g^2-1)\ne0$, $M_{{}_\IL}\le5$.
 By Theorem~\ref{theor:E1,E2}, from (\ref{GCD:1})  the systems (\ref{2CF_18}) possess the
following four distinct invariant affine lines:\ \ $
   y=0,\ \ x+1=0,\ \ x-y+1=0,\ \ gx+1=0.
$\ \ Thus we obtain the Config. 5.1.
%%%%%%%%%%%%%%%%

 {\bf 2)} For $ H_1=0$ we have $u=0$ and
the systems~(\ref{s1_kap_}) become
\be\label{2CF_20}
\dot  x=  gx^2,  \quad  \dot y=(g-1)xy +y^2,
\ee
with $g(g^2-1)\ne0$ and  we  calculate:\ \ $
 {\cal H}=\gcd\left({\cal E}_1,{\cal E}_2\right)= g X^2Y(X - Y).$
\ \ Hence $\deg {\cal H}=4$ and we obtain $M_{{}_\IL}\ge5$. Since
$N\ne0$ by Lemma \ref{gcd:4}, $M_{{}_\IL}$ cannot be equal to 6.
The systems~(\ref{2CF_20}) possess the invariant lines $x=0,$
$y=0$ and $x=y$. Moreover, according to Lemma \ref{lm3} the line
$x=0$ could  be of multiplicity two and the perturbations ({\it
V.8${}_\varepsilon$}) from Table 5  show this. Hence,  for $H_1=0$
we obtain Config. 5.8.
%%%%%%%%%%%%%%%%%%%%%%%%%%%%%%%

\smallskip\noindent
{\bf The subcase $\mu=0$.} The condition $\mu =32g^2=0$ yields
$g=0$, and for the systems~(\ref{s1_kap}) the condition $B_3=0$
yields $g=l=c(c+d)-k=0$. Thus, $g=l=0$,  $k=c(c+d)\ne0$, otherwise
we get degenerate systems (\ref{s1_kap}). Hence,  we may assume
$c=1$ via Remark~\ref{rem:transf} ($\gamma=c$, $s=1$) and  we
obtain the systems
\be\label{S_mu0}
  \dot x=d+1 +x+dy,  \quad \dot y=-xy +y^2.
\ee
Calculations yield:
$$
\bal
& {\cal E}_1= \left[-X^2 + 2 X Y + d(Y + Z)^2 + Z(2 Y +
Z)\right]{\cal H},\qquad {\cal H}=  Y Z (X - Y + Z + d Z),\\
&{\cal E}_2= (Y-X) (Y + Z) \big[X + Z + d (Y + Z)\big]{\cal
H},\quad \Res_X({\cal E}_1/{\cal H},\ {\cal E}_2/{\cal H})=-9d(d +
1)^2(Y + Z)^6.
\eal
$$
Hence $\deg\,{\cal H}=\gcd\left({\cal E}_1,{\cal E}_2\right)=3$
and the condition on the parameter $d$  so as to have an
additional common factor of ${\cal E}_1$, ${\cal E}_2$, according
to Lemma \ref{Trudi:2} is $\Res_X({\cal E}_1/{\cal H},\ {\cal
E}_2/{\cal H})\equiv0$.  Since $d+1\ne0$ (otherwise we get the
degenerate system \eqref{S_mu0}) this  condition   yields $d=0$.
Then we obtain the following system
\be\label{2CF_22}
  \dot x=1 +x,  \qquad \dot y=-xy +y^2
\ee
for which ${\cal H}=\gcd\left({\cal E}_1,{\cal
E}_2\right)=YZ(X+1)(X-Y+Z)$. We observe that  this system
possesses the invariant affine straight lines:\ \ $
   y=0,\ \ x+1=0,\ \ x-y+1=0.$\ \
Taking into account that $Z\mid{\cal H}$, we have  by Corollary
\ref{Mult:Z=0} that the line $Z=0$ could be of multiplicity two.
This is confirmed by the  perturbations ({\it
V.7${}_\varepsilon$}) from Table 5.
   On the other hand for the systems~(\ref{S_mu0}) calculations yields
$H_6=128dx^2(x^2-xy+y^2)(x^2-2xy-dy^2)$. Hence, the conditions
$g=0=d$ are equivalent to   $\mu=0$ and $H_6=0$. In this case we
obtain Config. 5.7.
%%%%%%%%%%%%%%%%%%%%%%%%%%%%%%%%%%%%
%%%%%%%%  New Subsection %%%%%%%%%%%%
%%%%%%%%%%%%%%%%%%%%%
\vspace{-3mm}
\subsection{Systems with the divisor\ $D_S(C,Z)=1\cdot w_1^c+1\cdot w_2^c+1\cdot w_3$}
\vspace{-3mm}
 We are in the case of the canonical form
$(\SSS_{I\!I})$.
%%%%%%%%  New Subsubsection %%%%%%%%%%%%
\vspace{-3mm}
\subsubsection{The case $N=0=B_2$}
\vspace{-3mm} It was shown above (see page~\pageref{S2:N=0}) that
the systems $(\SSS_{I\!I})$ with $N(\ab,x,y)=0$ can be brought by
an affine transformation to the systems~(\ref{S2:N=0}) for which
we have
\beq
&& B_2=648\,\big[(8efk-(f^2+e^2)^2\big]x^4-16k(e^2-f^2)xy(x^2-y^2)
-48efkx^2y^2  +8efky^4,\\
&& B_3=6\,\big[(ef-2k)x^4+(f^2-e^2)x^3y-(4k+ef)x^2y^2-2ky^4\big].
\eeq
If $B_3=0$ then $k=e=f=0$ and we obtain the systems~(\ref{CF_3})
for which $M_{{}_\IL}=6$ (see page~\pageref{CF_3}). Hence
$B_3\ne0$ and this implies $k\ne0$, otherwise  from $B_2=0$ we
obtain again $e=f=0$. Therefore  $k\ne0$ and we may consider $k>0$
via the change $x\to -x$ and  by  Remark \ref{rem:transf}
($\gamma=k,\ s=1/2$) we may assume $k=1$. Then the condition
$B_2=0$ yields $e=f=0$ and  we obtain the systems
\be\label{2CF_4}
  \dot x= 1 +2xy,\qquad
  \dot y= l  -x^2 + y^2
\ee
which  possess the invariant  lines\quad $ y+ix=\pm
\sqrt{-l-i}$,\quad $y-ix=\pm\sqrt{i-l}.$\quad This leads to the
Config. 5.6.
%%%%%%%%  New Subsubsection %%%%%%%%%%%%
\vspace{-3mm}
\subsubsection{The case $N\ne0,$ $\theta=B_3=0$}
\vspace{-3mm}
For the  systems $(\SSS_{I\!I})$ we calculate\\[-4mm]
$$
 \theta= 8(h+1)[(h-1)^2+g^2],\quad N=(g^2-2h+2)x^2 +2g(h+1)xy+(h^2-1)y^2
$$
and hence  by $N\ne0$, the condition $\theta=0$ yields $h=-1$.
Then we may assume $f=0$ due to a translation and   the systems
$(\SSS_{I\!I})$ become
\be\label{S_IIa}
\dot x=k+cx+dy+gx^2,\quad
\dot y= l+ex-x^2+gxy-y^2.
\ee
For these systems  calculations yield\  Coefficient$[B_3,\
y^4]=-3d^2g $\ and $ \mu= 32\,g^2$. We shall   examine two
subcases: $\mu\ne0$ and $\mu=0$.
%%%%%%%%%%%%%%%%%%%%%

\smallskip\noindent
{\bf The subcase $\mu\ne0$}. This yields  $g\ne0$ and then the
condition $B_3=0$ implies $d=0$. Moreover, we may assume $c=0$ via
the translation of the origin of coordinates to the point
$(-c/(2g),-c/4)$.  Thus, the systems (\ref{S_IIa}) become
$$
  \dot x= k +2gx^2,\quad
  \dot y= l +ex -x^2 +2gxy -y^2,
$$
for which we calculate $\ B_3 =
3\left[k(4-g^2)-4gl\right]x^2(x^2-y^2)+
      6\left[l(4-g^2)+4gk+e^2\right]x^3y.
$ Hence, the condition $B_3=0$ yields the following linear system
of equations with
respect to parameters $k$ and $l$:\\[-5mm]
$$
    k(4-g^2)-4gl=0, \quad 4gk+l(4-g^2)+e^2=0.
$$
Setting $e=u(g^2+4)$ ( $u$ is a new parameter)   we have the
following solution of this  system:  $k=-4gu^2$, $l=(g^2-4)u^2$.
Thus we obtain the  systems
\be\label{S_II_kap0}
  \dot x= -4gu^2 + gx^2,\ \,
  \dot y= (g^2-4)u^2 +u(g^2+4)x -x^2 +g xy -y^2,
\ee
for which $H_1= -2^{12}3^2u^2g^2$.
%%%%%%%%%%%%%%%%%%%%%%%%

{\bf 1)} If $H_1\ne0$ we have $u\ne0$ and we can assume $u=1$ via
the Remark~\ref{rem:transf} ($\gamma=u,$ $s=1$). Hence the
systems~(\ref{S_II_kap0}) become
\be\label{2CF_19}
  \dot x= g(x^2-4),\ \,
  \dot y= (g^2-4) +(g^2+4)x -x^2 +gxy -y^2.
\ee
and calculations yield:$\quad  {\cal H}= g(X - 2 Z)(X + 2
Z)(X^2+Y^2-4 X Z +  2\,g Y Z + 4Z^2 + g^2 Z^2). $\ Hence,
$M_{{}_\IL}\ge5$ and since $N\ne0$ by  Lemma \ref{gcd:4}
$M_{{}_\IL}$ cannot be equal to 6.
 By Theorem \ref{theor:E1,E2} the systems (\ref{2CF_19})
possess the following four distinct invariant straight lines:\ \
$y-ix+g+2i=0,$ \ \ $ y+ix-2i+g=0, $\ \ $ x=\pm2.$ \ Thus we obtain
the Config. 5.2.
%%%%%%%%%%%%%%%%%%%%%%%%

{\bf 2)} For $H_1=0$ we have $u=0$ and the
systems~(\ref{S_II_kap0}) become
\be\label{2CF_21}
  \dot x= gx^2,\quad
  \dot y= -x^2 +gxy -y^2,
\ee
with $g\ne0$. We calculate \ $  {\cal H}=\gcd\left({\cal
E}_1,{\cal E}_2\right)=  g X^2(X^2 + Y^2) $ \ and hence $\deg
{\cal H}=4$. Since  $N\ne0$ by Lemma \ref{gcd:4}  we obtain that
$M_{{}_\IL}$ equals exactly 5.
 The systems (\ref{2CF_20}) possess the following invariant
straight lines:\ \ $    x=0,\ \  y=\pm ix$ and the line $x=0$
could be of multiplicity two. This is confirmed by the
perturbations ({\it V.10${}_\varepsilon$}) from Table~5. Thus we
get   Config. 5.10.
%%%%%%%%%%%%%%%%%%%%%

\smallskip\noindent
{\bf The subcase $\mu=0$.} Then we obtain $g=0$ and we may assume
$e=0$ via a translation. Therefore the systems (\ref{S_IIa})
become $\quad
  \dot x= k+cx+dy,\quad  \dot y= l -x^2 -y^2
\quad $ and we calculate\quad $ B_3= 12kx^4 +
6(4l-c^2-d^2)x^3y-12kx^2y^2. $\ Hence the condition $B_3=0$ yields
$k=4l-c^2-d^2=0$. We replace $c$ by $2c$ and $d$ by $2d$ and then
we obtain $l=c^2+d^2$. This leads to  the systems:
\be\label{S_II_2}
  \dot x= 2cx+2dy,\quad  \dot y= c^2+d^2 -x^2 -y^2
\ee
for which  calculations yield:
$$
\bal
 &{\cal E}_1= \left[d X^2 - 2\,c XY - d Y^2 + 2\,
              (c^2+d^2)XZ  +d(c^2+d^2) Z^2\right]{\cal H},\\
 &{\cal E}_2=(c X + d Y)[X^2 + Y^2 + 2\,d X Z -
        2\,cY Z + (c^2 + d^2) Z^2]\,{\cal H},\\
& {\cal H}=2 Z[X^2 + Y^2 - 2\,d X Z + 2\,c Y Z + (c^2 + d^2) Z^2].
\eal
$$
Thus,  $\deg{\cal H}=3$ and we need an additional factor of ${\cal
E}_1$ and ${\cal E}_2$. Since $c^2+d^2\ne0$ we observe that such
common factor of the polynomials ${\cal E}_1/{\cal H}$ and ${\cal
E}_2/{\cal H}$ must depend on $X$. Hence,  by Lemma \ref{Trudi:2}
the following condition must hold:
$$
\bal
\Res_X({\cal E}_1/{\cal H},\ {\cal E}_2/{\cal H})&=4\,d(c^2 +
d^2)^2(Y - c Z)^6=0.
\eal
$$
Therefore the condition  $d=0$ must be satisfied and then $c\ne0$
(otherwise we get the degenerate system \eqref{S_II_2}).\ On the
other hand for the systems~(\ref{S_II_2}) we have
 $H_6=-2^{13}dx^3(3x^2-y^2)(dx^2-2cxy-dy^2)$ and hence the condition $d=0$ is
 are equivalent to $H_6=0$. We may assume $c=1$  via the
Remark~\ref{rem:transf} ($\gamma=c$, $s=1$) and then we obtain the system
\be\label{2CF_23}
  \dot x= 2x,\quad  \dot y= 1 -x^2 -y^2.
\ee
For this system we  calculate $\quad {\cal H}=\gcd\left({\cal
E}_1,{\cal E}_2\right)=4\, XZ(X^2 + Y^2 + 2Y Z + Z^2) $\ and
according to Theorem \ref{theor:E1,E2} the system (\ref{2CF_23})
possess the invariant affine lines: $ x=0$ and $y\pm ix +1=0 $.\
Moreover, by Corollary \ref{Mult:Z=0}  the line $l_\infty: Z=0$
could be of multiplicity two.
 This is confirmed by the   perturbations ({\it
 V.9${}_\varepsilon$}) from Table~5.
Therefore we obtain the  Config. 5.9.
%%%%%%%%%%%%%%%%%%%%%%%%%%%%%%%%%%%%%%%%
%%%%%%%%  New Subsection %%%%%%%%%%%%
%%%%%%%%%%%%%%%%%%%%%%%%%%%%%%%%%%%%%%%%%
\vspace{-3mm}
\subsection{Systems with the divisor\ $D_S(C,Z)=2\cdot w_1+1\cdot w_2$}
\vspace{-3mm}
 We are in the case of the canonical form
$(\SSS_{I\!I\!I})$. For this case we shall later  need the
following polynomial which is shown to be a $CT$-comitant in Lemma
\ref{Table:Propreties}.
 \bnot\label{not_N5} Let us denote\quad $N_5(\ab,x,y)=
\big((D_2,C_1)^{(1)} + D_1D_2\big)^2
-4\big(C_2,C_2\big)^{(2)}\big(C_0,D_2\big)^{(1)}$.
\enot
%%%%%%%%  New Subsubsection %%%%%%%%%%%%
\vspace{-5mm}
\subsubsection{The case $N=0=B_2$}
\vspace{-3mm} It was previously shown  (see page~\pageref{eq_H}) that to
examine the systems $(\SSS_{I\!I\!I})$  with $N(\ab,x,y)=0$  we
have to consider  two subcases: $H(\ab,x,y)\ne0$ and
$H(\ab,x,y)=0$.
%%%%%%%%%%%%%%%%%%%%%%

\smallskip\noindent
{\bf The subcase $H(\ab,x,y)\ne0.$} In this case the systems
$(\SSS_{I\!I\!I})$ with  $N=0$  can be brought by an affine
transformation to the systems (\ref{S3_NM_}) (see  page
\pageref{S3_NM_})  for which we have: $\ B_2= -648d(8clx^4+
16dlx^3y+ d^3 y^4).$\  Therefore, the condition $B_2=0$ yields
$d=0$ and we obtain the systems
\be\label{S3_B2}
  \dot x =k +cx -x^2,\qquad
  \dot y =l - 2xy,
\ee
for which calculations yield:  $\
 {\cal E}_1=\left[-2 X^2 Y + Z^2(2 k Y + c l Z)\right] {\cal H},
 $\quad
${\cal H}=(kZ^2+cXZ-X^2),$\quad
  $
 {\cal E}_2= (X^2 - c X Z - k Z^2)(2 X Y - l Z^2) {\cal H}.$
\ Thus, $\deg{\cal H}=2$ and to have $M_{{}_\IL}=5$  the
polynomials ${\cal E}_1/{\cal H}$ and ${\cal E}_2/{\cal H}$ must
have a common factor of degree two. We observe, that this common
factor necessarily  depends on $X$ and hence  by Lemma
\ref{Trudi:2} the following condition must hold:
$$
\Res_X({\cal E}_1/{\cal H},\ {\cal E}_2/{\cal H})=-2 c^2 Y Z^6(4k
Y^2 + 2cl YZ - l^2 Z^2)^2\equiv0.
$$
Herein we obtain either $c=0$ or $k=l=0$, but the second case
yields degenerate systems. So, we assume $c=0$ and then for the
systems (\ref{S3_B2}) we obtain $\quad {\cal E}_1= 2Y\,{\cal H},$
$\  {\cal H}= (kZ^2-X^2)^2,$\ $ {\cal E}_2= (-6 X Y + 3 l
Z^2)\,{\cal H}.$ \ We observe, that $\deg\,{\cal H}=4$ and that
the polynomials ${\cal E}_1$ and ${\cal E}_2$ do not have an
additional common factor  if and only if $l\ne0$. This condition
is equivalent to $B_3\ne0$,  since for the systems (\ref{S3_B2})
we have $B_3=-12lx^4$. By $l\ne0$ we may consider $l=1$ via the
rescaling $y\to ly$ and we obtain the systems
\be\label{2CF_7}
  \dot x =k  -x^2,\qquad
  \dot y =1 - 2xy.
\ee
Moreover, due to the re-scaling  $x\to|k|^{1/2}x$,
$y\to|k|^{-1/2}y$ and $t\to|k|^{-1/2}t$ (for $k\ne0$) we may
assume $k\in \{-1,0,1\}$. These systems possess  two invariant
lines $x=\pm\sqrt{k}$. By Lemma~\ref{lm3} for $k\ne0$ each one of
these lines could be of multiplicity two and for $k=0$ the
invariant line $x=0$ of the system (\ref{2CF_7}) is of the
multiplicity four. This is confirmed by the perturbations ({\it
V.22${}_\varepsilon$}) (respectively, ({\it
V.25${}_\varepsilon$}); ({\it V.29${}_\varepsilon$})) from Table~5
for $k=1$ (respectively, $k=-1$; $k=0$). Thus, we obtain Config.
5.22 (respectively, Config. 5.25; Config. 5.29).

On the other hand  for the systems~(\ref{S3_B2}) calculations
yield: $H_2= 16cx^2$ and $H_3= 32kx^2$. Hence, these $T$-comitants
capture exactly the conditions   $c=0$ and $k<0$ (respectively
$c=0$, $k=0$ or $c=0$, $k>0$) . It remains to observe, that the
condition $B_3\ne0$ implies $H\ne0$, since for $H=0$ the condition
$B_2=0$ implies $B_3=0$ (see the subcase $H(a,x,y)=0$ below).
%%%%%%%%%%%

\smallskip\noindent
{\bf The subcase $H(a,x,y)=0$.} It is previously shown  (see page
\pageref{S3_NH_0}) that if $N=H=0$ then  the systems
$(\SSS_{I\!I\!I})$    can be   brought by an affine transformation
to the systems (\ref{S3_NH_0}). For these systems we have $\  B_2=
-648d^4y^4,$ \ $ B_3= 6dxy^2(fx-dy) $\ and hence the condition
$B_2=0$ yields $d=0$. Therefore  the conditions $B_2=0$ and
$B_3=0$ are equivalent and since for any quadratic system
\eqref{2l1} the condition $B_3=0$ implies $B_2=0$ (see the formulas
\eqref{Comit:Bi} on page \pageref{Comit:Bi}),  we shall use in this case
the condition $B_3=0$.

 Assuming $d=0$ we obtain the systems (\ref{s4.5}) for which $D(x,y)=-f^2x^2y$
and we shall  consider two subcases: $D\ne0$ and $D=0$.
%%%%%%

{\bf 1)} For $D\not=0$ the systems (\ref{s4.5}) can  be brought by
an affine transformation to the systems (\ref{s4.6}) and
calculations yield  the values (\ref{val:Ei}) of the affine
comitants ${\cal E}_1$ and ${\cal E}_2$. We observe that
$\deg\,{\cal H}=3$ and taking into account that the polynomials
${\cal E}_1/{\cal H}$ and ${\cal E}_2/{\cal H}$ cannot have the
common factor $Z$, to have  an additional factor of these
polynomials  according to Lemma \ref{Trudi:2} at least one of  the
following two conditions must hold: $\ \Res_X({\cal E}_1/{\cal
H},\, {\cal E}_2/{\cal H})= -36 e(k+1)(4 Y^2 Z + e^2k Z^3)^2=0,\ $
$ \Res_Y({\cal E}_1/{\cal H},\, {\cal E}_1/{\cal H})=-6e(X^2 +
kZ^2)^2=0. $
 Thus we obtain either the condition $e=0$ or $k=-1$.\
On the other hand for  systems (\ref{s4.6})   we obtain
 $N_1=8ex^4$ and $N_2=16(k+1)x$ and we shall consider two
 subcases: $N_1=0$ and  $N_1\ne0,$\ $N_2=0$.

{\bf 1a)} Assume $N_1=0$. Then $e=0$  and
 the systems (\ref{s4.6})  become
\be\label{2CF_10}
\dot x=k+ x^2,\qquad \dot y= 2y.
\ee
Calculations yield: $ \
 {\cal E}_1=2(X-Z)\,{\cal H},\ $ $
{\cal E}_2= 3(X^2 + k Z^2)\,{\cal H},\ $ $ {\cal H}=4YZ(kZ^2+X^2),
$ $  \Res_X({\cal E}_1/{\cal H},\ {\cal E}_2/{\cal
H})=12(k+1)Z^2.$\ Hence  $\deg{\cal H}=4$ and we observe that in
order not to have an additional common factor of the polynomials
${\cal E}_1$ and ${\cal E}_2$ we must have   $k+1\ne0$ (i.e.
$N_2\ne0$). The systems (\ref{2CF_10}) possess the invariant
affine lines $y=0,\ x=\pm\sqrt{-k} $. According to Corollary
\ref{Mult:Z=0} the line $l_\infty: Z=0$ could be  of multiplicity
2 and the line $x=0$ also could be  of multiplicity 2 in the case
when $k=0$. Since for  systems (\ref{s4.6})   we have
$N_5=-64kx^2$, we obtain  Config. 5.11 for $N_5>0$, Config. 5.15
for $N_5<0$ and  Config. 5.17 for $N_5=0$.

Note that for $k<0$ (respectively, $k>0$) one can set $k=-g^2$
(respectively, $k=g^2$) and due to the substitution $x\to gx$ we
obtain the canonical system ({\sl V.13}) (respectively, ({\sl
V.15})) from Table 4.  It remains to observe that the
perturbations
  ({\it V.13${}_\varepsilon$}), ({\it V.15${}_\varepsilon$})
and ({\it V.17${}_\varepsilon$}) from Table 5  confirm the
validity of the Config. 5.13, 5.15 and 5.17, respectively.

{\bf 1b)} For $N_1\ne0$, $N_2=0$ we have $e\ne0$, $k=-1$  and then
for the  systems~(\ref{s4.6}) calculations yield: $\ {\cal E}_1=
\left[4 Y + e (X - Z)\right]{\cal H}, $ $
  {\cal E}_2=  (ce X + 2 Y) (X + Z)\,{\cal H}, $ $
 {\cal H}=  Z(X+Z)(X-Z)^2, $ $
 \Res_Y({\cal E}_1/{\cal H},\ {\cal E}_2/{\cal H})=-2 e(X + Z)^2.
$\ Hence  $\deg{\cal H}=\deg\gcd\left({\cal E}_1,{\cal
E}_2\right)= 4$ and  since $N_1\ne0$ (i.e. $e\ne0$) the
polynomials  ${\cal E}_1$ and ${\cal E}_2$ could not have an
additional common factor. Assuming   $e=1$ via the rescaling $y\to
ey$   the systems (\ref{s4.6}) become
\be\label{2CF_13}
\dot x=-1+ x^2,\qquad
\dot y= x+2y.
\ee
This  system  possesses the invariant lines $x=\pm1$. According to
Lemma \ref{lm3} and Corollary \ref{Mult:Z=0} the line $x=1$ as
well as the line $Z=0$ could  be  of multiplicity two. This is
confirmed by the perturbations ({\it V.21${}_\varepsilon$}) from
Table~5. Thus for $N_1\ne0$ and $N_2=0$   we obtain Config.~5.21.
%%%%%

 {\bf 2)} If  $D=0$  then we have $f=0$ and the systems
(\ref{s4.5}) become  the systems (\ref{s4.7}) (see page
\pageref{s4.7}) for which calculations yield the corresponding
expressions (\ref{val:Eia}) for the affine comitants ${\cal E}_1$
and ${\cal E}_2$. As  $\deg\,{\cal H}=3$  we need an additional
common factor of ${\cal E}_1$ and ${\cal E}_2$. Taking into
account that these polynomials   depend only on $X$ and $Z$,
according to  Lemma \ref{Trudi:2}  at least one of the following
two conditions must hold:
$$
\Res_X({\cal E}_1/{\cal H},\ {\cal E}_2/{\cal H})=-4\,ek(e^2 k +
l^2)^2 Z^6=0,\ \,\Res_Z({\cal E}_1/{\cal H},{\cal E}_2/{\cal
H})=-4\,ek(e^2 k + l^2)^2 X^6=0.
$$
Hence we obtain either $ek=0$ or $e^2 k + l^2=0$. Since the second
condition leads to   degenerate  systems, we must examine the
conditions $e=0$  and  $k=0$. For   systems (\ref{s4.7}) we have
$N_1=8ex^4$ and $N_2=16kx$ and  we shall consider two subcases:
$N_1=0$ and $N_1\ne0$, $N_2=0$.

{\bf 2a)} Assume $N_1=0$. Then $e=0$  and
 the systems (\ref{s4.7})  become
$ \  \dot x=k+ x^2,\  \dot y= l.\ $ Calculations yield $\
 {\cal E}_1= X\,{\cal H}, \ $ $ \
  {\cal E}_2= (X^2 + k Z^2)\,{\cal H},\ $ $\
 {\cal H}= \,l Z^2(X^2 + k Z^2),\ $ \mbox{$
\Res_X({\cal E}_1/{\cal H},\ {\cal E}_2/{\cal H})=4\,kZ^2.\ $}
Therefore  $\deg\,{\cal H}=4$ and the polynomials  ${\cal E}_1$
and ${\cal E}_2$ could not have an additional common factor if and
only if $k\ne0$ (i.e. $N_2\ne0$). Since $l\ne0$ (otherwise we get
degenerate systems) after the rescaling  $x\to |k|^{1/2} x$, $
y\to l|k|^{-1/2}y$ and $ t\to|k|^{-1/2}t$ we get the systems
\be\label{2CF_14}
\dot x=k+ x^2,\qquad
\dot y= 1
\ee
with $k\in\{-1,1\}$.  The systems (\ref{2CF_14}) possess two
invariant affine lines: $x=\pm\sqrt{-k}$ which are distinct due to
the condition $k\ne0$. Moreover, by Corollary \ref{Mult:Z=0} the
line $Z=0$ could be of multiplicity three.  This is confirmed by
the perturbations ({\it V.20${}_\varepsilon$}) (for $k<0$) and
({\it V.24${}_\varepsilon$}) (for $k>0$) from Table~5. On the
other hand for systems (\ref{2CF_14})   we have $N_5=-64kx^2$.
Therefore, for $N_1=0$ and $N_2\ne0$ we obtain the Config. 5.20 if
$N_5>0$ and the Config. 5.24 if $N_5<0$.

{\bf 2b)} Let $N_1\ne0$, $N_2=0$. In this case we have $e\ne0$,
$k=0$ and   systems (\ref{s4.7})  become
\be\label{S_III_a}
\dot x= x^2,\qquad
\dot y= l+ex.
\ee
For the systems (\ref{S_III_a}) calculations yield: $\ {\cal
E}_1=(e X + 2 l Z)\,{\cal H},\ $ $\
  {\cal E}_2=  X (e X + l Z)\,{\cal H},\ $ $
 {\cal H}= X^3Z.$\ Therefore  $\deg\,{\cal H}=4$ and since $l\ne0$ (otherwise we
get the degenerate systems (\ref{S_III_a})) and  $e\ne0$ (
$N_1\ne0$) we conclude that the polynomials  ${\cal E}_1$ and
${\cal E}_2$ could not have an additional common factor, i.e. each
non-degenerate system of the  family (\ref{S_III_a}) belongs to
$\QSL_{\bf5}$. Via the rescaling $x\to le^{-1}x$, $ y\to e\,y $
and $ t\to el^{-1}t$ systems (\ref{S_III_a}) become
\be\label{2CF_28}
\dot x= x^2,\qquad \dot y= 1+x.
\ee
This system  possesses the  invariant line $x=0$. Taking into
account the  polynomial ${\cal H}$, by Lemma \ref{lm3} and
Corollary \ref{Mult:Z=0} we obtain that  the line $x=0$ could be
of multiplicity three whereas  the line $Z=0$ could be of
multiplicity two. This is confirmed by the  perturbations ({\it
V.28${}_\varepsilon$}) from Table~5. Thus we obtain the Config.
5.28.
%%%%%%%%  New Subsubsection %%%%%%%
%%%%%
\vspace{-3mm}
\subsubsection{The case $N\ne0,$\ $\theta=0=B_3$}
\vspace{-3mm}
 Since for  the  systems~$(\SSS_{I\!I\!I})$
we have
\be\label{S_III_1}
 \theta=  -8h^2(g-1),\quad  \mu=32gh^2,\quad
N= (g^2-1)x^2+2h(g-1)xy+h^2y^2
\ee
we shall consider two cases: $\mu\ne0$ and $\mu=0$.
%%%%%%%%%%%%%%%%%%%%%%%%%

\smallskip\noindent
{\bf The subcase  $\mu\ne0$}. Then $gh\ne0$ and the condition
$\theta=0$ yields $g=1$.  Then the systems~$(\SSS_{I\!I\!I})$ with
$g=1$ by the transformation $\  x\to\ x -d/h,$ \ $ y\to (hy+
2d-ch)/h^2 $\ will be brought to the systems: $ \dot x=k+ x^2
+xy,\  \dot y= l+ex+fy+y^2, $ for which $ B_3=
-3e^2x^4+(3l-12k)x^2y^2-6kxy^3. $\ Hence the condition $B_3=0$
yields $e=k=l=0$ and we obtain the system
\be\label{2CF_24}
\dot x= x^2 +xy,\quad \dot y= fy+y^2,
\ee
for which we can assume  $f\in\{0,1\}$ via Remark \ref{rem:transf}
($\gamma=f,\ s=1$).
  For the systems~(\ref{2CF_24})
 calculations yield:\
$
 {\cal H}=\gcd\left({\cal E}_1,{\cal E}_2\right)= X^2Y(Y+fZ).
$\ Hence, $\deg\,{\cal H}=4$, i.e $M_{{}_\IL}\ge5$ and since
$N\ne0$ by  Lemma \ref{gcd:4} we have $M_{{}_\IL}<6$.
 By Theorem \ref{theor:E1,E2} the  systems~(\ref{2CF_24}) possess the invariant
lines $x=0,$ $y=0$ and $y+f=0$. Moreover, according to Lemma
\ref{lm3} the line $x=0$ of  could be of the multiplicity two, and
the lines $y=0$ and $y=-f$ are distinct if and only if $f\ne0$.
Since for the systems~(\ref{2CF_24}) we have $D= -f^2x^2y$, the
condition $f\ne0$ can be expressed by using this $T$-comitant.

If $D\ne0$ (then $f=1$) the perturbed systems ({\it
V.11${}_\varepsilon$}) from Table~5  show that the invariant line
$x=0$ is double one.
 Thus, for $D\ne0$ we obtain Config. 5.11.

Assume  $D=0$. Then $f=0$ and the   invariant line $x=0$ as well as the line
 $y=0$   is of  multiplicity two. This  is
confirmed by the  perturbed systems ({\it V.19${}_\varepsilon$})
 from Table~5. Therefore the case $D=0$ leads to the  Config. 5.19.
%%%%%%%%%%%%%%%%%%%%%

\smallskip\noindent
{\bf The subcase $\mu=0$.} In this case from (\ref{S_III_1}) we
obtain $h=0$ and the condition $N\ne0$ yields $g^2-1\ne0$. Then
the systems~$(\SSS_{I\!I\!I})$ with $h=0$ will be brought via the
translation\ $ x\to x+f/(1-g),$ $ y\to y+e/(1-g)$\  to the
systems:
\be\label{S_III_2}
\dot x=k+ cx +dy +gx^2,\quad
\dot y= l+(g-1)xy,
\ee
for which $\ B_3=
-3\,l\,(g-1)^2x^4-3\,c\,d(g-1)x^2y^2-6\,d^2g\,xy^3. $\ Hence, as
$N\ne0$ the condition $B_3=0$ yields $l=c\,d=d\,g=0$. We  claim
that if $d\ne0$ then for the systems (\ref{S_III_2}) we have
$M_{{}_\IL}<5$.
 Indeed, suppose $d\ne0$. Hence
the condition $B_3=0$ yields $l=c=g=0$. Thus we obtain the systems
$\ \dot x=k+ dy,\ \dot y= -xy,$\ for which calculations yield:
$$
{\cal E}_1=(-k X^2 + d^2 Y^2 + 2 d k Y Z + k^2 Z^2)\,{\cal
H},\quad
  {\cal E}_2= -X (d Y + k Z)^2{\cal H},\quad {\cal H}= YZ^2.
$$
Thus, $\deg\,{\cal H}=3$ and since $d\ne0$ to have an additional
common factor of ${\cal E}_1$ and ${\cal E}_2$   by Lemma
\ref{Trudi:2} the   following condition must hold: $\ \Res_Y({\cal
E}_1/{\cal H},\ {\cal E}_2/{\cal H})= d^4 k^2 X^6=0.$\ Therefore,
$dk=0$ and since  $d\ne0$ we obtain $k=0$. However this condition
leads to the degenerate systems. Our claim is proved.

Let us assume  $d=0$. Then the condition $B_3=0$ yields $l=0$ and
the systems~(\ref{S_III_2}) become
\be\label{S_III_3}
\dot x=k+ cx +gx^2,\quad
\dot y= (g-1)xy.
\ee
Calculations yield: $\
 {\cal E}_1=(X^2 - k Z^2)\,{\cal H},\ $ $
  {\cal E}_2= X(g X^2 + c X Z + k Z^2)\,{\cal H},\ $
  $  \Res_X({\cal
E}_1/{\cal H},\ {\cal E}_2/{\cal H})=k^2[c^2 - k(1 + g)^2] Z^6 $,\
where
 $ {\cal H}= (g-1)Y(g X^2 + c X Z + k Z^2). $\
Hence $\deg {\cal H}=3$ and we need  an additional common factor
of ${\cal E}_1$ and ${\cal E}_2$. For this, according to Lemma
\ref{Trudi:2}, the condition  $\Res_X({\cal E}_1/{\cal H},\ {\cal
E}_2/{\cal H})=0$ is necessary, i.e.
$k\left[(c^2-k(g+1)^2\right]=0$. As $k\ne0$ (otherwise we get
degenerate systems)  we obtain the condition
\mbox{$c^2-k(g+1)^2=0$}.

Assume  $c^2=k(g+1)^2$. Since $N\ne0$ (i.e. $g+1\ne0$) we may set
$c=u(g+1)$, where $u$ is a new parameter. Then $k=u^2\ne0$ and via
the Remark \ref{rem:transf} ($\gamma=u$, $s=1$) the
systems~(\ref{S_III_3}) will be brought to the form:
\be\label{2CF_26}
\dot x=1+ (g+1)x +gx^2,\quad
\dot y= (g-1)xy.
\ee
For systems~(\ref{2CF_26}) we obtain $\
 {\cal H}=\gcd\left({\cal E}_1,{\cal E}_2\right)= 2(g-1) Y(X + Z)^2(gX + Z).
$\ Hence $\deg\,{\cal H}=4$, i.e. $M_{{}_\IL}\ge5$  and since
$N\ne0$ by Lemma \ref{gcd:4}, $M_{{}_\IL}\ne6$ for any  system
(\ref{2CF_26}).

On the other hand the conditions $d=0$ and $c^2-k(g-1)^2=0$ are
equivalent to $B_3=H_6=0$. Indeed, the condition $B_3=0$ implies
$dg=0$ for system \eqref{S_III_2} (see above) and then $\quad
H_6=64(g-1)^2x^4\big[2(g-1)^2[k(g+1)^2-c^2]x^2-5cdxy-2d^2y^2\big].
$\ Hence as $N\ne0$ (i.e. $g-1\ne0$) the condition $H_6=0$ yields
$d=c^2-k(g-1)^2=0$.

We observe, that  $Z\mid {\cal H}$ if and only if  $g=0$. So,
since by $N\ne0$ the condition $g=0$ is equivalent to
$K=2g(g-1)x^2=0$, we shall examine two subcases: $K\ne0$ and
$K=0$.

%%%%%%%%%%%%%%%%%%%%%%%%%%

{\bf 1)} If $K\ne0$ then $g\ne0$ and according to Theorem
\ref{theor:E1,E2} the systems (\ref{2CF_26}) possess the following
invariant affine straight lines:\ $
   y=0,\  x+1=0,\  gx+1=0.
$\ By $g-1\ne0$ the invariant line $gx+1=0$ cannot coincide with
$x=-1$. Moreover, the line $x=-1$ could be  of multiplicity two
and this is confirmed by the  perturbations ({\it
V.14${}_\varepsilon$}) from Table~5. Thus for $K\ne0$ we obtain
the Config. 5.14.
%%%%%%%%%%%%%%%%%%%%%%%%%%

{\bf 2)} For $K=0$ we obtain $g=0$ and then
 the line $l_\infty:\ Z=0$ appears as a component of
a conic in the pencil of conics corresponding to
systems~(\ref{2CF_26}). Hence, the invariant line $x+1=0$ as well
as the line $Z=0$ is of multiplicity two, as is shown by the
perturbations ({\it V.18${}_\varepsilon$}) from Table~5. Thus for
$K=0$ we obtain  the Config. 5.18.
%%%%%%%%%%%%%%%%%%%%%%%%%%%%%%%%
%%%%%% New Subsection   %%%%%
%%%%%%%%%%%%%%%%%%%%%%%%%%%%%%%%
\vspace{-3mm}
\subsection{Systems with the divisor\ $D_S(C ,Z)=3\cdot w$}
\vspace{-3mm}
 We are in the case of the canonical form
$(\SSS_{I\!V})$ and we shall later need the following polynomial
 which is  shown to be a $CT$-comitant in  Lemma \ref{Table:Propreties}.
\bnot\label{not_N7} Let us denote\quad $ N_6(\ab,x,y)=
8D+C_2\left[8(C_0,D_2)^{(1)}-3(C_1,C_1)^{(2)}+2D_1^2\right]. $
\vspace{-2.6mm}
\enot
%%%%%%%%  New Subsubsection %%%%%%%%%%%%
\vspace{-3mm}
\subsubsection{The case $N=0=B_2$}
\vspace{-3mm} It was previously shown  (see page~\pageref{S4_N0}) that
for $N(\ab,x,y)=0$  we have to   examine the systems (\ref{S4_N0})
for which we have $\
  B_2= -648\,d^4 x^4,$\ $ B_3= 6dx^3(fx-dy).
$\ Thus the condition $B_2=0$ is equivalent to $B_3=0$ and this
yields $d=0$. Then we obtain the systems (\ref{s4.8}) for which
the expressions  of the affine comitants ${\cal E}_1$ and ${\cal
E}_2$ are given in  (\ref{val_Eib}). We observe that $\deg\,{\cal
H}=3$ and we need to have an additional common factor of ${\cal
E}_1$ and ${\cal E}_2$.  Since the polynomial ${\cal E}_2/{\cal
H}$ does not depend of $Y$, to have such a common factor by Lemma
\ref{Trudi:2} at least one of two following conditions must hold:
\beq
&& \Res_X({\cal E}_1/{\cal H},\ {\cal E}_2/{\cal H})= (c - f)^2 (c^2fY - k^2Z + c^2l Z)^2Z^4=0;\\
&& \Res_Z({\cal E}_1/{\cal H},\ {\cal E}_2/{\cal H})=(c - f)^2(c +
f) X^4(k^2 X - c^2l X + cf k Y)^2=0.
\eeq
Hence, we obtain either $(c-f)(c+f)=0$ or $k=cl=cf=0$, however the
second case leads to the degenerate systems (\ref{s4.8}). On the
other hand for  systems (\ref{s4.8})   we obtain $\
 N_3= 3(c-f)x^3,$\ $ D_1=c+f.$\
 So, the condition $(c-f)(c+f)=0$ is equivalent to $N_3D_1=0$ and we
shall consider two subcases: $N_3=0$ and $N_3\ne0$, $D_1=0$.
%%%%%%%%%%%%%%%%%%%%%%

\smallskip\noindent
{\bf The subcase  $N_3=0$}. Then
 $f=c$ and we get the systems:
\be\label{s4.8:f=c}
\dot x=k+ cx,\qquad
\dot y= l+cy-x^2
\ee
for which calculations yield: $\
 {\cal E}_1=2\,X\,{\cal H},$\quad $
  {\cal E}_2=Z (c X + k Z)\,{\cal H}, $\quad $ {\cal H}=
   Z^2(c X + k Z)^2.
$\ So $\deg\,{\cal H}=4$. It is easy to observe  that the
polynomials ${\cal E}_1$ and ${\cal E}_2$ do not have an
additional common factor  if and only if  $k\ne0$   and that the
polynomial ${\cal H}$ has a factor of multiplicity four if $c=0$.
On the other hand for systems (\ref{s4.8:f=c}) $D_1= 2c$ and
$N_4(\ab,x,y)=12kx^2$. Hence the conditions $c=0$ and $k\ne0$ can
be expressed by using the $CT$-comitants $D_1$ and $N_4$. So, we
 shall consider two cases: $D_1\ne0$ and $D_1=0$.
%%%%%%%%%%%

{\bf 1)} Assume $D_1\ne0$. Then $c\ne0$ and the
systems~(\ref{s4.8:f=c}) can be brought by the affine
transformation $\
  x=c^{-1}k x_1,\quad y=c^{-3}k^2 y_1-c^{-1}l,\quad t=c^{-1}t_1
$\ to the  system
\be\label{2CF_16}
\dot x=1+ x,\qquad \dot y= y-x^2
\ee
for which ${\cal H}=\gcd\left({\cal E}_1,{\cal
E}_2\right)=Z^2(X+Z)^2$. By Lemma \ref{lm3} and Corollary
\ref{Mult:Z=0} the line $x=-1$ could be of multiplicity 2, whereas
the line $Z=0$ could be of multiplicity 3. This is confirmed by
the   perturbations ({\it V.27${}_\varepsilon$}) from Table~5.
Thus, for $N_3=0,\ N_4\ne0$ and $D_1\ne0$ we obtain the Config.
5.27.
%%%%%%%%%%%

{\bf 2)} If $D_1=0$ for  systems (\ref{s4.8:f=c})
 we have $c=0$ and since $k\ne0$ we may consider $k=1$ via the rescaling
 $x\to kx$ and $y\to k^2y$. Thus we obtain the systems
\be\label{2CF_17}
\dot x=1,\qquad
\dot y= l-x^2,
\ee
and we can assume $l\in \{-1,0,1\}$ via the rescaling
$x\to|l|^{1/2}x$, $y\to|l|^{3/2}y$ and $t\to|l|^{1/2}t$. For the
systems (\ref{2CF_17}) we have ${\cal H}=\gcd\left({\cal
E}_1,{\cal E}_2\right)=Z^4$ and hence by Corollary \ref{Mult:Z=0}
the line $Z=0$ could be of multiplicity 5.  This is confirmed by
the  perturbation ({\it V.30${}_\varepsilon$}) from Table~5. Thus,
we obtain the Config. 5.30.
%%%%%%%%%%%%%%%%%

\smallskip\noindent
{\bf The subcase $N_3\ne0$, $D_1=0$}. These conditions yield
$f=-c\ne0$  and we may assume $c=-1$ via the Remark
\ref{rem:transf} ($\gamma=-c,$ $s=1$). Then  systems (\ref{s4.8})
will be brought by a translation to the systems:
\be\label{2CF_30}
\dot x=k-x,\qquad
\dot y=  y-x^2,
\ee
in which we can assume $k\in\{0,1\}$ due to the rescaling $x\to
kx$ and $y\to k^2y$. For  systems~(\ref{2CF_30}) we calculate: $\
{\cal E}_1=2( k X -  Y)\,{\cal H},\ $ $
  {\cal E}_2=(X - k Z)^2{\cal H},\ $ $
 {\cal H}= Z^3(k Z-X),\ $\ i.e. $\deg\,{\cal H}=4$. We observe that the polynomials
${\cal E}_1/{\cal H}$ and ${\cal E}_1/{\cal H}$ cannot have a
common factor. By Corollary \ref{Mult:Z=0} the line $Z=0$ could be
of multiplicity four. This  is confirmed by the perturbations
({\it V.26${}_\varepsilon$})from Table~5. Thus, for $D_1=0$ and
$N_3\ne0$ we obtain the Config. 5.26.
%%%%%% New Subsubsection   %%%%%
\vspace{-3mm}
\subsubsection{The case $N\ne0$, $\theta=0=B_3$}
\vspace{-3mm}
For the systems~$(\SSS_{I\!V})$ we  calculate: $\
 \theta= 8h^3,$\ $
N= (g^2-2h)x^2+2ghxy +h^2y^2. $\ From $\theta=0,$  $N\ne0$ we
obtain $h=0$, $g\ne0$ and we may assume $g=1$ via the rescaling
$x\to g^{-1}x$, $y\to g^{-2}y$. Then the systems~$(\SSS_{I\!V})$
with $h=0$ and $g=1$ will be brought by the translation $x\to
x-c/2$, \mbox{$y\to y-(c+e)$}
 to the systems:
$\quad \dot x=k+ dy+x^2,\ \ \dot y =l + fy -x^2 +xy. $\quad For
these systems we have $\ B_3=3(2df-k-f^2)x^4-
6d(d-4f)x^3y-9d^2x^2y^2. $\ Hence the condition $B_3=0$ yields
$d=0,$ $k=-f^2$ and  we obtain the systems
\be\label{S_IV_1}
\dot x=-f^2 + x^2,\quad \dot y =l + fy -x^2 +xy,
\ee
for which   calculations yield: $\ {\cal E}_1=(X^2 + 2f X Z + l
Z^2)\,{\cal H},\quad
  {\cal E}_2=3(X - f Z)(X + f Z)^2{\cal H},
$\ where  ${\cal H}= 2(X - fZ)^2(X + fZ)$.\ Hence $\deg\,{\cal
H}=3$ and to  have an additional common factor of ${\cal E}_1$ and
${\cal E}_2$,   according to Lemma~\ref{Trudi:2}, the condition \
$ \Res_X({\cal E}_1/{\cal H},\ {\cal E}_2/{\cal H})= 9(f^2 -
l)^2(3f^2 + l)Z^6=0$ \ must hold. We observe that for $l=f^2$ the
systems (\ref{S_IV_1}) become degenerate. Therefore $l=-3f^2$ and
since $f\ne0$ (otherwise we get the degenerate system) we may
assume $f=1$ via Remark \ref{rem:transf} ($\gamma=f$, $s=1$). Thus
we obtain the canonical system:
\be\label{2CF_29}
\dot x=-1 + x^2,\quad
\dot y =-3 + y -x^2 +xy
\ee
and calculations yield\ ${\cal H}=\gcd\left({\cal E}_1,{\cal
E}_2\right)= 2(X - Z)^3(X + Z)$. Therefore, according to Theorem
\ref{theor:E1,E2}, the system  (\ref{2CF_29}) possesses the
invariant straight lines $x=1$ and $x=-1$ and the line $x=1$ could
be of multiplicity three. This is confirmed by the
 perturbations ({\it V.23${}_\varepsilon$}) from Table~5.
 Thus we obtain  Config. 5.23. Since for systems (\ref{S_IV_1})
we have $N_6= 8(l+3f^2)x^3$, the condition $l+3f^2=0$ is
equivalent to $N_6=0$.

\smallskip
All the cases in  Theorem \ref{th_mil_5} are thus examined. To
finish the proof of the Theorem \ref{th_mil_5} it remains to show
that the conditions  occurring in the middle column of Table 4 are
affinely invariant. This follows from the proof of Lemma
\ref{Table:Propreties}.
  \EProof

\blm\label{Table:Propreties}
   The  polynomials which are used
   in  Theorems \ref{th_mil_6} or  \ref{th_mil_5}  have the properties
indicated in the Table 6. In the last column are indicated the
algebraic sets on which the  $GL$-comitants on the left are
$CT$-comitants.
\elm
\begin{table}[!htb]
\begin{center}
\begin{tabular}{|c|c|c|c|c|c|}
\multicolumn{6}{r}{\bf Table 6}\\[2mm]
\hline \raisebox{-0.7em}[0pt][0pt]{Case} &
\raisebox{-0.7em}[0pt][0pt]{$GL$-comitants}
 & \multicolumn{2}{c|}{Degree in }  & \raisebox{-0.7em}[0pt][0pt]{Weight} &  Algebraic subset \\
\cline{3-4}
  &  & $\ \ a\ \ $ & $\!x$ and $y\!$ &    &  $V(*)$   \\
\hline\hline \rule{0pt}{2ex}
 $1$ & $\eta(a)$,\ $\mu(a)$,\ $\theta(a)$ & $4$ &  $0$ & $ 2$  &  $V(0)$ \\[0.5mm]
\hline\rule{0pt}{2ex} $2$  & $C_2(a,x,y)$   & $1$ &  $3$ & $-1$  & $V(0)$\\[0.5mm]
\hline\rule{0pt}{2ex} $3$ & $H(a,x,y),\ K(a,x,y)$  & $2$ &  $2$ & $ 0$   &  $V(0)$ \\[0.5mm]
\hline\rule{0pt}{2ex} $4$ & $M(a,x,y),\ N(a,x,y)$   & $2$ &  $2$ & $ 0$   &  $V(0)$ \\[0.5mm]
\hline\rule{0pt}{2ex} $5$ & $D(a,x,y)$  & $3$ &  $3$ & $-1$   &  $V(0)$ \\[0.5mm]
\hline\rule{0pt}{2ex} $6$ & $B_1(a)$  & $12$ &  $0$ & $3$   &  $V(0)$ \\[0.5mm]
\hline\rule{0pt}{2ex} $7$ & $B_2(a,x,y)$  & $8$ &  $4$ & $0$   &  $V(0)$ \\[0.5mm]
\hline\rule{0pt}{2ex} $8$ & $B_3(a,x,y)$  & $4$ &  $4$ & $-1$   &  $V(0)$ \\[0.5mm]
\hline\rule{0pt}{2ex} $9$ & $H_1(a)$  & $6$ &  $0$ & $2$   &  $V(0)$ \\[0.5mm]
\hline\rule{0pt}{2ex} $10$ & $H_2(a,x,y))$  & $3$ &  $2$ & $0$   &  $V(0)$ \\[0.5mm]
\hline\rule{0pt}{2ex} $11$ & $H_3(a,x,y)$  & $4$ &  $2$ & $0$   &  $V(0)$ \\[0.5mm]
\hline\rule{0pt}{2ex} $12$ & $H_4(a)$  & $6$ &  $0$ & $2$   &  $V(0)$ \\[0.5mm]
\hline\rule{0pt}{2ex} $13$ & $H_5(a)$  & $8$ &  $0$ & $2$   &  $V(0)$ \\[0.5mm]
\hline\rule{0pt}{2ex} $14$ & $H_6(a,x,y))$  & $8$ &  $6$ & $0$   &  $V(0)$ \\[0.5mm]
\hline \rule{0pt}{2ex} $15$ & $N_1(a,x,y)$  & $3$ &  $4$ & $-1$   & $V(\eta,H)$ \\[0.5mm]
\hline\rule{0pt}{2ex}  $16$ & $N_2(a,x,y)$  & $3$ &  $1$ & $0$   &  $V(\eta,H,B_3)$ \\[0.5mm]
\hline\rule{0pt}{2ex}  $17$ & $N_3(a,x,y)$  & $2$ &  $3$ & $-1$   &  $V(M,N)$ \\[0.5mm]
\hline\rule{0pt}{2ex}  $18$ & $N_4(a,x,y)$  & $2$ &  $2$ & $-1$   &  $V(M,N,N_3)$ \\[0.5mm]
\hline\rule{0pt}{2ex}  $19$ & $N_5(a,x,y)$  & $4$ &  $2$ & $0$   &  $V(\eta,H,B_3)$ \\[0.5mm]
\hline\rule{0pt}{2ex}  $20$ & $N_6(a,x,y)$  & $3$ &  $3$ & $-1$   &  $V(M,\theta,B_3)$ \\[0.5mm]
\hline\rule{0pt}{2ex} $21$ & $D_1(a)$  & $1$ &  $0$ & $0$   &  $V(M,N)$ \\[0.5mm]
\hline
\end{tabular}
\end{center}
\end{table}
\BProof  {\it I. Cases  1,\ldots,14}.
Let us consider the action of the translation group $T(2,\R)$ on
systems in $\widehat\QS$. It $\tau\in T(2,\R)$, i.e. $\tau:\  x =
\tilde x+x_0,$ $y=\tilde y+y_0$ and $S_{\ab}$ is a system in
$\widehat\QS$ of coefficients $\ab\in \R^{12}$, then applying this
action to $S_{\ab}$ we obtain the system $S_{\tilde \ab}$ of
coefficients $\tilde\ab\in \R^{12}$, i.e.
$$
S_{\tilde\ab}:\quad \left\{\ba{l}\dot {\tilde x}=
P(\ab,x_0,y_0)+P_x(\ab,x_0,y_0)\tilde x
           +P_y(\ab,x_0,y_0)\tilde y+p_2(\ab,\tilde x,\tilde y), \\
\dot {\tilde y}= Q(\ab,x_0,y_0)+Q_x(\ab,x_0,y_0)\tilde
x+Q_y(\ab,x_0,y_0)\tilde y+q_2(\ab,\tilde x,\tilde y).\ea\right.
$$
Then calculations yield:
$$
\bal
  & U (\tilde\ab)=U (\ab)\quad
  \text{for each}\quad U\in \{\eta,\mu,\theta,B_1,H_1,H_4,H_5\}, \\
 & W (\tilde\ab,\tilde x,\tilde y)=W (\ab,\tilde x,\tilde y)\quad
  \text{for each}\quad W\in \{C_2,K,H,M,N,D,B_2,B_3,H_2,H_3,H_6\}. \\
\eal
$$
Since this holds for every $\ab\in\R^{12}$, according to
Definition \ref{def:T-com} we conclude that the GL- comitants
indicated in the lines 1--15 of Table 6 are $T$-comitants for
systems
\eqref{2l1}.

{\it II. Cases 15,\ldots,21}. {\bf 1)} We consider firstly the
$GL$-comitants $N_1(a,x,y)$, $N_2(a,x,y)$ and $N_5(a,x,y)$ which
according to Tables 2 and 4 were used only when the conditions
$\eta=0=H$ are satisfied. According to Lemma \ref{lm_3:2} for
$\eta=0$ there correspond three canonical forms: $(\SSS_{I\!I})$,
$(\SSS_{I\!I\!I})$ and $(\SSS_{V})$. Since for the systems
$(\SSS_{V})$ we have $H=-x^2\ne0$, we need to consider the
following cases: $(i)$ $\eta=0$ and $M\ne0$;\quad $(ii)$ $\eta=0$
and $M=0$ and $C_2\ne0$.

$(i)$ For $\eta=0$ and $M\ne0$ we are in the  class of systems
$(\SSS_{I\!I\!I})$, for which the condition
$H=-(g-1)^2x^2-2h(g+1)xy-hy^2=0$ yields $h=g-1=0$ and this leads
to  systems \eqref{S3_NH_0} (see page \pageref{S3_NH_0}). On the
other hand  for any system corresponding to a point $\tilde\ab\in
\R^{12}$ in the orbit under the translation group action of a
system \eqref{S3_NH_0}    calculations yield:
$$
\bal
& N_1(\tilde \ab,\tilde x,\tilde y)=8\tilde x^2(e\tilde
x^2-2d\tilde y^2),\quad N_2(\tilde \ab,\tilde x,\tilde
y)=4(f^2+4k)\tilde x + 4df  \tilde y+8d (x_0\tilde y +2y_0\tilde x),\\
& N_5(\tilde \ab,\tilde x,\tilde y)=-16(4k \tilde x^2 - d^2
 \tilde y^2) + 64d\tilde x (x_0\tilde y -y_0\tilde x),\quad B_3(\tilde \ab,\tilde x,\tilde y)
  = 6d\tilde x\tilde y^2(f\tilde x-d\tilde y).
\eal
$$
 Hence the value of $N_1$
does not depend of the vector defining the translation  and for
$B_3=0$ the same occurs for $N_2$ and $N_5$. Therefore we conclude
that for $M\ne0$ the polynomial $N_1$  is a $CT$-comitant modulo
$\langle\eta,H\rangle$, whereas the  polynomials $N_2$ and $N_5$
are $CT$-comitants modulo $\langle\eta,H,B_3\rangle$.

$(ii)$ Assume now that $M=0$ and $C_2\ne0$. Then we are in the
class of systems $(\SSS_{I\!V})$, for which the condition
$H=-(g^2+4h)x^2-2ghxy-hy^2=0$ yields $g=h=0$. In this case using
an additional translation (see page \pageref{S4_N0}) we obtain the
systems \eqref{S4_N0}. Then for any system corresponding to a
point $\tilde\ab\in \R^{12}$ in the orbit under the translation
group action of a system \eqref{S4_N0}    calculations yield: $\
N_1(\tilde \ab,\tilde x,\tilde y)=-24 d\tilde x^4,\quad N_2(\tilde
\ab,\tilde x,\tilde y)=12d(c+f)\tilde x,\quad N_5(\tilde
\ab,\tilde x,\tilde y)=0. $\ Since the condition $M=0$ implies
$\eta=0$, considering the case $(i)$ above we conclude that
independently of either $M\ne0$ or $M=0$, the $GL$-comitant $N_1$
is a $CT$-comitant modulo $\langle\eta,H\rangle$ and $N_2$ and
$N_5$ are $CT$-comitants modulo $\langle\eta,H,B_3\rangle$.

{\bf 2)} Let us now consider  the $GL$-comitants $N_3(a,x,y)$,
$N_4(a,x,y)$, $N_6(a,x,y)$ and $D_1(a)$. According to Tables 2 and
4 the polynomials $N_3$, $N_4$ and $D_1$ (respectively $N_6$) were
used only when the conditions $M=N=0$ (respectively
$M=\theta=0,N\ne0$)  are satisfied.  In both cases we are in the
class of systems $(\SSS_{I\!V})$ and we shall consider the two
subcases: $N=0$ and $N\ne0$, $\theta=0$.

$(i)$ If for  the system $(\SSS_{I\!V})$  the condition $N=0$ is
fulfilled then as it was shown on the page \pageref{S4_N0} we
obtain systems \eqref{S4_N0}. Then for any system corresponding to
a point $\tilde\ab\in \R^{12}$ in the orbit under the translation
group action of a system \eqref{S4_N0} calculations yield:
$$
\bal
& N_3(\tilde \ab,\tilde x,\tilde y)=3(c-f)\tilde x^3+2d\tilde x^2\tilde y,\quad
B_3(\tilde \ab,\tilde x,\tilde y)=6d\tilde x^3(f\tilde x-d\tilde y),\\
& N_4(\tilde \ab,\tilde x,\tilde y)=12 k\tilde x^2 +
3(f^2-c^2)\tilde x\tilde y
               -3d(c+f)\tilde y^2+6\tilde x^2[(c-f)x_0+2dy_0],\\
& N_6(\tilde \ab,\tilde x,\tilde y)=8c(c-f)\tilde x^3 + 16 df
\tilde x^2 \tilde y-8d^2\tilde x\tilde y^2-48d x_0\tilde x^3,\quad
D_1(\tilde \ab)=c+f.
\eal
$$
These relations show us that: $(\alpha)$ the $GL$-comitants $N_3$
and $D_1$ are  $CT$-comitants modulo $\langle M,N\rangle$;
$(\beta)$ the $GL$-comitant $N_4$ is a $CT$-comitant modulo
$\langle M,N,N_3\rangle$; $(\gamma)$ the $GL$-comitant $N_6$ is a
$CT$-comitant modulo $\langle M,N,B_3\rangle$.

$(ii)$ Assume that for the system $(\SSS_{I\!V})$  the conditions
$\theta=0$ and $N\ne0$ are fulfilled. As it was shown on the page
\pageref{S_IV_1} for $B_3=0$ we obtain systems
\eqref{S_IV_1}. For any system corresponding to a point
$\tilde\ab\in \R^{12}$ in the orbit under the translation group
action of a system \eqref{S_IV_1} we have $N_6(\tilde \ab,\tilde
x,\tilde y)=8(l+3f^2)\tilde x^3$. Therefore, since the condition
$N=0$ implies $\theta=0$, considering the case
 $(i)$ above we conclude that independently of either
$N\ne0$ or $N=0$, the $GL$-comitant $N_6$  is a $CT$-comitant
modulo $\langle M,\theta,B_3\rangle$.

The Table 6 shows us that all the conditions indicated in the
middle column of Tables 2 and 4 are affinely invariant. Indeed,
the $CT$-comitants $N_i$, $i=1,\ldots,...,6$ and $D_1$ are used in
Table 2 only for the varieties indicated in the last column of
Table 6. This complete the proof of the Theorems \ref{th_mil_6}
and \ref{th_mil_5}.

 \EProof

%%%%%%%%%%%%%%%%%%%%%%     REFERENCES     %%%%%%%%%%%%%%%%%%%%%%%%


{\footnotesize
\begin{thebibliography}{99}
\bibitem{Art_Llib1} { J. Artes, J.Llibre},
\emph{On the number of slopes of invariant straight lines for
polynomial differential systems}. {J. of Nanjing University
\textbf{ 13} (1996), 143--149. }
\bibitem{Art_Llib2} {J. Artes, B. Gr\"unbaum, J.Llibre},
\emph{On the number of invariant straight lines for polynomial
differential systems}. Pacific Journal of Mathematics \textbf{
184}, (1998), 207--230.

\bibitem{Bul_Tim} D. Boularas, Iu. Calin, L. Timochouk, N. Vulpe.
\emph{T-comitants of qudratic systems: A study via the
          translation invariants.}
 Delft University of Technology, Faculty of Technical
          Mathematics and Informatics, Report no. 96-90, 1996;
        (URL: {\small \tt
ftp://ftp.its.tudelft.nl/publications/tech-reports/1996/
DUT-TWI-96-90.ps.gz}

\bibitem{Cris_Llib} { C. Christopher, J. Llibre}.
 \emph{Integrability via   invariant algebraic
curves for planar polynomial differential systems.} Ann.
Differential equations   {\bf 16} (2000), no. 1, 5-19.

\bibitem{Cris_Llib_Per} { C. Christopher, J. Llibre, J. V. Pereira},
\emph{Multiplicity of invariant algebraic curves}. Preprint 2002.

\bibitem{Cris_Llib_Pant} { C. Christopher, J. Llibre, C. Pantazi,
X. Zhang}, \emph{Darboux integrability and invariant algebraic
curves for planar polynomial systems.} J. Phys. {\bf A 35} (2002),
no. 10, 2457-2476.


\bibitem{Darb} {G. Darboux}, \emph{M\'emoire sur les \'equations
    diff\'erentielles du premier ordre et du premier
    degr\'e}. Bulletin de Sciences Math\'ematiques, 2me s\'erie,
  {\bf 2} (1) (1878), 60-96; 123-144; 151-200.

\bibitem{Druzhkova} {T.A. Druzhkova},
\emph{Quadratic differential systems with algebraic integrals}.
Qualitative theory of differential equations, Gorky Universitet
\textbf{ 2} (1975), 34--42 (Russian).

\bibitem{Fult}
 W. Fulton, {\it Algebraic curves. An introduction to Algebraic
 Geometry}.
W.A. Benjamin, Inc., New York, 1969.

\bibitem{Gr_Yng} J. H. Grace, A. Young,
\emph{The algebra of invariants}. New York: Stechert, 1941.


\bibitem{Lib_DS}  J. Llibre, D. Schlomiuk,
\emph{ The geometry of quadratic systems with a weak focus of
third order}. To appear in the Canadian J. of Math. (A previous
version of this paper appeared as Preprint, n\'um. 486, Nov. 2001.
CRM, Barcelona, 48 pp.)

\bibitem{Lib_Vul}  J. Llibre, N. Vulpe,
\emph{ Planar cubic  polynomial differential systems with
     the maximum  number of invariant straight lines}.
 Report, n\'um. 34,  2002, Universitat Aut\`onoma de Barcelona,
54 pp.

\bibitem{Lyubim1} R.A. Lyubimova, \emph{ On some differential equation which  possesses invariant
lines}.  Differential and integral eequations, Gorky Universitet,
{\bf 1}, 1977 (Russian).

\bibitem{Lyubim2} R.A. Lyubimova, \emph{On some differential equation which possesses invariant
lines}.  Differential and integral eequations, Gorky Universitet,
{\bf  21 }, 1984 (Russian).


\bibitem{Olver} {P.J. Olver},
\emph{ Classical Invariant Theory}.  London Mathematical Society
student texts: \textbf{ 44}, Cambridge University Press, 1999.

\bibitem{Po1}
{H. Poincar\'e}, \emph{ M\'emoire sur les courbes d\'efinies par
les \'equations diff\'erentielles}. J. Math. Pures Appl. (4) {\bf
1} (1885), 167--244; O\!euvres de Henri Poincar\'e, Vol. {\bf 1},
Gauthier--Villard, Paris, 1951, pp 95--114.

\bibitem{Po2}
{ H. Poincar\'e}, \emph{ Sur l'int\'egration alg\'ebrique des
  \'equations diff\'erentielles}. C. R. Acad. Sci. Paris, {\bf 112} (1891),
761--764.

\bibitem{Po3}
{ H. Poincar\'e}, \emph{ Sur l'int\'egration alg\'ebrique des
  \'equations diff\'erentielles du premier ordre et du premier degr\'e
  } I. Rend. Circ.Mat. Palermo {\bf5} (1891), 169-191.

\bibitem{Po4}
{ H. Poincar\'e}, \emph{ Sur l'int\'egration alg\'ebrique des
  \'equations diff\'erentielles du premier ordre et du premier
  degr\'e}
II. Rend. Circ.Mat. Palermo {\bf11} (1897), 169-193-239.


\bibitem{Popa2} M.N. Popa,  \emph {Application of invariant
processes to the study of homogeneous linear particular integrals
of a differential system}. Dokl. Akad. Nauk SSSR, {\bf 317}, no.
4, 1991 (Russian); translation in Soviet Math. Dokl. {\bf 43}
(1991), no. 2.


\bibitem{Popa4} M.N. Popa,  \emph {Aplications of algebras to differential
systems}. Academy of Science of Moldova  2001 (Russian).


\bibitem{Popa_Sib1} M.N. Popa and K. S. Sibirskii, \emph{Conditions
for the existence of a homogeneous linear partial integral of a
differential system}. Differentsial'nye Uravneniya, {\bf 23}, no.
8, 1987 (Russian).


\bibitem {Dana1}
{ D. Schlomiuk}, \emph{ Elementary first integrals and algebraic
invariant curves of differential equations}. Expo. Math. {\bf11}
(1993),  433--454.

\bibitem {Dana2}
{ D. Schlomiuk}, \emph{ Algebraic and Geometric Aspects of the
Theory of Polynomial Vector Fields}. In Bifurcations and Periodic
Orbits of Vector Fields, D. Schlomiuk (ed.), 1993, 429--467.

\bibitem{Pal_DS2}  D. Schlomiuk, J. Pal,
\emph{ On the Geometry in the Neighborhood of Infinity of
Quadratic Differential Systems with a
 Weak Focus}.  QualitativeTheory of Dynamical
Systems, $\mathbf 2$~(2001), no. 1, 1-43

\bibitem{Dana_Vlp1}  D. Schlomiuk,  N. Vulpe,
\emph{Planar quadratic differential systems
         with invariant straight lines of at least five total multiplicity}, CRM  Report no. 2922,
 Universit\'e de Montr\'eal, 2003, 42 pp.

\bibitem{Sib1}  K. S. Sibirskii.
\emph{  Introduction to the algebraic theory of invariants of
differential equations.} Translated from the Russian. Nonlinear
Science: Theory and Applications. Manchester University Press,
Manchester, 1988.

\bibitem{Sib2}  K. S. Sibirskii.
\emph{ Method of invariants in the qualitative theory of
differential equations}. Kishinev: RIO AN Moldavian SSR, 1968.

\bibitem{Sib3}  K. S. Sibirskii,
\emph{Conditions for the presence of a straight integral line of a
quadratic differential system in the case of a center or a focus}.
Mat. Issled. No. 106, Differ. Uravneniya i Mat. Fizika, 1989
(Russian).

\bibitem{Sokulski} J. Sokulski,
\emph{ On the number of invariant lines for polynomial vector
fields}. Nonlinearity, {\bf  9 } 1996.

\bibitem{Vlp1} N.I.Vulpe.
\emph{ Polynomial bases of comitants of differential systems and
their applications in qualitative theory}. (Russian)
``Shtiintsa'', Kishinev, 1986.

\bibitem{Walker} R.J.Walker.
\emph{ Algebraic Curves}. Dover Publications, Inc., New York,
1962.


\bibitem{ZX} { Zhang Xiang},
\emph{ Number of integral lines of polynomial systems of degree
three and four}. J. of Nanjing University, Math. Biquarterly
\textbf{ 10} (1993), 209--212.

\end{thebibliography}
}
\end{document}